\documentclass[11pt,fullpage,doublespace]{article}
\usepackage{graphicx} 
\usepackage{setspace} 
\usepackage{latexsym} 
\usepackage{amsfonts} 
\usepackage{amsmath} 
\usepackage{amssymb}
\usepackage{accents}
\usepackage{textcomp}
\usepackage{enumerate}
\usepackage{mathrsfs}
\usepackage{bm}
\usepackage{stmaryrd}
\usepackage[margin=1in]{geometry}
\usepackage{amsthm}
\usepackage{amssymb,latexsym,amsmath}
\usepackage{graphics}
\usepackage{titlesec}
\usepackage{changepage}
\usepackage{lipsum}
\usepackage{thmtools}
\usepackage{comment}
\usepackage[colorlinks=true,linkcolor=blue]{hyperref}
\usepackage{centernot}
\usepackage{mathtools}
\usepackage{tikz}
\usepackage{scrextend}
\usetikzlibrary{matrix,arrows,positioning,scopes}
\declaretheorem[style=definition,qed=$\dashv$,numberwithin=section]
{definition}

\declaretheorem[style=plain,sibling=definition]{theorem}
\declaretheorem[style=plain,sibling=definition]{lemma}

\declaretheorem[style=plain,sibling=definition]{question}
\declaretheorem[style=plain,sibling=definition]{conjecture}
\declaretheorem[style=plain,sibling=definition]{fact}
\declaretheorem[style=plain,sibling=definition]{proposition}
\declaretheorem[style=plain,sibling=definition]{corollary}
\declaretheorem[style=definition,sibling=definition]{remark}
\declaretheorem[style=plain,sibling=definition]{claim}
\declaretheorem[style=plain,sibling=definition]{claim*}

\titleformat{\section}{\normalsize\centering}{\thesection.}{1em}{}
\titleformat{\subsection}{\normalsize\centering}{\thesubsection.}{1em}{}
\titleformat{\subsubsection}{\normalsize}{\thesubsubsection.}{1em}{}
\numberwithin{equation}{section}

\DeclareMathOperator{\bfEnv}{\mathbf{Env}}

\newcommand{\gTheta}{\mathsf{G}}

\renewcommand{\gg}{\gamma}

\newcommand{\rest}{\restriction}
\newcommand{\la}{\langle}
\newcommand{\ra}{\rangle}

\newcommand{\param}{\mathfrak{P}}
\newcommand{\dom}{{\rm dom}}

\newcommand{\cf}{{\rm cf}}
\newcommand{\lh}{{\rm lh}}

\newcommand{\crit}{{\rm crt}}

\def\k{\kappa}
\def\a{\alpha}

\def\d{\delta}

\def\l{\lambda}

\renewcommand{\models}{\vDash}
\newcommand{\powerset}{{\wp }}

\def\P{{\mathcal{P} }}
\def\W{{\mathcal{W} }}
\def\Q{{\mathcal{ Q}}}

\def\L{{\rm{L}}}

\def\R{{\mathcal R}}

\def\H{{\rm{HOD}}}
\def\M{{\mathcal{M}}}
\def\N{{\mathcal{N}}}
\def\F{{\mathcal{F}}}
\def\T {{\mathcal{T}}}
\def\U{{\mathcal{U}}}
\def\S{{\mathcal{S}}}
\def\G{{\mathcal{G}}}

\def\VT{{\vec{\mathcal{T}}}}

\newcommand{\rthm}[1]{Theorem~\ref{#1}}

\newcommand{\rfig}[1]{Figure~\ref{#1}}

\newcommand{\rsec}[1]{Section~\ref{#1}}

\newcommand{\ins}{\trianglelefteq}

\newcommand{\pins}{\triangleleft}

\newcommand{\OR}{\text{OR}}
\newcommand{\J}{\mathcal J}
\newcommand{\un}{\cup}
\newcommand{\sub}{\subseteq}
\newcommand{\om}{\omega}

\newcommand{\es}{\mathbb{E}}
\newcommand{\Tt}{\mathcal{T}}

\newcommand{\Uu}{\mathcal{U}}

\newcommand{\conc}{\ \widehat{\ }\ }

\newcommand{\Lp}{\mathrm{Lp}}
\newcommand{\sats}{\vDash}
\newcommand{\cut}{\backslash}

\newcommand{\pow}{\mathcal{P}}
\newcommand{\inter}{\cap}

\renewcommand{\l}{\mathit{l}}

\renewcommand{\OR}{\textit{o}}

\newcommand{\Ord}{\text{Ord}}

\newcommand{\HC}{\mathrm{HC}}

\newcommand{\BBB}{\mathfrak{B}}
\newcommand{\Ll}{\mathcal{L}}
\newcommand{\all}{\forall}
\newcommand{\ex}{\exists}

\newcommand{\id}{\mathrm{id}}
\newcommand{\cross}{\times}
\newcommand{\Mbar}{\bar{\M}}
\newcommand{\gammabar}{\bar{\gamma}}
\newcommand{\Ttbar}{\bar{\Tt}}
\newcommand{\bbar}{\bar{b}}
\newcommand{\Bbar}{\bar{B}}
\newcommand{\Rbar}{\bar{\R}}

\newcommand{\C}{\mathcal{C}}

\newcommand{\rank}{\mathrm{rank}}

\newcommand{\MFsharp}{\mathfrak{M}}
\newcommand{\Jop}{\J^{\mathrm{m}}}
\newcommand{\Pbar}{\bar{\P}}
\newcommand{\abar}{\bar{a}}

\newcommand{\Fop}{\mathcal{F}}

\newcommand{\g}{\mathrm{g}}

\newcommand{\trancl}{\mathrm{trancl}}
\newcommand{\ahat}{\hat{a}}

\newcommand{\wh}{\mathrm{wh}}

\newcommand{\ord}{\mathrm{ON}}
\newcommand{\PR}{\powerset_{\omega_1}(\mathbb{R})}
\newcommand{\PI}{\mathbb{P}_{\mathcal{I}}}

\newcommand*{\TitleFont}{%
      \usefont{\encodingdefault}{\rmdefault}{b}{n}%
      \fontsize{12}{16}%
      \selectfont}

 \input xy
 \xyoption{all}
\begin{document}
\title{\TitleFont IDEALS AND STRONG AXIOMS OF DETERMINACY}

\author{
Dominik Adolf\footnote{Department of Mathematics, University of North Texas, Denton, TX, USA. Email: dominikt.adolf@googlemail.com} , Grigor Sargsyan\footnote{Institute of Mathematics of the Polish Academy of Sciences, Poland. Email: gsargsyan@impan.pl}, Nam Trang\footnote{Department of Mathematics, University of North Texas, Denton, TX, USA. Email: Nam.Trang@unt.edu} , Trevor M. Wilson\footnote{Miami University, Oxford, Ohio, USA. Email: twilson@miamioh.edu}, Martin Zeman\footnote{UC Irvine, Irvine, CA, USA. Email: mzeman@math.uci.edu}} 

\date{}
\maketitle
\begin{abstract}

$\Theta$ is the least ordinal $\a$ with the property that there is no surjection $f:\mathbb{R}\rightarrow \a$. ${\sf{AD}}_{\mathbb{R}}$ is the ${\sf{Axiom\ of\ Determinacy}}$ for games played on the reals. It asserts that every game of length $\omega$ of perfect information in which players take turns to play reals is determined. For a sentence $\phi$ in the language of set theory, we say that $M$ is the minimal model of ${\sf{ZF}}+{\sf{AD}}_\mathbb{R}+\phi$ if $M$ is a transitive model of ${\sf{ZF}}+{\sf{AD}}_\mathbb{R}+\phi$ containing all reals and ordinals, and whenever $N$ is a transitive model of ${\sf{ZF}}+{\sf{AD}}_\mathbb{R}+\phi$ containing all reals and ordinals then $M\subseteq N$. We consider the theories, where $\sf{CH}$ stands for the Continuum Hypothesis,\\\\
(${\sf{T_1}}$) $\sf{ZFC} + \sf{CH}  +``$There is an $\omega_1$-dense ideal on $\omega_1$."  \\
(${\sf{T_2}}$) ${\sf{ZF}}+{\sf{AD}}_\mathbb{R} + {}$``$\Theta$ is a regular cardinal.''\\\\
The main result of this paper is that ${\sf{T_1}}$ implies that the minimal model of ${\sf{T_2}}$ exists.   
Woodin, in unpublished work, showed that the consistency of ${\sf{T_2}}$ implies the consistency of ${\sf{T_1}}$. We will also give a proof of this result, which, together with our main theorem, establish the equiconsistency of ${\sf{T_1}}$ and ${\sf{T_2}}$. 

As a consequence, this resolves part of question 12 in \cite{Woodin}; in particular, it shows that the theories (b) and (c) in \cite[Question 12]{Woodin} are equiconsistent. Thus, our work completes the work that started by Woodin and Ketchersid in \cite{ketchersid2000toward} some 25 years ago.  We also establish other theorems of similar nature in this paper, showing the equiconsistency of ${\sf{T_2}}$ and the statement that   the nonstationary ideal on $\powerset_{\omega_1}(\mathbb{R})$ is strong and pseudo-homogeneous. The aforementioned results are the only known equiconsistency results at the level of $\sf{AD}_\mathbb{R} +{}$``$\Theta$ is a regular cardinal.''



\end{abstract}

\section{INTRODUCTION}

This paper studies the consistency of strong determinacy theories, specifically the theory 
\begin{center} 
$``\sf{ZF+AD}_{\mathbb{R}} + $$\Theta$ is regular" 
\end{center}
and the consistency of strong ideals on $\omega_1$, specifically $\omega_1$-dense ideals on $\omega_1$. The main theorems of the paper, Theorem \ref{thm:DI} and Corollary \ref{cor:equi2}, resolve a long-standing conjecture by Woodin in \cite{Woodin}. The work in this paper contributes to our understanding of and helps further establish the close connections between ideals and determinacy, two very seemingly different areas in set theory. 
\\\\
\textbf{Some background} \\\\
Famously, Ulam's investigations of the ${\sf{Measure\ Problem}}$, which asks whether there is a measure on $[0, 1]$, led  him to prove that there is no countably complete 0-1 measure, that is an \textit{ultrafilter}, on $\omega_1$ (e.g. \cite[Chapter 10]{Jech}). Ulam's theorem is often presented as showing that $\omega_1$ is not a measurable cardinal, where we say that $\k$ is a measurable cardinal if there is a  $\kappa$-complete ultrafilter $U$ on $\k$. 

Ulam's theorem and the ${\sf{Measure\ Problem}}$ in general have been a source of great ideas in set theory, and one of these ideas has been the study of $\textit{ideals}$ that could induce nice ultrafilters on uncountable cardinals. Suppose, for example, that $\mathcal{I}\subseteq \powerset(\k)$ is an ideal on $\k$. 
Let $\mathbb{P}_{\mathcal{I}} = \powerset(\omega_1)\slash \mathcal{I}$ be the corresponding boolean algebra induced by $\mathcal{I}$. One can also think of $\mathbb{P}_{\mathcal{I}}$ as a poset ordered by inclusion. It is not hard to see that if $U$ is a generic object for $\mathbb{P}_{\mathcal{I}}$\footnote{I.e. intersects all dense open subsets of $\mathbb{P}_{\mathcal{I}}$.} then the function $U^*:( \powerset(\k))^V\rightarrow \{0, 1\}$ given by $U^*(A)=0 \iff A\not \in U$ satisfies many of the properties of being a 0-1 measure with two major deficiencies. First $U^*$ may not measure all subsets of $\k$ that exist in $V[U]$, and second, $U^*$ may not be countably complete. It is then unclear exactly in what way this approach could lead to a reasonable study of the ${\sf{Measure\ Problem}}$.

The concept of ultrapower introduced the necessary formalism to eliminate the aforementioned issues. It is a well-known fact that a cardinal $\kappa$ is a measurable cardinal if and only if there is an elementary embedding $j: V\rightarrow M$ such that $M$ is a transitive class of $V$, $j\not=id$, $j\rest \k=id$ and $j(\k)>\k$. If $\k$ is a measurable cardinal then one obtains the $M$ above as an ultrapower of $V$ by a $\k$-complete ultrafilter on $\k$. The same can be done with our generic $U$ above, and for the start one can only demand the well-foundedness of $Ult(V, U)$. 

An ideal $\mathcal{I}$ is called \textit{precipitous} if whenever $U\subseteq \mathbb{P}_{\mathcal{I}}$ is a generic ultrafilter, the generic ultrapower of $V$ by $U$, $Ult(V, U)$, is well-founded. This approach to the ${\sf{Measure\ Problem}}$ has been incredibly fruitful and has lead to many great discoveries. The story has been partially told in Foreman's long manuscript \cite{Foreman}. The study of precipitous ideals has led to solutions of problems considered not just by set theorists but by wider mathematical community. For example, Theorem 5.42 of \cite{Foreman} states that  the existence of a certain nice ideal implies among other things that every projective set of reals is Lebesgue measurable.

Let $\mathcal{I}$ be an ideal on $\omega_1$. We write $\mathcal{I}^+$ for the collection of $\mathcal{I}$-positive sets and $\mathcal{F}_\mathcal{I}$ for the dual filter of $\mathcal{I}$.  $\mathcal{I}$ is \textit{$\kappa$-saturated} if there is no family $(S_i : i<\kappa)$ of sets in $\powerset(\omega_1)\backslash \mathcal{I}$ such that $S_i \cap S_j \in \mathcal{I}$ for all $i\neq j$; in other words, $\mathcal{I}$ is $\kappa$-saturated if there is no antichain in $\mathbb{P}_{\mathcal{I}}$ of size $\kappa$. A $\kappa$-complete ideal $\mathcal{I}$ is \textit{saturated} if it is $\kappa^+$-saturated. \textit{Presaturation} is a technical weakening of saturation. $\mathcal{I}$ is \textit{presaturated} if for any $A\in \powerset(\omega_1)\backslash \mathcal{I}$, any sequence of maximal antichains $(\mathcal{A}_i : i<\omega)$ in $\powerset(\omega_1)\slash \mathcal{I}$, there is $B\subseteq A$ such that $B\notin \mathcal{I}$ and such that for each $i<\omega$, $\{X\in \mathcal{A}_i : X\cap B\notin \mathcal{I}\}$ has cardinality at most $\omega_1$.  The reader can consult Foreman's paper \cite{Foreman} and Woodin's book \cite{Woodin} for more comprehensive discussions on the topic of ideals, which is an important area of research in modern set theory.

Shelah, Jensen, and Steel have established the following famous theorem, which is one of the first equiconsistency results that connects ideals and large cardinals.

\begin{theorem}
The following theories are equiconsistent.
\begin{enumerate}
\item $\sf{ZFC}+{}$There is a pre-saturated ideal on $\omega_1$.
\item $\sf{ZFC}+{}$There is a saturated ideal on $\omega_1$.
\item $\sf{ZFC}+{}$There is a Woodin cardinal.
\end{enumerate}
\end{theorem}

Shelah proves the consistency of 1 and 2 from the consistency of 3 by forcing techniques. Jensen and Steel prove the converse using inner model theoretic techniques, in particular core model theory. See, for example, \cite{CMIP, jensen2013k}.

Claverie and Schindler \cite{claverie77woodin} have improved the above result and shown that in fact theory 3 above is equiconsistent with the theory ``there is a strong ideal on $\omega_1$".\footnote{The property of being strong is weaker than being pre-saturated. Strong ideals are precipitous.}

Below we say that the ideal $\mathcal{I}\subseteq \powerset(\omega_1)$ is \textit{$\omega_1$-dense} if $\mathbb{P}_{\mathcal{I}}$ has a dense set of size $\omega_1$. $\omega_1$-density is a stronger property than saturation. The consistency question of $\omega_1$-dense ideals has been studied extensively in the last 25 years, starting with Woodin \cite{Woodin} and by various other authors in \cite{ketchersid2000toward, CMI}. Unlike saturation and presaturation of ideals on $\omega_1$, which can be forced from a relatively mild large cardinal like a Woodin cardinal, there is no known traditional forcing construction of an $\omega_1$-dense ideal from a large cardinal significantly weaker than an almost huge cardinal. \textit{The Axiom of Determinacy $(\sf{AD})$} comes into the picture in a rather surprising and dramatic fashion via the following remarkable theorem of Woodin.  

Recall that $\sf{AD}$ states that every infinite-length, two-person game of perfect information where players take turns to play integers is determined, i.e. one of the players has a winning strategy.  It is worth noting that $\sf{AD}$ is equiconsistent with ``$\sf{ZFC} + $ there are $\omega$ many Woodin cardinals" and the latter theory is much weaker than an almost huge cardinal. A (stronger) variation of $\sf{AD}$ is $\sf{AD}_{\mathbb{R}}$, which is like $\sf{AD}$ except the players are allowed to play reals. The theory ``$\sf{ZF+AD}_{\mathbb{R}} + $$\Theta$ is regular" is strictly stronger than $\sf{ZF+AD}_{\mathbb{R}}$; here $\Theta$ is the supremum of ordinals $\alpha$ for which there is a surjection from $\mathbb{R}$ onto $\alpha$. If the axiom of choice holds, then $\Theta = \mathfrak{c}^+$, the successor of the continuum. If $\sf{AD}$ holds, then $\Theta$ is a limit of measurable cardinals and more. In the following theorem and this paper, $L(\mathbb{R})$ is the minimal model of ${\sf{ZF}}$ that contains all the ordinals and the real numbers (see \cite[Theorem 2.11.1]{CMI}).

\begin{theorem}[Woodin]\label{woodin} The following theories are equiconsistent.
\begin{enumerate}
\item ${\sf{ZFC}}\   + $ ``There is an $\omega_1$-dense ideal on $\omega_1$".
\item ${\sf{AD}}$ holds in $L(\mathbb{R})$. 
\end{enumerate}
\end{theorem}

Woodin introduces two very important sets of techniques in the proof of the above theorem. In one direction, to show the consistency of ${\sf{ZFC}}$+``There is an $\omega_1$-dense ideal on $\omega_1$", he develops powerful and general forcing techniques over models of determinacy, i.e. $\mathbb{P}_{\textrm{max}}$ and its cousins (cf. \cite{Woodin}).\footnote{This work was partly inspired by previous work by Steel and Van Wesep \cite{steel1982two}.} To prove the other direction, \textit{the core model induction ($\sf{CMI}$)} technique was introduced. $\sf{CMI}$ is a general technique for obtaining lower-bound consistency by inductively proving determinacy in canonical models like $L(\mathbb{R})$. $\sf{CMI}$ has since then been developed further into a very powerful and versatile method for proving lower-bound consistency and equiconsistency results (see for example \cite{CMI, Trang2015PFA,wilson2012contributions,steel2014determinacy}) from a variety of hypotheses.  Part of this development is in understanding determinacy models beyond $L(\mathbb{R})$\footnote{More precisely, $\sf{AD}^+$ models. $\sf{AD}^+$ is a technical strengthening of $\sf{AD}$.} and their canonical inner models of large cardinals (like the HOD).

The aforementioned theorems of Shelah, Jensen, Steel, Woodin, and others demonstrate intimate connections between different branches of set theory, namely the study of precipitous ideals and the study of models of ${\sf{AD}}$. It seems that the connections that theorems like above establish are rooted in the naturalness of the constructions that produce the models of one theory given a model of another, and this naturalness --- the feeling of having no barriers to naturally drift from one theory to another as if they were one and the same theory --- is not fully expressed in the formal statement of the theorem, namely that the two theories are equiconsistent. We discuss this a bit more later in this section.

The main theorem of this paper, \rthm{thm:DI}, has the same spirit as Woodin's theorem above.

\begin{definition}[$\sf{DI}$]\label{dfn:DI}
Let $\sf{DI}$ be the conjunction of 
\begin{itemize}
\item $\sf{CH}$,
\item there is an $\omega_1$-dense ideal $\mathcal{I}$ on $\omega_1$.
\end{itemize}
\end{definition}

\begin{theorem}\label{thm:DI}
$\sf{ZFC + DI}$ implies that the minimal model of ${\sf{ZF}}+{\sf{AD}}_\mathbb{R} + ``\Theta$ is a regular cardinal" exists.
\end{theorem}


As was mentioned in the abstract, a theory $T$ extending ${\sf{ZF}}+{\sf{AD}}_\mathbb{R}$ has a minimal model if it has a transitive model $M$ containing the reals and ordinals such that it is contained in any other transitive model of $T$ containing the ordinals and the reals. The proof of \cite[Theorem 6.26]{ATHM} explicitly establishes that the existence of divergent models of $\sf{AD^+}$ implies their common part is beyond a model of ${\sf{ZF}}+{\sf{AD}}_\mathbb{R} + ``\Theta$ is a regular cardinal". Thus, if there is a model of ${\sf{ZF}}+{\sf{AD}}_\mathbb{R} + ``\Theta$ is a regular cardinal" then there is a minimal one.  As a result of this theorem and Woodin's unpublished work, which we will present in Section \ref{sec:upper_bound}, we obtain the following equiconsistency result. 
\begin{theorem}\label{thm:HI}
The following theories are equiconsistent.
\begin{enumerate}
\item ${\sf{ZFC}}+\sf{DI}$
\item ${\sf{ZF}}+{\sf{AD}}_\mathbb{R} + ``\Theta$ is a regular cardinal."
\end{enumerate}
\end{theorem}

Let $NS_{\omega_1}$ be the nonstationary ideal on $\omega_1$ and $(T)$ be the theory 
\begin{center}
$\sf{ZFC + CH \ + }$ ``$NS_{\omega_1}|S$ is $\omega_1$-dense for a dense set of $S\in \powerset(\omega_1)\slash NS_{\omega_1}".$ 
\end{center}
Woodin, unpublished, has shown that Con($(T)$) follows from Con($\sf{ZF} + \sf{AD}_\mathbb{R} + ``$$\Theta$ is regular."). This result and Theorem \ref{thm:HI} immediately show
\begin{corollary}\label{cor:equi2} 
The following theories are equiconsistent.
\begin{enumerate}
\item $\sf{ZF} + \sf{AD}_\mathbb{R} + $$``\Theta$ is regular".
\item (T).
\end{enumerate}
\end{corollary}
\noindent This confirms that theories (b) and (c) in \cite[Question 12]{Woodin} are indeed equiconsistent.  Below we give some more motivations for proving such theorems.\\\\
\noindent \textbf{Some definitions and more results.}\\\\
For any set $X$, let $\powerset_{\omega_1}(X)$ be the set of countable subsets of $X$. Let $\mathcal{I}$ be an ideal on $\powerset_{\omega_1}(\mathbb{R})$. We let $\mathcal{I}^+$ and $\mathcal{F}_\mathcal{I}$ be as before and let $\mathbb{P}_\mathcal{I}$ be the boolean algebra $\powerset(\powerset_{\omega_1}(\mathbb{R}))\slash \mathcal{I}$. Let $\mathfrak{c}$ denote the size of the continuum.

\begin{definition}\label{dfn:prec_ideal}
An ideal $\mathcal{I}$ on $\omega_1$ or on $\powerset_{\omega_1}(\mathbb{R})$ is \textit{precipitous} if whenever $G\subseteq \mathbb{P}_{\mathcal{I}}$ is a $V$-generic ultrafilter, the generic ultrapower Ult$(V,G)$ induced by $G$ is well-founded.
\end{definition}

\begin{definition}\label{dfn:strong_ideal}
An ideal $\mathcal{I}$ on $\powerset_{\omega_1}(\mathbb{R})$ is \textit{strong} if
\begin{enumerate}[(a)]
\item $\mathcal{I}$ is precipitious, and
\item whenever $G\subseteq\mathbb{P}_\mathcal{I}$ is $V$-generic, letting $j_G: V\rightarrow \textrm{Ult}(V,G)$ be the ultrapower map, then $j_G(\omega_1) = \mathfrak{c}^+$.
\end{enumerate}
\end{definition}

\begin{definition}\label{dfn:pseudo_ideal}
An ideal $\mathcal{I}$ on $\powerset_{\omega_1}(\mathbb{R})$ is \textit{pseudo-homogeneous} if
for every $\alpha \in \ord$, $s\in \ord^\omega$, $\lambda < \mathfrak{c}^+$, and formula $\theta$ in the language of set theory, letting $G\subseteq \mathbb{P}_\mathcal{I}$ be a $V$-generic filter and $j_G:V\rightarrow \textrm{Ult}(V,G)$ the corresponding ultrapower map, the truth of the statement 
\begin{center}
Ult$(V,G)\models \theta[\alpha, j_G(s), j_G[\lambda^\omega]]$
\end{center}
is independent of the choice of $G$.
\end{definition}

%
%
We obtain an equiconsistency regarding strong, pseudo-homogeneous ideals on $\powerset_{\omega_1}(\mathbb{R})$.
\begin{theorem}\label{cor:equiconsistency}
The following are equiconsistent.
\begin{enumerate}
\item $\sf{ZFC} +{}$``The nonstationary ideal on $\powerset_{\omega_1}(\mathbb{R})$ is strong and pseudo-homogeneous.''
\item $\sf{ZF} + \sf{AD}_\mathbb{R} + $``$\Theta$ is a regular cardinal.''
\end{enumerate}
\end{theorem}

\noindent \textbf{Motivations}\\\\
Motivated by the success of the generic elementary embeddings induced by ideals or other similar structures, Foreman has suggested them as a possible foundational framework, and exposited his ideas in \cite[Chapter 11]{Foreman}. As is well known, the basic foundational issue that set theory is facing is its inability to produce a single foundational framework that is accepted by all and at the same time solves all fundamental problems including the ${\sf{Continuum\ Hypothesis}}$. Several successful foundational frameworks, such as ${\sf{Forcing\ Axioms}}$, ${\sf{Canonical\ Inner\ Models}}$ and ${\sf{Generic\ Embeddings}}$, have been proposed and developed, but they all seem to disagree on basic questions such as whether the  ${\sf{Continuum\ Hypothesis}}$ is true or whether the universe is a ground (i.e., cannot be obtained as a non-trivial forcing extension of an inner model) and on many other such fundamental questions.

One of the main goals of ${\sf{CMI}}$ is to unify all of these frameworks by showing that each can be naturally interpreted in another. Given such bi-interpretations, disagreements on fundamental questions can be traced to subjective preferences in one framework over another, or preferences in one type of formalism over another. 

For example, Woodin's theorem (\rthm{woodin}) and \rthm{thm:HI} show how to interpret natural ideas occurring in  the study of generic embeddings in models of determinacy and vice versa. The reason is that, in both cases, the forcing notion used to obtain the models carrying such ideals are natural forcing notions, and in the other direction, the models of determinacy built in both cases are natural canonical models of ${\sf{AD}}$. This sort of bi-interpretability demonstrates that one cannot have scientifically objective reasons for preferring generic embeddings over, say, determinacy axioms, as they are deeply interconnected: commitment to one entails commitment to the other. A bias towards a particular formalism can be justified by other more pragmatic ways, for example by insisting on the shortest or clearest or most natural possible proofs of certain desired theorems. The ideas exposited above are the motivational ideas behind proving theorems like the main theorem of this paper.\footnote{The authors first learned about these ideas from John Steel.} For a more detailed discussion of ${\sf{CMI}}$ and its role in set theory, readers may consult \cite{SargTrang}.\\\\
\textbf{The history behind the paper.}\\\\
The first written presentation of ${\sf{CMI}}$ is Ketchersid's PhD thesis \cite{ketchersid2000toward}, which motivated  Ralf Schindler and John Steel to work on a book presenting the ${\sf{Core\ Model\ Induction}}$ (see \cite{CMI}). In 2006 they organized a seminar in Berlin covering the basics of ${\sf{CMI}}$. As one can see by flipping through \cite{CMI}, one of the main directions pursued by the community at this time was to complete Ketchersid's project. See John Steel's \cite{DMATM} for a conjecture along the same vein.

One of the main reasons this was believed to be important was that it was not known and still is not known how to force ${\sf{DI}}$, clause 1 of \rthm{thm:HI}, from conventional large cardinals that are weaker than supercompact cardinals. Woodin forced ${\sf{DI}}$ both over the models of ${\sf{AD}}_{\mathbb{R}}+``\Theta$ is a regular cardinal'' and from an almost huge cardinal (see \cite[Chapter 7.14]{Foreman}). In \cite{Woodin}, Woodin also forced ${\sf{MM}}^{++}(\mathfrak{c})$, Martin's Maximum for forcing posets of size at most the continuum, over a model of ${\sf{AD}}_{\mathbb{R}}+``\Theta$ is a regular cardinal'' (see \cite[Theorem 9.40]{Woodin}), and just like with ${\sf{DI}}$, it is not known how to force ${\sf{MM}}^{++}(\mathfrak{c})$ from conventional large cardinals much weaker than a supercompact cardinal. These and other results of Woodin from \cite{Woodin} seem to suggest that the theory ${\sf{AD}}_{\mathbb{R}}+``\Theta$ is a regular cardinal'' is in the region of supercompact cardinals, and the project of getting a model of it via ${\sf{CMI}}$ seemed to be equivalent to getting canonical inner models that could have supercompact cardinals in it, which has been one of the Holy Grails of set theory. 

However, \cite{ATHM} showed that in fact the theory  ${\sf{AD}}_{\mathbb{R}}+``\Theta$ is a regular cardinal'' is much weaker than a supercompact cardinal: it is weaker than a Woodin cardinal that is a limit of Woodin cardinals (see \cite[Theorem 6.26]{ATHM}). This theorem seems to suggest the existence of a gap  in our understanding of models of set theory. On the one hand, the conventional forcing and large cardinal technology that is needed to force statements such as ${\sf{DI}}$ or ${\sf{MM}}^{++}(\mathfrak{c})$ requires the complexity of a supercompact cardinal or beyond, and on the other hand, equally natural but different technologies based on \cite{Woodin} place the complexity far below a supercompact cardinal. This phenomenon has not yet found a proper explanation.

While \cite{ATHM} did show that finishing Ketchersid's project will not lead to one of the Holy Grails of set theory, the importance of the project didn't diminish, as it was perceived to be one of the main guiding problems for developing the ${\sf{CMI}}$ to a technique for producing models of ${\sf{AD}}_{\mathbb{R}}+``\Theta$ is a regular cardinal'' and beyond\footnote{See for example \cite{LSA} for an analysis of determinacy models stronger than those of ${\sf{AD}}_{\mathbb{R}}+``\Theta$ is a regular cardinal'' and core model induction techniques for constructing such models from strong theories like $\sf{PFA}$.}. In this direction, the last chapter of the second author's thesis \cite{ATHMT} gave a rough outline of producing models of  ${\sf{AD}}_{\mathbb{R}}+``\Theta$ is a regular cardinal'' from a strengthening of $\sf{DI}$,\footnote{\label{lab:amenable}The strengthening is $\sf{DI}$ plus the statement: letting $\mathcal{I}$ be an $\omega_1$-dense ideal $\omega_1$, the generic embedding induced by any generic $G\subseteq \mathbb{P}_{\mathcal{I}}$ when restricted to the ordinals is amenable to $V$.} but later on a substantial error was discovered in the proof by Steel and the third author. The concept of \textit{embeddings with condensation} introduced in \cite{sargsyanCovering2013} (see \cite[Definition 11.14, Lemma 11.15]{sargsyanCovering2013}) and further developed in \cite{Trang2015PFA} (see \cite[Definition 3.81, Lemma 3.82]{Trang2015PFA}) and \cite{LSA} seemed good enough for correcting the aforementioned error, which is what we will do in this paper (see \rthm{thm:cond}). However, to obtain Theorem \ref{thm:DI}, more substantial ideas beyond this need to be developed.

Furthermore, the fourth author, in his thesis \cite{wilson2012contributions}, developed techniques for handling the successor stages of ${\sf{CMI}}$ that avoid the famous ``$A$-iterability" proofs (see \cite[Theorem 5.4.8]{CMI} or \cite[Theorem 1.46]{PFA}) and various other complicated arguments originally due to Woodin. We adapt the third author's arguments to our current context (see \rsec{sec:suc}). The second, third, and fourth authors established the consistency of $\sf{ZF + AD_\mathbb{R} \ + } ``\Theta$ is regular" from the aforementioned strengthening of $\sf{DI}$ in 2020-2021; this completes the project started by Ketchersid in \cite{ketchersid2000toward}. The obvious question is how to get rid of the technical assumption used in the second author's thesis, as mentioned in Footnote \ref{lab:amenable}. The first and fifth authors joined the ongoing work in 2022 and finished the project. The result of these collaborations is Theorem \ref{thm:DI}.

As mentioned above, it is a well-known unpublished theorem of Woodin that one can force ${\sf{DI}}$ over models of ${\sf{AD}}_{\mathbb{R}}+``\Theta$ is a regular cardinal". The fourth author forced some more general statements about ideals in his thesis, and we will use his argument to give a proof of this theorem of Woodin in Subsection \ref{subsec:upper} below. Thus, this paper presents a self-contained proof of \rthm{thm:HI}, giving the proof of both directions in as much detail as it is possible to do in a research article.  \rthm{thm:HI} and \rthm{cor:equiconsistency} are currently the only known equiconsistency results at the level of ${\sf{AD}}_{\mathbb{R}}+``\Theta$ is a regular cardinal". 


In Section \ref{sec:upper_bound}, we summarize basic facts about ideals and $\sf{AD}^+$ we need in this paper and show that $\sf{DI}$ and the existence of a strong, pseudo-homogeneous ideal on $\powerset_{\omega_1}(\mathbb{R})$ are consistent relative to $\sf{AD}_\mathbb{R} +{}$``$ \Theta$ is a regular cardinal.''\footnote{We adapt the proof given in the fourth author's thesis here. We note the result that Con$(\sf{ZFC + DI})$ follows from Con$(\sf{AD}_\mathbb{R} + {}$``$\Theta$ is regular") is due to Woodin.} In Section \ref{sec:prelim}, we summarize preliminaries and basic notions we need for ${\sf{CMI}}$. Section \ref{sec:outline} outlines the proof of Theorem \ref{thm:DI}. Sections \ref{sec:suc} and \ref{sec:lim} fill in the details of the outline and complete the proof of Theorem \ref{thm:DI}, obtaining models of ``$\sf{AD}_\mathbb{R} + $$\Theta$ is a regular cardinal" from $\sf{ZFC + DI}$. In Section \ref{sec:lower_bound_2}, we outline the argument obtaining models of ``$\sf{AD}_\mathbb{R} + $$\Theta$ is a regular cardinal" from the assumption that the nonstationary ideal on $\powerset_{\omega_1}(\mathbb{R})$ is strong and pseudo-homogeneous. Since the argument is very similar to the argument from $\sf{DI}$, we simply focus on the main changes, leaving the details to the reader. In the following, we will often write ``$\Theta$ is regular" for ``$\Theta$ is a regular cardinal."\\

\noindent \textbf{Acknowledgments.} The work here is greatly influenced by Ketchersid's work in his thesis \cite{ketchersid2000toward}, which in turn is greatly influenced by Woodin's early work in the $\sf{CMI}$. We are grateful to them for their inspiring work in this direction. We are also grateful to Woodin for his permission to include the proof of his unpublished work which shows that Con($\sf{ZF} + \sf{AD}_\mathbb{R} + $``$\Theta$ is a regular cardinal") implies Con$(\sf{ZFC}+\sf{DI}$). The third author is grateful to the NSF for its generous support via Career Award DMS-1945592.

\section{DENSE IDEALS AND STRONG PSEUDO-HOMOGENEOUS IDEALS FROM MODELS OF $\sf{AD}_\mathbb{R} + \Theta$ IS REGULAR}\label{sec:upper_bound}
In this section, we show the consistency of $\sf{ZFC + DI}$ and of the existence of a strong, pseudo-homogeneous ideal on $\PR$ from $\sf{AD}_\mathbb{R} + {}$``$\Theta$ is regular.''\footnote{In fact we show the nonstationary ideal on $\powerset_{\omega_1}(\mathbb{R})$ has these properties.} We first review basic facts about $\sf{AD}^+$ and ideals. In Subsection \ref{subsec:upper}, we will give the consistency proof.

\subsection{Basic facts about $\sf{AD}^+$}
We start with the definition of Woodin's theory of $\sf{AD}$$^+$. In this paper, we identify $\mathbb{R}$ with $\omega^{\omega}$. We use $\Theta$ to denote the sup of ordinals $\alpha$ such that there is a surjection $\pi: \mathbb{R} \rightarrow \alpha$. Under $\textsf{AC}$, $\Theta$ is just the successor cardinal of the continuum. In the context of $\sf{AD}$, the cardinal $\Theta$ is shown to be the supremum of $w(A)$\footnote{$w(A)$ is the Wadge rank of $A$.} for $A\subseteq \mathbb{R}$ (cf.\ \cite{solovay1978independence}). The definition of $\Theta$ relativizes to any determined pointclass $\Gamma$ with sufficient closure properties, and we may write $\Theta^\Gamma$ for the supremum of ordinals $\alpha$ such that there is a surjection from $\mathbb{R}$ onto $\alpha$ coded by a set of reals in $\Gamma$.
\begin{definition}
\label{AD+}
$\sf{AD}^+$ is the theory $\sf{ZF} + \sf{AD} + \sf{DC}_{\mathbb{R}}$ plus the following two statements: 
\begin{enumerate}
\item For every set of reals $A$, there are a set of ordinals $S$ and a formula $\varphi$ such that $x\in A \iff L[S,x] \vDash \varphi[S,x]$. The pair $(S,\varphi)$ is called an $\infty$-Borel code\index{$\infty$-Borel code} for $A$.
\item For every $\lambda < \Theta$, every continuous $\pi: \lambda^\omega \rightarrow \omega^{\omega}$, and every set of reals $A$, the set $\pi^{-1}[A]$ is determined.
\end{enumerate}
\end{definition}
\noindent $\sf{AD}$$^+$ is equivalent to $\sf{AD} + {}$``the set of Suslin cardinals is closed below $\Theta$.'' Another, perhaps more useful, characterization of $\sf{AD}^+$ is $\sf{AD} + {}$``$\Sigma_1$ statements reflect into the Suslin co-Suslin sets'' (see \cite{steelderived} for the precise statement). 
\\
\indent For $A\subseteq \mathbb{R}$, we let $\theta_A$\index{$\theta_A$} be the supremum of all $\alpha$ such that there is an $OD(A)$ surjection from $\mathbb{R}$ onto $\alpha$. If $\Gamma$ is a determined pointclass and $A\in \Gamma$, we write $\Gamma\rest A$ for the set of all $B\in\Gamma$ that are Wadge reducible to $A$. If $\alpha < \Theta^\Gamma$, we write $\Gamma\rest \alpha$ for the set of all $A\in \Gamma$ with Wadge rank strictly less than $\alpha$.
\begin{definition}[$\sf{AD}^+$]
\label{Solovaysequence}
The \textbf{Solovay sequence} is the sequence $\langle\theta_\alpha \ | \ \alpha \leq \lambda\rangle$ where
\begin{enumerate}
\item $\theta_0$ is the supremum of ordinals $\beta$ such that there is an $OD$ surjection from $\mathbb{R}$ onto $\beta$;
\item if $\alpha>0$ is limit, then $\theta_\alpha = \sup\{\theta_\beta \ | \ \beta<\alpha\}$;
\item if $\alpha =\beta + 1$ and $\theta_\beta < \Theta$ (i.e. $\beta < \lambda$), fixing a set $A\subseteq \mathbb{R}$ of Wadge rank $\theta_\beta$, $\theta_\alpha$ is the sup of ordinals $\gamma$ such that there is an $OD(A)$ surjection from $\mathbb{R}$ onto $\gamma$, i.e. $\theta_\alpha = \theta_A$.
\end{enumerate}
\end{definition}

Note that the definition of $\theta_\alpha$ for $\alpha = \beta+1$ in Definition \ref{Solovaysequence} does not depend on the choice of $A$.  One can also make sense of the Solovay sequence of pointclasses that may not be constructibly closed. Such pointclasses show up in core model induction applications. The Solovay sequence $(\theta_\alpha: \alpha <\gamma)$ of a pointclass $\Omega$ with the property that if $A\in\Omega$, then $L(A,\mathbb{R})\models \sf{AD}^+$ and $\powerset(\mathbb{R})\cap L(A,\mathbb{R})\subseteq \Omega$ is defined as follows. First, $\theta_0$ is the supremum of all $\alpha$ such that there is some $A\in\Omega$ and some $OD^{L(A,\mathbb{R})}$ surjection $\pi:\mathbb{R}\rightarrow \alpha$. If $\lambda<\gamma$ is limit, then $\theta_\gamma = \textrm{sup}_{\alpha<\lambda} \theta_\alpha$. If $\theta_\alpha$ has been defined and $\alpha+1<\gamma$, then letting $A\in\Omega$ be of Wadge rank $\theta_\alpha$, $\theta_{\alpha+1}$ is the supremum of $\beta$ such that there is some $B\in\Omega$ and some $OD(A)^{L(B,\mathbb{R})}$ surjection $\pi:\mathbb{R}\rightarrow \beta$.

Roughly speaking, the longer the Solovay sequence is, the stronger the associated $\sf{AD}$$^+$-theory is. The minimal model of $\sf{AD}^+$ is $L(\mathbb{R})$, which satisfies $\Theta = \theta_0$. The theory $\sf{AD}^+ + \sf{AD}_\mathbb{R}$ implies that the Solovay sequence has limit length. The theory $\sf{AD}$$_\mathbb{R} + \textsf{DC}$ is strictly stronger than $\sf{AD}$$_\mathbb{R}$ since by \cite{solovay1978independence}, $\textsf{DC}$ implies cof$(\Theta)>\omega$ whereas the minimal model\footnote{From here on, whenever we talk about ``models of $\textsf{AD}^+$", we always mean transitive models of $\sf{AD}^+$ that contain all reals and ordinals.} of $\sf{AD}$$_\mathbb{R}$ satisfies $\Theta = \theta_\omega$. The theory  ``$\sf{AD}$$_\mathbb{R} + \Theta$ is regular" is much stronger still, as it implies the existence of many models of $\sf{AD}$$_\mathbb{R} + \textsf{DC}$. We end this section with a theorem of Woodin, which produces models with Woodin cardinals from $\sf{AD}$$^+$. The theorem is important in the HOD analysis of such models.

\begin{theorem}[Woodin, see \cite{koellner2010large}]
Assume $\sf{AD}$$^+$. Let $\langle \theta_\alpha \ | \ \alpha \leq \Omega\rangle$ be the Solovay sequence. Suppose $\alpha = 0$ or $\alpha = \beta+1$ for some $\beta < \Omega$. Then $\rm{HOD} $ $ \vDash \theta_\alpha$ is Woodin.
\end{theorem}

\subsection{Basic properties of ideals}\label{sec:ideal_properties}

We summarize standard facts about ideals that we will need in this paper. See for example \cite{Woodin} and \cite{Jech} for a more detailed discussion.

Suppose $\mathcal{I}$ is an ideal on a set $X$. We say that $\mathcal{I}$ is \textit{countably complete} if whenever $\{A_n : n<\omega\}$ are sets in $\mathcal{I}$ then $\bigcup_{n<\omega} A_n \in \mathcal{I}$. Supposing $X$ is a cardinal (e.g. $X=\omega_1$), we say $\mathcal{I}$ is \textit{normal} if whenever $\{A_x : x\in X\}\subset \mathcal{I}$ then the diagonal union $\{x\in X: \exists y\in x (x\in A_y) \}\in \mathcal{I}$. All ideals $\mathcal{I}$ on a cardinal considered in this paper will be assumed countably complete and normal.

Suppose $\mathcal{I}$ is an $\omega_1$-dense ideal on $\omega_1$. The following are standard facts; see \cite[Definition 6.19]{Woodin} and the discussion after it.
\begin{fact}\label{fact:dense_ideal_facts}
\begin{enumerate}[(i)]
\item $\PI$ is a homogeneous forcing.\footnote{A forcing $\mathbb{P}$ is homogeneous if whenever $p,q\in\mathbb{P}$, there is an automorphism $\sigma: \mathbb{P}\rightarrow \mathbb{P}$ such that $\sigma(p)$ is compatible with $q$.}
\item  There is a boolean isomorphism $\pi:\mathbb{P}_{\mathcal{I}} \rightarrow \textrm{RO}(\textrm{Coll}(\omega,\omega_1))$\footnote{$\textrm{RO}(\textrm{Coll}(\omega,\omega_1))$ is the regular open algebra of Coll$(\omega,\omega_1)$.}.  In particular, $\mathbb{P}_{\mathcal{I}}$ is forcing equivalent to Coll$(\omega,\omega_1)$.
\item For any $V$-generic filter $G\subset \textrm{Coll}(\omega,\omega_1)$, $\pi$ induces a $V$-generic filter $H\subset \mathbb{P}_{\mathcal{I}}$, and letting $j: V\rightarrow M =_{def} Ult(V,H) \subset V[H]$ be the associated generic ultrapower map, we have:
\begin{enumerate}[(a)]
\item $j(f)(\omega_1^V) = G$ for some $f:\omega_1\rightarrow H_{\omega_1}$; in particular, $V[H] = V[G]$.
\item $j(\omega_1^V) = \omega_2^V$.
\item $M$ is well-founded and $M^\omega \subset M$ in $V[H]$.
\end{enumerate}
\end{enumerate}
\end{fact}

Let $\mathcal{I}$ be an $\omega_1$-dense ideal on $\omega_1$. For any $V$-generic $g\subset \mathbb{P}_{\mathcal{I}}=_{def} \powerset(\omega_1)\slash \mathcal{I}$, let $j_g: V\rightarrow M = \textrm{Ult}(V,g)$ be the associated ultrapower map. We fix a Boolen isomorphism $\pi:\mathbb{P}_{\mathcal{I}} \rightarrow \textrm{RO}(\textrm{Coll}(\omega,\omega_1))$ as in Fact \ref{fact:dense_ideal_facts} and let $G\subset Coll(\omega,\omega_1)$ be such that $g$ is induced from $G$ via $\pi$. When $g$ is clear from the context, we will write $j$ for $j_g$. 

We say that a set of reals $A$ is \textit{$\omega_1$-universally Baire} (or $\omega_1$-UB) if there is some ordinal $\gamma$ and a pair of trees $T, U$ on $\omega\times \gamma$ such that $A = p[T] = \mathbb{R} - p[U]$ and for any forcing $\mathbb{P}$ of size $\leq \omega_1^V$, for any $V$-generic $h\subset\mathbb{P}$, in $V[h]$, $p[T] = \mathbb{R} - p[U]$. Here $p[T] = \{x\in\mathbb{R} : \exists f\in \gamma^\omega \ (x,f)\in [T]\}$.

\begin{lemma}\label{lem:UBprojection}
Let $\mathbb{P}_{\mathcal{I}}, g, G, M$ be as above. Suppose $A\subset \mathbb{R}$ is $\omega_1$-UB as witnessed by trees $(T,U)$, then in $V[G]$, $p[T] =  p[j(T)]$ and $p[U] = p[j(U)]$.
\end{lemma}
\begin{proof}
We write $j$ for $j_g$. Clearly, $p[T]\subseteq p[j(T)]$ and $p[U] \subseteq p[j(U)]$. In $M$, equivalently in $V[G]$, 
\begin{center}
$p[j(T)] = \mathbb{R} - p[j(U)]$.
\end{center}
This follows from elementarity of $j$, the fact that in $V$, $p[T] = \mathbb{R} - p[U]$, and property (c) of Fact \ref{fact:dense_ideal_facts}. 

By the fact that $T,U$ witness $A$ is $\omega_1$-UB and $Coll(\omega,\omega_1)$ has size $\omega_1$, in $V[G]$, $p[T] = \mathbb{R} - p[U]$. We must then get $p[T] = p[j(T)]$ and $p[U]  = p[j(U)]$.
\end{proof}

Suppose  $X = \powerset_{\omega_1}(Y)$, where  $\powerset_{\omega_1}(Y)$ is the collection of all countable subsets of $Y$, for some set $Y$ (e.g. $Y = \mathbb{R}$).  We say $\mathcal{I}$ is \textit{fine} if for any $y\in Y$, the set $\{\sigma \in \powerset_{\omega_1}(Y): y\notin \sigma\}\in\mathcal{I}$. We say $\mathcal{I}$ is \textit{normal} if whenever $\{A_y : y\in Y\}\subset \mathcal{I}$, the diagonal union $\{\sigma\in \powerset_{\omega_1}(Y): \exists y\in \sigma \ (\sigma\in A_y)\}\in \mathcal{I}$. $\mathcal{I}$ is $|Y|$-dense if there is a dense subset of $\PI$ of size $|Y|$. All ideals on sets of the form $\powerset_{\omega_1}(Y)$ considered in this paper will be assumed countably complete, normal, and fine.

\begin{lemma}\label{lem:quasi_hom}
Suppose $\mathcal{I}$ is a pseudo-homogeneous ideal on $\PR$. Let $G\subset \PI$ be $V$-generic and let $j_G:V \rightarrow Ult(V,G)$ be the associated generic embedding. Then:
\begin{enumerate}[(a)]
\item For any ordinal $\alpha$, $j_G\rest \alpha$ does not depend on $G$; in particular, $j_G\rest \alpha\in V$.
\item If $\lambda<\mathfrak{c}^+$, then $j_G[\lambda^\omega]$ does not depend on $G$ and $j_G[\lambda^\omega]\in V$. 
\item If $A$ is a set of ordinals that is definable in $V$ from a countable sequence of ordinals, then $j_G(A)$ does not depend on $G$ and $j_G(A)\in V$.
\end{enumerate}
\end{lemma}
\begin{proof}
We give the proof for (a). The other items are similar. Let $\theta(u,v,w)$ be the formula ``$u = v(0)$". Let $\alpha$ be an ordinal. Let $s:\omega \rightarrow \textrm{Ord}$ be the constant function $s(n) = \alpha$ for all $n\in\omega$. For each ordinal $\beta$ the truth of the statement Ult$(V,G)\models \theta[\beta, j_G(s), \emptyset]$ is independent of $G$ by pseudo-homogeneity, so the value of $j_G(\alpha)$ is independent of $G$.
\end{proof}

\subsection{Ideals from determinacy}\label{subsec:upper}
We assume $\sf{AD}_\mathbb{R} + {}$``$\Theta$ is regular'' and $V = L(\powerset(\mathbb{R}))$.
Let $\mathbb{P}$ be a poset with the following properties:
\begin{itemize}
\item $\mathbb{P}$ is coded by a set of reals.
\item $\mathbb{P}$ is $\sigma$-closed.
\item $\mathbb{P}$ is homogeneous.
\item $1 \Vdash_{\mathbb{P}} \mathbb{R}$ is wellorderable.
\item $1 \Vdash_{\mathbb{P}} \mathfrak{c}$-$\mathsf{DC}$, dependent choices for $\mathfrak{c}$-sequences. 
\end{itemize}
Examples of such $\mathbb{P}$ are $Coll(\omega_1,\mathbb{R})$ and $\mathbb{P}_{\textrm{max}}$.

Let $G\subseteq \mathbb{P}$ be $V$-generic and let $H\subset Coll(\Theta,\powerset(\mathbb{R}))^{V[G]}$. Note that by the properties of $\mathbb{P}$ and the assumption $V=L(\powerset(\mathbb{R}))$, in $V[G][H]$, $\sf{ZFC}$ holds and $\Theta= \mathfrak{c}^+$.

\begin{definition}\label{dfn:ord_covering}
In $V[G][H]$ an ideal $\mathcal{I}$ on $\powerset_{\omega_1}(\mathbb{R})$ is said to have the \textit{ordinal covering property} with respect to $V$ if for every function $F: \powerset_{\omega_1}(\mathbb{R})\rightarrow \textrm{Ord}$ and every $\mathcal{I}$-positive set $S$, there is some $\mathcal{I}$-positive set $S_0\subseteq S$ and some $F_0: \powerset_{\omega_1}(\mathbb{R})\rightarrow \textrm{Ord}$ in $V$ such that $F\rest S_0 = F_0\rest S_0$.
\end{definition}

We will show that in $V[G][H]$, there is an ideal $\mathcal{I}$ with the ordinal covering property with respect to $V$. Let $\mu$ be the Solovay measure on $\powerset_{\omega_1}(\mathbb{R})^V$, so $A\in \mu$ if and only if $A$ contains a club set in $\powerset_{\omega_1}(\mathbb{R})$. A set $A$ is club if and only if there is a function $F: \mathbb{R}^{<\omega}\rightarrow \mathbb{R}$ such that 
\begin{center}
$\sigma\in A \Leftrightarrow F[\sigma]\subseteq \sigma$.  
\end{center}
We say that $A$ is the \textit{club set generated by $F$.}

The measure $\mu$ induces an ultrapower map on the ordinals, $j_\mu: \textrm{Ord}\rightarrow \textrm{Ord}$. By the basic theory of $\sf{AD}^+$,
\begin{equation}\label{eqn:j_mu}
j_\mu(\omega_1) = \Theta.
\end{equation}
See, for example, \cite[Section 1.2]{wilson2012contributions} for a proof of this fact.

\begin{lemma}\label{lem:agreement}
Suppose $V,G,H$ are as above. Suppose $\mathcal{I}$ is an ideal on $\powerset_{\omega_1}(\mathbb{R})$ with the ordinal covering property with respect to $V$. Let $K\subset\mathbb{P}_{\mathcal{I}}$ be a $V[G][H]$-generic filter. Then:
\begin{enumerate}[(a)]
\item The generic embedding $j_K\rest \textrm{Ord} = j_\mu\rest \textrm{Ord}$. In particular, $j_K\rest \alpha\in V[G][H]$ for every ordinal $\alpha$ and doesn't depend on the choice of $K$.
\item $\mathcal{I}$ is strong.
\end{enumerate}
\end{lemma}
\begin{proof}
For (a), for any $F: \powerset_{\omega_1}(\mathbb{R})\rightarrow \textrm{Ord}$ in $V[G][H]$, the covering property gives some $S\in K$ and $F_0\in V$ such that $F\rest S = F_0\rest S$. Also, $K\cap V = \mu$ since $K$ is normal; this gives 
\begin{center}
$\{F: \powerset_{\omega_1}(\mathbb{R})\rightarrow \textrm{Ord}\}^{V[G][H]}\slash K = \{F: \powerset_{\omega_1}(\mathbb{R})\rightarrow \textrm{Ord}\}^V \slash \mu$
\end{center}
and $j_K\rest \textrm{Ord} = j_\mu \rest \textrm{Ord}$. Part (b) follows from (a) and \eqref{eqn:j_mu}.
\end{proof}

\begin{lemma}\label{lem:ord_covering}
In $V[G][H]$, if $\mathcal{I}$ has the ordinal covering property relative to $V$, then $\mathcal{I}$ is pseudo-homogeneous.
\end{lemma}
\begin{proof}
Let $K\subset \mathbb{P}_\mathcal{I}$ be a $V[G][H]$-generic filter.  Let  $\alpha \in \textrm{Ord}$, $s\in \textrm{Ord}^\omega$, $\lambda < \mathfrak{c}^+$, and let $\theta$ be a formula in the language of set theory. It suffices to show that the statement Ult$(V[G][H],K)\models \theta[\alpha, j_K(s), j_K[\lambda^\omega]]$ is independent of $K$.
 By the ordinal covering property, we can find $F_0\in V$ that represents $\alpha$ in both Ult$(V,\mu)$ and Ult$(V[G][H],K)$. In both ultrapowers, $j(s)$ is represented by the constant function $F_1(\sigma) = s$ for all $\sigma\in \powerset_{\omega_1}(\mathbb{R})$. Fix a surjection $\pi: \mathbb{R}\rightarrow \lambda^\omega$ in $V$.  Then $j_K[\lambda^\omega]$ is represented by the function $F_2\in V$ given by $F_2(\sigma) = \pi[\sigma]$. So we have Ult$(V[G][H],K)\models \theta[\alpha, j_K(s), j_K[\lambda^\omega]]$  if and only if the set
\begin{center}
$S = \{\sigma: V[G][H] \vDash \theta[F_0(\sigma), F_1(\sigma), F_2(\sigma)\}$
\end{center}
is in $K$.
By homogeneity of $\mathbb{P}$, $S\in V$. But then we have
$S\in K$ if and only if $S\in \mu$, as desired.
\end{proof}

\begin{theorem}\label{thm:ord_covering}
In $V[G][H]$, the nonstationary ideal $\mathcal{I} = NS_{\omega_1,\mathbb{R}}$ on $\powerset_{\omega_1}(\mathbb{R})$  has the ordinal covering property with respect to $V$.
\end{theorem}

To establish the covering property of $\mathcal{I}$ in $V[G][H]$, or equivalently in $V[G]$, we will need the following lemma.

\begin{lemma}
Let $\dot{S}$ be a $\mathbb{P}$-name for a subset of $\powerset_{\omega_1}(\mathbb{R})$. The following statements are equivalent for any given $p\in\mathbb{P}$:
\begin{enumerate}[(a)]
\item $p \Vdash {}$``$\dot{S}$ contains a club.''
\item For a club of $\sigma\in \powerset_{\omega_1}(\mathbb{R})$,
\begin{center}
$(\dag) \ \ \ \ \ \ \forall^* g\subset \mathbb{P}\rest \sigma \textrm{ containing } p \ \forall q \leq g \ q \Vdash \sigma \in \dot{S}$.
\end{center}
Here $\forall^* g$ stands for ``for a comeager set of filters $g$"\footnote{By $\mathbb{P}\rest \sigma$, we mean the set of conditions in $\mathbb{P}$ coded by a real in $\sigma$. Note that $\mathbb{P}\rest \sigma$ is countable, so the category quantifier over the set of all filters on it makes sense.} and $q\leq g$ means $\forall r\in g \ q \leq r$.
\end{enumerate}
\end{lemma}
\begin{proof}
Fix $p\in\mathbb{P}$. Assume (a) holds for $p$. Let $\dot{f}$ be a $\mathbb{P}$-name for a function from $\mathbb{R}^{<\omega}$ into $\mathbb{R}$ such that $p$ forces $\dot{S}$ to contain the club set generated by $\dot{f}$. We may assume $\mathbb{P}\subseteq \mathbb{R}$. To see (b), note that there is a club set of $\sigma$ such that for all $t\in \sigma^{<\omega}$, the set
\begin{center}
$D_t = \{q \in \mathbb{P}\cap \sigma : (\exists x\in \sigma) \ (q\Vdash \dot{f}(t) = x) \}$.
\end{center}
is dense below $p$ in $\mathbb{P}\cap \sigma$. This easily gives $(\dag)$ for $\sigma$ as there are countably many dense sets $D_t$ and hence there is a comeager set of filters $g\subset \mathbb{P}\cap \sigma$ meeting all the $D_t$'s.

Assume (b) holds for $p$. Let 
\begin{center}
$A = \{(q,x) : x \textrm{ codes } \sigma\in\powerset_{\omega_1}(\mathbb{R}) \textrm{ and } q \Vdash \sigma\in \dot{S}\}$.
\end{center}
Take $N = L_\alpha(P_{\beta}(\mathbb{R}))$ satisfying $\sf{ZF}^- + \sf{AD}_\mathbb{R} + {}$``$\Theta$ is regular", containing $A$, and admitting a surjection $F:\mathbb{R}\rightarrow N$.\footnote{Here $P_\beta(\mathbb{R})$ is the set $\{B\subset\mathbb{R} : B \textrm{ has Wadge rank less than } \beta \}$.} Let $B\subset \mathbb{R}$ code the first order theory of the structure $(V_{\omega+1},\in, A)$. Because $\sf{AD}_\mathbb{R}$ implies that every set of reals is $\mathbb{R}$-universally Baire (see e.g.\ \cite[Section 1.2]{wilson2012contributions},) in particular $A$ and $B$ are $\mathbb{R}$-universally Baire. There is then a club $C$ of $\sigma\in \powerset_{\omega_1}(N)$ having the following properties:
\begin{itemize}
\item $(\dag)$ holds for $\sigma\cap \mathbb{R}$.
\item $\sigma \prec N$.
\item Defining $\pi_\sigma: \sigma \rightarrow N_\sigma$ as the transitive collapse of $\sigma$, we have
\begin{center}
$(V_{\omega+1}\cap N_\sigma[h],\in, A\cap N_\sigma[h])\prec (V_{\omega+1},\in, A)$
\end{center}
for any $N_\sigma$-generic filter $h\subset Coll(\omega,\sigma\cap \mathbb{R})$.
\end{itemize}
The last item follows from the $\mathbb{R}$-universal Baireness of $B$.

All $\sigma\in C$ have the following property:
\begin{equation}\label{eqn:forced}
N_\sigma \vDash p \Vdash^g_{\mathbb{P}\rest(\mathbb{R}\cap\sigma)} (1 \Vdash^h_{Coll(\omega,\mathbb{R}\cap\sigma)} (\forall q\leq g)((q,\sigma_h)\in \pi_\sigma(A)_{g\times h})).
\end{equation}
In \eqref{eqn:forced}, $\sigma_h$ denotes the real generally coding $\sigma\cap \mathbb{R}$ relative to $h$ and $\pi_\sigma(A)_{g\times h}$ denotes the unique extension of $\pi_\sigma(A)$ to a set of reals in $N_\sigma[g][h]$, which can be construed as a generic extension of $N_\sigma$ by $Coll(\omega,\sigma\cap\mathbb{R})$; the extension is given by the universal Baireness of $A$.

Now suppose $G\subset \mathbb{P}$ is $V$-generic and $p\in G$. There is a club set $D$ of $\sigma\in C$ such that $\sigma[G]\prec N[G]$ and $\sigma[G]\cap V = \sigma$. Take a $\sigma$ in this club and $g = G\cap \sigma$. Note that any lower bound $q\leq g$ forces $\sigma\in \dot{S}$ by \eqref{eqn:forced} and there is $q\leq g$ in $G$; so $\sigma\cap \mathbb{R}\in \dot{S}_G$. Therefore, the club set $\{\sigma\cap \mathbb{R}: \sigma\in D\}$  witnesses (a).
\end{proof}

\begin{proof}[Proof of Theorem \ref{thm:ord_covering}]
Suppose $p_0$ forces ``$\dot{F}:\dot{S} \rightarrow \textrm{Ord}$ and $\dot{S}\subseteq \powerset_{\omega_1}(\mathbb{R})$ is stationary." Using $(\dag)$, the latter part of this statement is equivalent to the following statement. For stationary many (equivalently by $\sf{AD}_\mathbb{R}$, for club many) countable $\sigma\subset \mathbb{R}$, 
\begin{center}
$\exists^* g\subset \mathbb{P}\rest \sigma \textrm{ containing } p_0 \ \exists q \leq g \ q \Vdash \sigma\in \dot{S}$.
\end{center}
Under $\sf{AD}$, a well-ordered union of meager sets is meager, so let $F_0(\sigma)$ be the least $\alpha$ such that 
\begin{center}
$\exists^* g\subset \mathbb{P}\rest \sigma \textrm{ containing } p_0 \ \exists q \leq g \ q \Vdash \dot{F}(\sigma)=\alpha$.
\end{center}
By the above, $p_0$ forces that the set of $\sigma\in \dot{S}$ such that $F(\sigma)=F_0(\sigma)$ is stationary.
\end{proof}

Theorem \ref{thm:ord_covering} and Lemmas \ref{lem:agreement} and \ref{lem:ord_covering} immediately give one direction of Theorem \ref{cor:equiconsistency}.

\begin{corollary}\label{cor:pseudo_consistency}
Con$(\sf{ZF} + \sf{AD}_\mathbb{R} + {}$``$\Theta$ is regular'') implies Con($\sf{ZFC} + {}$``the nonstationary ideal on $\PR$ is strong and pseudo-homogeneous'').
\end{corollary}

Now we proceed to prove one direction of Theorem \ref{thm:HI}. We show Con$(\sf{AD}_\mathbb{R} + {}$``$\Theta$ is regular'') implies Con$(\sf{ZFC + DI}$). We fix objects $V,\mathbb{P}, G, H$ as before. The following is the main theorem.

\begin{theorem}\label{thm:dense_ord_covering}
In $V[G][H]$, there is a $\mathfrak{c}$-dense ideal on $\powerset_{\omega_1}(\mathbb{R})$ with the ordinal covering property relative to $V$.
\end{theorem}
We review some facts regarding generic ultrapowers by $Coll(\omega,\mathbb{R})$-generics. See \cite{wilson2012contributions} for a more detailed discussion. Let $h\subset Coll(\omega,\mathbb{R})$ be $V$-generic and 
\begin{center}
$U_h = \{A\subseteq \mathbb{R}^\omega : A \textrm{ is weakly comeager below some } p\in h\}$.
\end{center}
Here $A\subseteq \mathbb{R}^\omega$ is weakly comeager below a condition $p\in Coll(\omega,\mathbb{R})$ if for a club set of $\sigma\in \powerset_{\omega_1}(\mathbb{R})$, $A\cap \sigma^\omega$ is comeager below $p$ in $\sigma^\omega$.\footnote{We equip $\sigma^\omega$ with the product of the discrete topologies on $\sigma$, so it is homeomorphic to the Baire space.} $U_h$ is the generic ultrafilter on $\mathbb{R}^\omega$ induced by $h$. $U_h$ gives a generic embedding $j_h: V\rightarrow \textrm{Ult}(V,U_h) \subset V[h]$. Using the fact that $\sf{AD}_\mathbb{R}+{}$``$\Theta$ is regular'' holds in $V$, we can prove \L o\'{s}'s theorem for $j_h$ and hence  $j_h$ is elementary. We can show that the map $[F_0]_\mu \mapsto [F_0\circ \textrm{ran}]_{U_h}$ is an isomorphism from $\textrm{Ult}(\textrm{Ord},\mu)$ to $\textrm{Ult}(\textrm{Ord},U_h)$, $\mathbb{R}^{V[h]} = \mathbb{R}^{\textrm{Ult}(V,U_h)}$, and $j_\mu\rest \textrm{Ord} = j_{U_h}\rest \textrm{Ord}$.

\begin{proof}[Proof of Theorem \ref{thm:dense_ord_covering}]
We first prove the following claim.
\begin{claim}\label{claim:extending_generic}
If $h\subset Coll(\omega,\mathbb{R})$ is a $V[H]$-generic filter such that $G\in V[h]$, then letting $j_h:V\rightarrow \textrm{Ult}(V,U_h)\subset V[h]$ denote the corresponding elementary embedding, in $V[h][H]$, there is an Ult$(V,U_h)$-generic filter $G'\subset j_h(\mathbb{P})$ extending $j_h``G$.
\end{claim}
\begin{proof}
The poset
$j_h(\mathbb{P})$ is countably closed in Ult$(V,U_h)$ and is coded by a set of reals there. 
In $V[h]$, because $\mathbb{R}\cap V[h] = \mathbb{R}\cap \textrm{Ult}(V,U_h)$ the poset $j_h(\mathbb{P})$ remains countably closed, and because $j_h``G$ is countable there is a lower bound $p\in j_h(\mathbb{P})$ for $j_h``G$.

Now note that in $V[h]$, there is a surjection $f$ from $\powerset(\mathbb{R})^V$ onto $\powerset(j_h(\mathbb{P}))^{\textrm{Ult}(V,U_h)}$; this is because every subset of $j_h(\mathbb{P})$ in Ult$(V,U_h)$ is represented by a function $\mathbb{R}^\omega\rightarrow \powerset(\mathbb{R})$ in $V$, which can be coded by a set of reals in $V$. In $V[G][H]$, there is a surjection $k$ from $\omega_1^{V[h]} = \Theta^V$ onto $\powerset(\mathbb{R})$ whose proper initial segments are in $V[G]\subset V[h]$; this follows from the fact that the forcing $Coll(\Theta^V,\powerset(\mathbb{R})^V)^{V[G]}$ is $\mathfrak{c}^+$-closed and $V[G]$ satisfies $\mathfrak{c}$-$\sf{DC}$. Then the surjection $k\circ f: \omega_1^{V[h]} \rightarrow \powerset(j_h(\mathbb{P}))^{\textrm{Ult}(V,U_h)}$ has the property that its proper initial segments are in $V[h]$.\footnote{We need this property for the following argument because this is the model in which $j_h(\mathbb{P})$ is countably closed.} Using this surjection, we recursively define a decreasing $\omega_1$-sequence of conditions $(p_\alpha: \alpha<\omega_1)$ in $j_h(\mathbb{P})$ below $p$ whose proper initial segments are in $V[h]$ and which generates the desired filter $G'$.
\end{proof}

By the assumptions on $\mathbb{P}$, $\mathbb{P}\times Coll(\omega,\mathbb{R})$ is forcing equivalent to $Coll(\omega,\mathbb{R})$; therefore, we can find an $h$ satisfying the hypothesis of Claim \ref{claim:extending_generic}. By Claim \ref{claim:extending_generic}, forcing with $Coll(\omega,\mathbb{R})$ adds an Ult$(V,U_h)$-generic filter $G'\subset j_h(\mathbb{P})$ extending $j``G$. We can then extend $j_h$ to an elementary embedding
\begin{center}
$j^*_h: V[G] \rightarrow \textrm{Ult}(V,U_h)[G']$
\end{center}
by defining $j^*_h(\tau_G) = j_h(\tau)_{G'}$. 

Now in $V[G][H]$, define an ideal $\mathcal{I}$ on $\powerset_{\omega_1}(\mathbb{R})$ by
\begin{center}
$S\in\mathcal{I}  \iff  \emptyset  \Vdash_{Coll(\omega,\mathbb{R})} \check{\mathbb{R}} \notin j^*_h(\check{S})$.
\end{center}
So $\mathbb{P}_{\mathcal{I}}$ is isomorphic to the subalgebra $\mathcal{B} = \{|| \check{R}\in j^*_h(\check{S}) || : S\subseteq \powerset_{\omega_1}(\mathbb{R})\}$ of the regular-open algebra RO$(Coll(\omega,\mathbb{R}))$.

$\mathcal{I}$ is fine: for any $x\in\mathbb{R}$, the set $T_x = \{\sigma: x\notin \sigma\}\in \mathcal{I}$ because clearly $\emptyset  \Vdash_{Coll(\omega,\mathbb{R})} \check{\mathbb{R}} \notin j^*_h(\check{T_x})$. $\mathcal{I}$ is normal: suppose $(S_x: x\in\mathbb{R})$ is a family of subsets of $\powerset_{\omega_1}(\mathbb{R})$ and $S$ is the diagonal union, i.e. $\sigma\in S$ if and only if there is some $x\in \sigma$ such that $\sigma\in S_x$. Then
\begin{center}
$||\check{\mathbb{R}}\in j^*_h(S)|| = ||\exists x\in\check{\mathbb{R}} \ (\check{\mathbb{R}}\in j^*_h(S_x))|| = \textrm{sup}_x ||\check{\mathbb{R}}\in j^*_h(S_x)||$.
\end{center}
This verifies normality of $\mathcal{I}$ and also verifies $\mathcal{B}$ is a $\mathfrak{c}$-complete subalgebra of RO$(Coll(\omega,\mathbb{R}))$. Since in $V[G][H]$, RO$(Coll(\omega,\mathbb{R}))$ has size $\mathfrak{c}^+$, has the $\mathfrak{c}^+$-chain condition, and is $\mathfrak{c}$-dense, $\mathcal{B}$ is $\mathfrak{c}$-dense and is a complete subalgebra of RO$(Coll(\omega,\mathbb{R}))$.

We now show $\mathcal{I}$ has the covering property relative to $V$. In $V[G][H]$, suppose $F:S\rightarrow \textrm{Ord}$ where $S\in \mathcal{I}^+$. Note that $F\in V[G]$. Let $p\in Coll(\omega,\mathbb{R})$ force ``$\check{\mathbb{R}}\in j_h^*(S)$" and $q\leq p$ force ``$j^*_h(F)(\check{\mathbb{R}})=\alpha$" for some ordinal $\alpha$. In $V$, let $F_0: \powerset_{\omega_1}(\mathbb{R})\rightarrow \textrm{Ord}$ such that $[F_0]_\mu = \alpha$. By the discussion above, before the proof of the theorem, 
\begin{center}
$\emptyset \Vdash_{Coll(\omega,\mathbb{R})} [F_0]_\mu = j_h(F_0)(\check{\mathbb{R}}) = j^*_h(F_0)(\check{\mathbb{R}})$.
\end{center}
Therefore,
\begin{center}
$q \Vdash_{Coll(\omega,\mathbb{R})} j^*_h(F_0)(\check{\mathbb{R}}) = j^*_h(F)(\check{\mathbb{R}})$.
\end{center}
This means the set $\{\sigma\in S: F(\sigma) = F_0(\sigma)\}$ is $\mathcal{I}$-positive.
\end{proof}

Now, let $\mathbb{P}$ be such that $\sf{CH}$ holds in $V[G][H]$. For example, we can take $\mathbb{P} = Coll(\omega_1,\mathbb{R})$. So in $V[G][H]$, $\mathfrak{c} = \omega_1$ and $\Theta^V=\omega_2$. By Theorem \ref{thm:dense_ord_covering}, in $V[G][H]$, there is an $\omega_1$-dense ideal $\mathcal{I}$ on $\powerset_{\omega_1}(\mathbb{R})$ that has the covering property with respect to $V$. Since $|\powerset_{\omega_1}(\mathbb{R})| = \omega_1$ in $V[G][H]$, we easily obtain an $\omega_1$-dense ideal on $\omega_1$ with the ordinal covering property. This and Lemma \ref{lem:agreement} give us one direction of Theorem \ref{thm:HI}.

\begin{corollary}\label{cor:DI}
Con$(\sf{ZF} + \sf{AD}_\mathbb{R} + {}$``$\Theta$ is regular'') implies Con($\sf{ZFC} +\sf{DI}$).
\end{corollary}

\begin{remark}\label{rem:additional}
We note that the $\omega_1$-dense ideal constructed above has the covering property with respect to $V$, so in fact, it satisfies the strengthening of $\sf{DI}$ in Footnote \ref{lab:amenable}, by Lemma \ref{lem:agreement}.
\end{remark}

\section{PRELIMINARIES}\label{sec:prelim}

This section, consisting of several subsections, develops some terminology and framework for the core model induction. The first subsection gives a brief summary of the theory of $\F$-premice and strategy premice developed in \cite{trang2013}. For a full development of these concepts, the reader should consult \cite{trang2013}. These concepts and notations will be used in the next subsection, which defines core model induction operators, which are the operators that we construct during the course of the core model induction in this paper. The next two sections briefly summarize the theory of hod mice and the HOD analysis in $\sf{AD}^+$ models (see \cite{ATHM} for a more detailed discussions of these topics). The reader who wishes to see the main argument can skip them on the first read, and go back when needed. Section \ref{sec:pullback} proves several important properties for reasonable hod pairs, defined in \ref{def:reasonable}, that we need for the proof of Theorem \ref{thm:DI}. The key result of this section is Lemma \ref{lem:pullback}, whose proof uses substantially Lemmata \ref{claim:bcshc}, \ref{lem:omega2inaccessible}. Lemma \ref{claim:bcshc} appears to be a new fact in the theory of hod mice at the level of ``$\sf{AD}_\mathbb{R} \ + $ $\Theta$ is regular". The last section reviews the technique of boolean valued comparisons for such hod pairs. Throughout this paper, we will identify a set $A\subset HC$ with $Code[A]\subset \mathbb{R}$, where $Code$ is a simple coding of elements of $HC$ by reals.

\subsection{$\F$-premice and strategy premice}\label{sec:strat_premice}

\begin{definition}
\label{model}
Let $\mathcal{L}_0$ be the language of set theory expanded by unary predicate
symbols $\dot{E}, \dot{B}, \dot{S}$, and constant symbols $\dot{a}$, 
$\dot{\param}$. Let 
$\Ll_0^-=\Ll_0\cut\{\dot{E},\dot{B}\}$.

Let $a$ be transitive. Let $\varrho:a\to\rank(a)$ be the rank function. 
We write $\ahat=\trancl(\{(a,\varrho)\})$. Let $\param\in\J_1(\ahat)$.

A \textbf{$\J$-structure over $a$ (with parameter $\param$) (for $\Ll_0$)} 
is a structure $\M$ for $\Ll_0$ such that $a^\M=a$, ($\param^\M=\param$), and 
there is
$\lambda\in[1,\Ord)$ such that $|\M|=\J_\lambda^{S^\M}(\ahat)$.

Here we also let $\l(\M)$ denote 
$\lambda$, the \textbf{length} of $\M$, and let $\ahat^\M$ denote $\ahat$.

For $\alpha\in[1,\lambda]$ let $\M_\alpha=\J_\alpha^{S^\M}(\ahat)$.
We say that $\M$ is \textbf{acceptable} iff for each 
$\alpha<\lambda$ and $\tau<\OR(\M_\alpha)$, if 
\[ 
\pow(\tau^{<\om}\cross\ahat^{<\om})\inter\M_\alpha\neq\pow(\tau^{<\om}
\cross\ahat^{<\om})\inter\M_{\alpha+1},\]
then there is a surjection 
$\tau^{<\om}\cross\ahat^{<\om}\to\M_\alpha$ in $\M_{\alpha+1}$.

A \textbf{$\J$-structure (for $\Ll_0$)} is a $\J$-structure over $a$, for some 
$a$.
\end{definition}

As all $\J$-structures we consider will be for $\Ll_0$, we will omit the 
phrase ``for $\Ll_0$''. We also often omit the phrase ``with parameter 
$\param$''. Note that if $\M$ is a $\J$-structure over $a$ then 
$|\M|$ is transitive and rud-closed, $\hat{a}\in M$, and $\OR\inter 
M=\rank(M)$. This last point is because we construct from $\hat{a}$ instead of 
$a$.

$\Fop$-premice will be $\J$-structures of the following form.

\begin{definition}\label{dfn:model}
A \textbf{$\J$-model over $a$ (with parameter $\param$)} 
is an acceptable
$\J$-structure over $a$ (with parameter $\param$), of the form
\[ \M = (M; E, B, S, a, \param) \]
where $\dot{E}^\mathcal{M}=E$, etc., and letting $\lambda=\l(\M)$, the following 
hold.
\begin{enumerate}
 \item $\M$ is amenable.
 \item $S = \langle S_\xi \ | \ \xi\in[1,\lambda) \rangle$
is a sequence of $\J$-models over $a$ (with parameter $\param$).
\item For each $\xi\in[1,\lambda)$, $\dot{S}^{S_\xi} = S\rest \xi$ and $\M_\xi = |S_\xi|$.
\item\label{item:extender} Suppose $E\neq\emptyset$. Then $B=\emptyset$ and 
there is an extender $F$ over $\M$ which is
($\hat{a}\times\gamma$)-complete for all $\gamma < \crit(F)$ and such that 
the premouse axioms \cite[Definition 2.2.1]{wilson2012contributions} hold for $(\M,F)$, and
$E$ codes $\tilde{F}\un\{G\}$ where: (i) $\tilde{F}\sub M$ is the 
amenable code for $F$ (as in \cite{steel2010outline}); and (ii) if $F$ is not 
type 2 then $G=\emptyset$, and otherwise $G$ is the ``longest'' non-type Z 
proper segment of $F$ in $\M$.\footnote{We use $G$ explicitly, instead of the 
code $\gamma^\M$ used for $G$ in \cite[Section 2]{FSIT}, because $G$ does not depend 
on which (if there is any) wellorder of $\M$ we use. This ensures that certain 
pure mouse operators are forgetful.}
\qedhere
\end{enumerate}
\end{definition}

Our notion of a ``$\J$-model over $a$" is a bit different from the notion of ``model with parameter $a$" in \cite{CMI} or \cite[Definition 2.1.1]{wilson2012contributions} in that we build into our notion some fine structure and we do not have the predicate $l$ used in \cite[Definition 2.1.1]{wilson2012contributions}. Note that with notation as above, if $\lambda$ is a successor ordinal then 
$M=J(S^\M_{\lambda-1})$, and otherwise, 
$M=\bigcup_{\alpha<\lambda}|S_\alpha|$.
The predicate $\dot{B}$ will be used to code extra information such as a (partial) branch of a tree in $M$.

\begin{definition}
\label{SomeNotations}
Let $\M$ be a $\J$-model over $a$ (with parameter $\param$). Let $E^\M$ denote $\dot{E}^\M$, 
etc. Let $\lambda=\l(\M)$, $S^\M_0=a$, $S^\M_\lambda=\M$, and
$\M|\xi=S^\M_\xi$ for all $\xi\leq\lambda$.
An \textbf{(initial) segment} of $\M$ is just a structure of 
the form $\M|\xi$ for some $\xi\in[1,\lambda]$. We write $\P\ins\M$ iff $\P$ is 
a segment of $\M$, and $\P\pins\M$ iff $\P\ins\M$ and $\P\neq\M$.
Let $\M||\xi$ be the structure having the same universe
and predicates as $\M|\xi$, except that $E^{\M||\xi}=\emptyset$.
We say that $\M$ is \textbf{$E$-active} iff
$E^\M\neq\emptyset$, and \textbf{$B$-active} iff
$B^\M\neq\emptyset$. \textbf{Active} means either
$E$-active or $B$-active; \textbf{$E$-passive} means not
$E$-active; \textbf{$B$-passive} means not $B$-active; and \textbf{passive}
means not active.

Given a $\J$-model $\M_1$ over $b$ and a $\J$-model $\M_2$ over
$\M_1$, we write $\M_2\downarrow b$ for the $\J$-model $\M$ over $b$,
such that $\M$ is ``$\M_1\conc\M_2$''. That is,
$|\M|=|\M_2|$, $a^\M=b$, $E^\M=E^{\M_2}$, $B^\M=B^{\M_2}$,
and $\P\pins\M$ iff $\P\ins\M_1$ or there is $\Q\pins\M_2$ such that 
$\P=\Q\downarrow b$, when such an $\M$ exists. Existence depends on whether the $\J$-structure $\M$ is acceptable.
\end{definition}

In the following, the variable $i$ should be interpreted as 
follows. When $i=0$, we ignore history, and so $\P$ is treated as a coarse 
object when determining $\Fop(0,\P)$. When $i=1$ we respect the history 
(given it exists).

\begin{definition}\label{operatorWithParama}
An \textbf{operator $\Fop$ with domain $D$} is a function with domain 
$D$, such that for some cone $C=C_\Fop$, possibly self-wellordered (sword),\footnote{$C$ is a cone if there are a cardinal $\kappa$ and a transitive set $a\in H_\kappa$ such that $C$ is the set of $b\in H_\kappa$ such that $a\in L_1(b)$; $a$ is called the base of the cone. A set $a$ is self-wellordered if there is a well-ordering of $a$ in $L_1(a)$. A set $C$ is a self-wellordered cone if $C$ is the restriction of a cone $C'$ to its own self-wellordered elements.} $D$ is the set of 
pairs $(i,X)$ such that either:
\begin{itemize}
 \item $i=0$ and $X\in C$, or
 \item $i=1$ and $X$ is a $\J$-model over $X_1\in C$,
\end{itemize}
and for each $(i,X)\in D$, $\Fop(i,X)$ is a $\J$-model over $X$ 
such that for each $\P\ins \Fop(i,X)$, $\P$ is fully 
sound. (Note that $\P$ is a $\J$-model over $X$, so soundness is in this 
sense.)

Let $\Fop,D$ be as above. We say $\Fop$ is \textbf{forgetful} iff 
$\Fop(0,X)=\Fop(1,X)$ whenever $(0,X),(1,X)\in D$, and whenever $X$ is a 
$\J$-model over $X_1$, and $X_1$ is a $\J$-model over $X_2\in C$, we have 
$\Fop(1,X)=\Fop(1,X\downarrow X_2)$. 
Otherwise we say $\Fop$ is \textbf{historical}. Even when $\Fop$ is historical, 
we often just write $\Fop(X)$ instead of $\Fop(i,X)$ when the nature of $\Fop$ is clear from the context.
We say $\Fop$ is \textbf{basic} iff for all 
$(i,X)\in D$ and $\P\ins\Fop(i,X)$, we have $E^\P=\emptyset$.
We say $\Fop$ is \textbf{projecting} iff for all $(i,X)\in D$, we have 
$\rho_\om^{\Fop(i,X)}=X$.
\end{definition}

Here are some illustrations. Strategy 
operators (to be explained in more detail later) are basic, and as usually 
defined, projecting and historical. Suppose we have an iteration strategy $\Sigma$ and we want to build a $\J$-model $\N$ (over some $a$) that codes a fragment of $\Sigma$ via its predicate $\dot{B}$. We feed $\Sigma$ into $\N$ by always providing $b = \Sigma(\T)$, for the $<$-$\N$-least tree $\T$ for which this information is required. So given a reasonably closed level $\P\lhd \N$, the choice of which tree $\T$ should be processed next will usually depend on the information regarding $\Sigma$ already encoded in $\P$ (its history). Using an operator $\Fop$ to build $\N$, then $\Fop(i,\P)$ will be a structure extending $\P$ and over which $b = \Sigma(\T)$ is encoded. The variable $i$ should be interpreted as follows. When $i = 1$, we respect the history of $\P$ when selecting $\T$. When $i = 0$ we ignore history when selecting $\T$ . The operator $\Fop(X)=X^\#$ is forgetful 
and projecting, and not basic; here $\Fop(X) = \Fop(0,X)$.

\begin{definition}
For any $P$ and any ordinal $\alpha\geq 1$, the operator $\Jop_\alpha(\ \cdot\ ;P)$ 
is defined as follows.\footnote{The ``$\mathsf{m}$'' is 
for 
``model''.} For $X$ such that $P\in\J_1(\hat{X})$, let
$\Jop_{\alpha}(X;P)$ be the $\J$-model $\M$ over $X$, with parameter $P$, such 
that
$|\M|=\J_\alpha(\hat{X})$ and for each $\beta\in[1,\alpha]$, $\M|\beta$ is 
passive. Clearly $\Jop_\alpha(\ \cdot\ ;P)$ is basic and forgetful.
If $P=\emptyset$ or we wish to supress $P$, we just 
write $\Jop_\alpha(\ \cdot\ )$.

\end{definition}

\begin{definition}[Potential $\Fop$-premouse, 
$\C_\Fop$]\label{potentialJpremouse}
Let $\Fop$ be an operator with domain $D$ of self-wellordered sets. Let $b\in C_\Fop$, so there is a well-ordering of $b$ in $L_1[b]$. A \textbf{potential
$\Fop$-premouse over $b$} is an acceptable $\J$-model $\M$ over $b$
such that there is an ordinal $\iota>0$ and an increasing, closed sequence
$\langle \zeta_\alpha \rangle_{\alpha\leq \iota}$ of ordinals such that for each
$\alpha\leq \iota$, we have:
\begin{enumerate}
\item $0=\zeta_0\leq\zeta_\alpha\leq\zeta_\iota=l(\M)$ (so $\M|\zeta_0=b$ and
$\M|\zeta_\iota=\M$).
\item If $1<\iota$ then $\M|\zeta_1=\Fop(0,b)$.
\item If $1=\iota$ then $\M\ins\Fop(0,b)$.
\item If $1<\alpha+1<\iota$ then
$\M|\zeta_{\alpha+1}=\Fop(1,\M|\zeta_\alpha)\downarrow b$.
\item If $1<\alpha+1=\iota$, then $\M\ins
\Fop(1,\M|\zeta_\alpha)\downarrow b$.
\item Suppose $\alpha$ is a limit. Then $\M|\zeta_\alpha$ is $B$-passive, and 
if $E$-active, then $\crit(E^{\M|\zeta_\alpha})>\rank(b)$.
\end{enumerate}
We say that $\M$ is \textbf{($\Fop$-)whole} iff $\iota$ is a limit or else,
$\iota=\alpha+1$ and $\M = \Fop(\M|\zeta_\alpha)\downarrow b$.

A \textbf{(potential) $\Fop$-premouse} is a (potential) $\Fop$-premouse over 
$b$, for some $b$.
\end{definition}
\begin{definition}
Let $\Fop$ be an operator and $b\in C_\Fop$. Let $\N$ be a whole
$\Fop$-premouse over $b$. A \textbf{potential continuing $\Fop$-premouse over 
$\N$} is a $\J$-model $\M$ over $\N$ such that $\M\downarrow b$ is a potential 
$\Fop$-premouse over $b$. (Therefore $\N$ is a whole strong cutpoint of $\M$.)

We say that $\M$ (as above) is \textbf{whole} iff $\M\downarrow 
b$ is whole.

A \textbf{(potential) continuing $\Fop$-premouse} is a (potential) continuing 
$\Fop$-premouse over $b$, for some $b$.
\end{definition}

\begin{definition}
\label{Lps}
$\Lp^\Fop(a)$ for an operator $\Fop$ denotes the stack of all 
countably $\Fop$-iterable $\Fop$-premice $\M$ over $a$ such that $\M$ is fully sound 
and 
projects to $a$.\footnote{Countable substructures of $\M$ are $(\omega,\omega_1+1)$-$\Fop$-iterable, i.e. all iterates are $\Fop$-premice. See \cite[Section 2]{trang2013} for more details on $\Fop$-iterability.}

Let $\N$ be a whole $\Fop$-premouse over $b$, for $b\in C_\Fop$. Then $\Lp^\Fop_+(\N)$ denotes the stack 
of all countably $\Fop$-iterable (above $o(\N)$) continuing $\Fop$-premice $\M$ over $\N$ such 
that $\M\downarrow b$ is fully sound and projects to $\N$.\footnote{Often times in this paper, when the context is clear, we will use the notation Lp for Lp$_+$.}

We say that $\Fop$ is
\textbf{uniformly $\Sigma_1$} iff there are $\Sigma_1$ formulas $\varphi_1$ and 
$\varphi_2$ in $\Ll_0^-$ such that whenever $\M$ is a (continuing) 
$\Fop$-premouse, 
then
the set of whole proper segments of $\M$ is defined over $\M$ by $\varphi_1$ 
($\varphi_2$). For 
such an operator $\Fop$, let $\varphi^\Fop_\wh$ denote the least such 
$\varphi_1$.
\end{definition}
\begin{definition}[Mouse operator]\label{MouseOperator}
Let $Y$ be a projecting, uniformly $\Sigma_1$ operator.
A \textbf{$Y$-mouse operator $\Fop$ with domain $D$}
is an operator with domain $D$ such for each $(0,X)\in D$, 
$\Fop(0,X)\pins\Lp^Y(X)$, and for each $(1,X)\in 
D$, $\Fop(1,X)\pins\Lp^Y_+(X)$.\footnote{This restricts the usual notion defined in \cite{CMI}.} (So any $Y$-mouse 
operator is an operator.) A $Y$-mouse operator $\Fop$ is called \textbf{first-order} if there are formulas $\varphi_1$ and $\varphi_2$ in the language of $Y$-premice such that $\Fop(0,X)$ ($\Fop(1,X)$) is the first $\M\lhd \textrm{Lp}^Y(X)$ ($\textrm{Lp}^Y_+(X)$) satisfying $\varphi_1$ ($\varphi_2$).

A \textbf{mouse operator} is a $\Jop_1$-mouse operator.
\end{definition}

We can then define $\Fop$-solidity, the $L^{\Fop}[\es]$-construction etc. as usual (see \cite{trang2013} for more details). We now define the kind of condensation that mouse operators need to satisfy to ensure for example that the $L^{\Fop}[\es]$-construction converges. We define the coarse version of condensation (condense coarsely) here for illustrative purposes. The finer version (condense finely), which is more technical, is discussed in detail in \cite{trang2013}. The core model induction operators, which form a subclass of the $Y$-mouse operators, will have these condensation properties.

\begin{definition}\label{dfn:condenses_coarsely}
 Let $Y$ be an operator. We say that $Y$ 
\textbf{condenses coarsely} iff 
for all $i\in\{0,1\}$ and $(i,\bar{X}),(i,X)\in\dom(Y)$, and all $\J$-models 
$\M^+$ over $\bar{X}$, if $\pi:\M^+\to Y_i(X)$ is fully elementary and fixes the parameters in the definition of $Y$, then 
\begin{enumerate}
 \item if $i=0$ then $\M^+\ins Y_0(\bar{X})$; and
 \item if $i=1$ and $X$ is a sound whole $Y$-premouse, then $\M^+\ins 
Y_1(\bar{X})$.\qedhere
\end{enumerate}
\end{definition}

We now proceed to defining $\Sigma$-premice, for an iteration strategy
$\Sigma$. We first define the operator to be used to feed in $\Sigma$.

\begin{definition}[$\BBB(a,\Tt,b)$, $b^\N$]\label{dfn:B(M,T,b)}
Let $a,\P$ be transitive, with $\P\in\J_1(\ahat)$. Let
$\lambda>0$ and let
$\Tt$ be an iteration tree\footnote{We formally take an \emph{iteration tree} to include
the entire sequence $\left<M^\Tt_\alpha\right>_{\alpha<\lh(\Tt)}$ of models. So it is 
$\Sigma_0(\Tt,\param)$ to assert that ``$\Tt$ is an iteration tree on
$\param$''.} on $\P$, of length $\om\lambda$, with $\Tt\rest\beta\in a$ for 
all $\beta\leq\om\lambda$. Let $b\sub\om\lambda$. We
define
$\N=\BBB(a,\Tt,b)$ recursively on $\lh(\Tt)$, as the $\J$-model $\N$ over
$a$ with parameter $\P$\footnote{$\P=M^\Tt_0$ is determined by $\Tt$.} such that:
\begin{enumerate}
 \item $\l(\N)=\lambda$,
 \item for each $\gamma\in(0,\lambda)$,
$\N|\gamma=\BBB(a,\Tt\rest\om\gamma,[0,\om\gamma]_\Tt)$,
 \item $B^\N$ is the set of ordinals $\OR(a)+\gamma$ such that $\gamma\in
b$,
\item $E^\N=\emptyset$.
\end{enumerate}
We also write $b^\N=b$.
\end{definition}

It is easy to see that every initial segment of $\N$ is sound, so $\N$ is 
acceptable and is indeed a $\J$-model (not just a $\J$-structure). 

In the context of a $\Sigma$-premouse $\M$ for an iteration strategy $\Sigma$,
if $\Tt$ is the $<_\M$-least tree for which $\M$ lacks instruction regarding
$\Sigma(\Tt)$,
then $\M$ will already have been instructed regarding $\Sigma(\Tt\rest\alpha)$
for all
$\alpha<\lh(\Tt)$. Therefore if $\lh(\Tt)>\om$ then
$\BBB(\M,\Tt,\Sigma(\Tt))$ codes redundant information (the branches already in
$\Tt$) before coding $\Sigma(\Tt)$. This redundancy seems to allow one to
prove slightly stronger condensation properties, given that $\Sigma$ has nice
condensation properties (see \cite{trang2013}). It also
simplifies the definition.

\begin{definition}
 Let $\Sigma$ be a partial iteration strategy. Let
$C$ be a class of iteration trees, closed under
initial segment. We say that $(\Sigma,C)$ is \textbf{suitably condensing} iff
for every $\Tt\in C$ such that $\Tt$ is via $\Sigma$ and $\lh(\Tt)=\lambda+1$
for some limit $\lambda$, either (i) $\Sigma$ has hull condensation with respect
to $\Tt$, or (ii) $b^\Tt$ does not drop and $\Sigma$ has branch
condensation with respect to $\Tt$, that is, any hull $\U^\smallfrown c$ of $\T^\smallfrown b$ is according to $\Sigma$.
\end{definition}

When $C$ is the class of all iteration trees according to $\Sigma$, we simply omit it from our notation. 

\begin{definition}\label{dfn:J_premice}
 Let $\varphi$ be an $\Ll_0$-formula. Let 
$\P$
be transitive. Let $\M$ be a $\J$-model (over some $a$), with parameter $\P$. 
Let $\Tt\in\M$. We
say that $\varphi$ \textbf{selects $\Tt$ for $\M$}, and write 
$\Tt=\Tt^\M_\varphi$, iff
\begin{enumerate}[(a)]
 \item $\Tt$ is the unique 
$x\in\M$ such that $\M\sats\varphi(x)$,
 \item $\Tt$ is an iteration tree on $\P$ of
limit length,
 \item for every $\N\pins\M$, we have $\N\not\sats\varphi(\Tt)$, and
 \item for every limit $\lambda<\lh(\Tt)$, there is $\N\pins\M$ such that
$\N\sats\varphi(\Tt\rest\lambda)$.\qedhere
\end{enumerate}
\end{definition}

One instance of $\phi(\P,\T)$ is, in the case $a$ is self-wellordered, the formula ``$\T$ is the least tree on $\P$ that doesn't have a cofinal branch", where least is computed with respect to the canonical well-order of the model.

\begin{definition}[Potential $\P$-strategy-premouse, $\Sigma^\M$]
\label{PotStrPremouse}
Let $\varphi\in\Ll_0$. Let $\P,a$ be transitive with $\P\in\J_1(\ahat)$. A
\textbf{potential $\P$-strategy-premouse (over $a$, of type
$\varphi$)} is a $\J$-model $\M$ over $a$, with parameter $\P$, 
such that
the $\BBB$ operator is used to 
feed in
an iteration strategy for trees on $\P$, using the sequence of trees
naturally determined by $S^\M$ and selection by $\varphi$. We let $\Sigma^\M$ 
denote the
partial strategy coded by the predicates $B^{\M|\eta}$, for
$\eta\leq\l(\M)$.

In more detail, there is an
increasing, closed sequence of ordinals
$\left<\eta_\alpha\right>_{\alpha\leq\iota}$ with the following properties.
We will also define $\Sigma^{\M|\eta}$ for all $\eta\in[1,\l(\M)]$ and
$\Tt_\eta=\Tt^\M_\eta$ for
all $\eta\in[1,\l(\M))$.\begin{enumerate}
 \item $1=\eta_0$ and 
$\M|1=\Jop_1(a;\P)$ and 
$\Sigma^{\M|1}=\emptyset$.
 \item $\l(\M)=\eta_\iota$, so $\M|\eta_\iota=\M$.
 \item Given $\eta\leq\l(\M)$ such that $B^{\M|\eta}=\emptyset$, we set
$\Sigma^{\M|\eta}=\bigcup_{\eta'<\eta}\Sigma^{\M|\eta'}$.
\end{enumerate}

Let $\eta\in[1,\l(\M)]$. Suppose
there is $\gamma\in[1,\eta]$ and $\Tt\in\M|\gamma$ such that
$\Tt=\Tt^{\M|\gamma}_\varphi$, and $\Tt$ is via 
$\Sigma^{\M|\eta}$, but no
proper extension of $\Tt$ is via
$\Sigma^{\M|\eta}$. Taking $\gamma$ minimal such, let
$\Tt_\eta=\Tt^{\M|\gamma}_\varphi$.
Otherwise let $\Tt_\eta=\emptyset$.
\begin{enumerate}\setcounter{enumi}{3}
 \item Let $\alpha+1\leq\iota$. Suppose $\Tt_{\eta_\alpha}=\emptyset$.
Then $\eta_{\alpha+1}=\eta_\alpha+1$ and
$\M|\eta_{\alpha+1}=\Jop_1(\M|\eta_\alpha;\P)\downarrow a$.
\item Let $\alpha+1\leq\iota$.
Suppose $\Tt=\Tt_{\eta_\alpha}\neq\emptyset$.
Let $\om\lambda=\lh(\Tt)$.
Then for some $b\sub\om\lambda$, and
$\mathcal{S}=\BBB(\M|\eta_\alpha,\Tt,b)$, we have:
\begin{enumerate}
\item $\M|\eta_{\alpha+1}\ins\mathcal{S}$.
\item If $\alpha+1<\iota$ then $\M|\eta_{\alpha+1}=\mathcal{S}$.
\item If $\mathcal{S}\ins\M$ then $b$ is a $\Tt$-cofinal branch.\footnote{We allow
$\M^\Tt_b$ to be illfounded, but then
$\Tt\conc b$ is not an iteration tree, so is not continued by $\Sigma^\M$.}
\item For $\eta\in[\eta_\alpha,\l(\M)]$ such that
$\eta<\l(\mathcal{S})$,
$\Sigma^{\M|\eta}=\Sigma^{\M|\eta_\alpha}$.
\item If $\mathcal{S}\ins\M$
then $\Sigma^{\mathcal{S}}=\Sigma^{\M|\eta_\alpha}\un\{(\Tt,b^\mathcal{S})\}$.
\end{enumerate}
\item For each limit $\alpha\leq\iota$, $B^{\M|\eta_\alpha}=\emptyset$.\qedhere
\end{enumerate}
\end{definition}

\begin{definition}[Whole]\label{dfn:whole_strategy}
Let $\M$ be a potential $\P$-strategy-premouse of type $\varphi$. We say $\P$
is \textbf{$\varphi$-whole} (or just \textbf{whole} if $\varphi$ is fixed) iff
for every $\eta<\l(\M)$, if
$\Tt_{\eta}\neq\emptyset$ and $\Tt_\eta\neq\Tt_{\eta'}$ for all
$\eta'<\eta$, then for some $b$,
$\BBB(\M|\eta,\Tt_{\eta},b)\ins\M$.\footnote{\emph{$\varphi$-whole} depends on
$\varphi$ as the definition of $\Tt_\eta$ does.}
\end{definition}

\begin{definition}[Potential $\Sigma$-premouse]
Let $\Sigma$ be a (partial) iteration strategy for a transitive
structure $\P$.
A \textbf{potential $\Sigma$-premouse (over $a$, of type $\varphi$)}
is a
potential
$\P$-strategy premouse $\M$ (over $a$, of type $\varphi$) such that
$\Sigma^\M\sub\Sigma$.\footnote{If $\M$ is a model all of whose proper segments
are potential $\Sigma$-premice, and the rules for potential $\P$-strategy
premice
require that $B^\M$ code a $\Tt$-cofinal branch, but $\Sigma(\Tt)$ is not
defined,
then $\M$ is not a potential $\Sigma$-premouse, whatever its predicates are.}
\end{definition}

\begin{definition}\label{dfn:strategy_op}
Let $\P$ be transitive and $\Sigma$ a partial iteration strategy for $\P$.
Let $\varphi\in\Ll_0$. Let $\Fop=\Fop_{\Sigma,\varphi}$ be 
the operator such that:
\begin{enumerate}
 \item $\Fop_0(a)=\Jop_1(a;\P)$, for all transitive 
$a$ such that $\P\in\J_1(\ahat)$;
 \item Let $\M$ be a sound branch-whole $\Sigma$-premouse of 
type $\varphi$. Let $\lambda=\l(\M)$ and with notation as in \ref{PotStrPremouse}, let 
$\Tt=\Tt_\lambda$. If $\Tt=\emptyset$ then $\Fop_1(\M)=\Jop_1(\M;\P)$. If 
$\Tt\neq\emptyset$ then $\Fop_1(\M)=\BBB(\M,\Tt,b)$ where $b=\Sigma(\Tt)$.
\end{enumerate}
We say that $\Fop$ is a \textbf{strategy operator}.
\end{definition}

\begin{lemma}\label{lem:FSigma_condenses}
Let $\P$ be countable and transitive. Let $\varphi$ be a formula of $\Ll_0$.
Let $\Sigma$ be a partial strategy for $\P$. Let $D_\varphi$ be the class of 
iteration trees $\Tt$ on $\P$ such that for some $\J$-model $\M$, with 
parameter $\P$, we have $\Tt=\Tt^\M_\varphi$. Suppose that $(\Sigma,D_\varphi)$ 
is suitably condensing. Then 
$\Fop_{\Sigma,\varphi}$ is uniformly $\Sigma_1$, projecting, and condenses finely.
\end{lemma}

\begin{definition}\label{dfn:MFsharp}
Let $a$ be transitive and let $\Fop$ be an operator. 
We say that 
\textbf{$\M_1^{\Fop,\#}(a)$ exists}
iff there is a $(0,|a|,|a|+1)$-$\Fop$-iterable, non-$1$-small 
$\Fop$-premouse over $a$. We write 
$\M_1^{\Fop,\#}(a)$ for the least such 
sound structure. For $\Sigma,\P,a,\varphi$ as in Definition \ref{dfn:strategy_op}, we 
write 
$\M_1^{\Sigma,\varphi,\#}(a)$ for $\M_1^{\Fop_{\Sigma,\varphi},\#}(a)$.

Let $\Ll_0^+$ be the language $\Ll_0\un\{\dot{\prec},\dot{\Sigma}\}$, 
where $\dot{\prec}$ is the binary relation defined by ``$\dot{a}$ is 
self-wellordered, with ordering $\prec_{\dot{a}}$, and $\dot{\prec}$ is the 
canonical wellorder of the universe extending $\prec_{\dot{a}}$'', and 
$\dot{\Sigma}$ is the partial function defined by ``$\dot{\param}$ is a transitive 
structure and the universe is a potential $\dot{\param}$-strategy premouse over 
$\dot{a}$ and $\dot{\Sigma}$ is the associated partial putative iteration 
strategy for $\dot{\param}$''. Let $\varphi_\textrm{all}(\Tt)$ be the 
$\Ll_0$-formula 
``$\Tt$ is the $\dot{\prec}$-least limit length iteration tree $\Uu$ on 
$\dot{\param}$ such that $\Uu$ is via $\dot{\Sigma}$, but no proper extension 
of $\Uu$ is via $\dot{\Sigma}$''. Then for $\Sigma,\P,a$ as in 
Definition \ref{dfn:strategy_op}, we sometimes write $\M_1^{\Sigma,\#}(a)$ for 
$\M_1^{\Fop_{\Sigma,\varphi_{\textrm{all}}},\#}(a)$.

Let $\kappa$ be a cardinal and suppose that $\MFsharp=\M_1^{\Fop,\#}(a)$ 
exists and is $(0,\kappa^++1)$-iterable.
We write $\Lambda_\MFsharp$ for the unique $(0,\kappa^++1)$-iteration strategy 
for $\MFsharp$ (given that $\kappa$ is fixed).
\end{definition}

\subsection{Core model induction operators}\label{sec:cmioperators}

In core model induction applications, we often have a pair $(\P,\Sigma)$
where $\P$ is a hod premouse and $\Sigma$ is $\P$'s strategy with branch
condensation and is fullness preserving (relative to mice with strategies in some pointclass) or
$\P$ is a sound (hybrid) premouse projecting to some countable set $a$ and
$\Sigma$ is the unique (normal) ($\omega_1+1$)-strategy for $\P$.  Let $\Fop$ be the operator corresponding to 
$\Sigma$ (using the formula $\varphi_{\rm{all}}$) and suppose $\M_1^{\Fop,\sharp}$ exists. Then \cite[Lemma 4.8]{trang2013} shows that $\Fop$
condenses finely and $\M_1^{\Fop,\sharp}$ generically interprets $\F$. Also, the core model induction will give us that $\Fop\rest \mathbb{R}$ is self-scaled (defined below). In the following, we will write $\M_1^{\Sigma,\sharp}$ for $\M_1^{\Fop,\sharp}$.

In this section, our main goal is to introduce the main concepts that one uses in the core model induction through the hierarchy Lp$^{^\gTheta\Sigma}(\mathbb{R}, \Sigma\rest \textrm{HC})$\footnote{An equivalent way to define this is to first fix a canonical coding function Code$:\textrm{HC}\rightarrow \mathbb{R}$ and consider Lp$^{^\gTheta\Sigma}(\mathbb{R},\textrm{Code}(\Sigma\rest\mathrm{HC}))$.} \footnote{Instead of feeding $\Sigma$ into the hierarchy, we feed in $\Lambda$, the canonical strategy of $\M_1^{\Sigma,\sharp}$, into the hierarchy. Roughly speaking, the trees according to $\Lambda$ that we feed into Lp$^{^\gTheta\Sigma}(\mathbb{R},\textrm{Code}(\Sigma\rest\mathrm{HC}))$ are those making the local HOD of Lp$^{^\gTheta\Sigma}(\mathbb{R},\textrm{Code}(\Sigma\rest\mathrm{HC}))|\alpha$ generically generic, for appropriately chosen ordinals $\alpha$. See \cite{trang2013}.}.  Here Lp$^{^\gTheta\Sigma}(\mathbb{R}, \Sigma\rest \textrm{HC})$ is the union of all sound, $\Theta$-$g$-organized $\Sigma$-premice $\M$ over $(\mathbb{R}, \Sigma\rest \textrm{HC})$ such that $\rho_\omega(\M) = \mathbb{R}$ and whenever $\pi:\M^*\rightarrow \M$ is sufficiently elementary and $\M^*$ is countable and transitive, then $\M^*$ has a unique ($\omega_1+1$)-$\Sigma$-iteration strategy $\Lambda$.\footnote{This means whenever $\T$ is an iteration tree according to $\Lambda$ with last model $\N$, then $\N$ is a $\Sigma$-premouse.} See \cite{trang2013} for a precise definition of $g$-organized $\Sigma$-premice, $\Theta$-$g$-organized $\Sigma$-premice, Lp$^{^\g\Sigma}(x)$, Lp$_+^{^\g\Sigma}(x)$ and other related concepts like operators. When we write Lp$^{^\g\Sigma}$ or Lp$_+^{^\g\Sigma}$, we refer to the hierarchy of $g$-organized $\Sigma$-mice; when we write Lp$^{^\gTheta\Sigma}$ or Lp$_+^{^\gTheta\Sigma}$, we refer to the hierarchy of $\Theta$-$g$-organized $\Sigma$-mice. The $g$-organized hierarchy of $\Sigma$-mice is considered (instead of the traditional ``least branch" hierarchy of $\Sigma$-mice) because the $S$-constructions (cf. \cite{schindler2009self}, where they are called $P$-constructions) work out nicely for this hierarchy.\footnote{It is not clear how one can perform $S$-constructions over the least branch hierarchy.} The $\Theta$-$g$-organized hierarchy, which is a slight modification of the $g$-organized hierarchy, is considered because the scales analysis under optimal hypotheses can be carried out in Lp$^{^\gTheta\Sigma}(\mathbb{R}, \Sigma\rest \textrm{HC})$ in much the same manner as the scales analysis in Lp$(\mathbb{R})$.\footnote{\cite{trang2013} generalizes Steel's scales analysis in \cite{Scalesendgap, K(R)} to Lp$^{^\gTheta\Sigma}(\mathbb{R}, \Sigma\rest \textrm{HC})$ for various classes of nice strategies $\Sigma$. It is not clear that one can carry out the full scales analysis for the hierarchy Lp$^{^\g\Sigma}(\mathbb{R}, \Sigma\rest \textrm{HC})$.} For the purpose of this paper, it will not be important to go into the detailed definitions of these hierarchies. Whenever it makes sense to define Lp$^\Sigma(x)$ and Lp$^{^\g\Sigma}(x)$, \cite{trang2013} shows that $\powerset(x)\cap \textrm{Lp}^\Sigma(x) = \powerset(x)\cap \textrm{Lp}^{^\g\Sigma}(x)$ (and similarly for Lp$^{^\gTheta\Sigma}(x)$); also in the case it is not clear how to make sense of Lp$^\Sigma(x)$ (say for instance when $x=\mathbb{R}$), it still makes sense to define $\textrm{Lp}^{^\g\Sigma}(x)$ and $\textrm{Lp}^{^\gTheta\Sigma}(x)$ and in that case, \cite{trang2013} shows that $\powerset(x)\cap \textrm{Lp}^{^\g\Sigma}(x) = \powerset(x)\cap \textrm{Lp}^{^\gTheta\Sigma}(x)$. In the paragraph below, we briefly remark on how the $S$-constructions work for the $g$-organized hierarchy and for the $\Theta$-$g$-hierarchy.

Suppose $\Fop$ is a nice operator (with parameter $\param$)\footnote{\label{footnote:nice}Nice is defined in \cite[Definition 3.8]{trang2013}. Roughly speaking, these are operators that condense well and determine themselves on generic extensions. CMI operators defined in this section are nice.}  and suppose $\M$ is a $\mathcal{G}$-mouse (over some transitive $a$), where $\mathcal{G}$ is either $^\g\Fop$ or $^\gTheta\Fop$. Suppose $\delta$ is a cutpoint of $\M$ and suppose $\N$ is a transitive structure such that $\delta\subseteq \N\subseteq \M|\delta$ and $\param\in \N$. Suppose $\mathbb{P}\in \J_\omega[\N]$ is such that $\M|\delta$ is $\mathbb{P}$-generic over $\J_\omega[\N]$ and suppose whenever $\Q$ is a $\mathcal{G}$-mouse over $\N$ such that $H^\Q_\delta = \N$ then $\M|\delta$ is $\mathbb{P}$-generic over $\Q$. Then the $S$-constructions (or $P$-constructions) from \cite{schindler2009self} give a $\mathcal{G}$-mouse $\R$ over $\N$ such that $\R[\M|\delta] = \M$. The $S$-constructions give the sequence $(\R_\alpha : \delta < \alpha\leq \lambda)$ of $\mathcal{G}$-premice over $\N$, where 
\begin{enumerate}[(i)]
\item $\R_{\delta+1} = \Jop_\omega(\N)$;
\item if $\alpha$ is limit then let $\R_\alpha^* = \bigcup_{\beta<\alpha}\R_\beta$. If $\M|\alpha$ is passive, then let $\R_\alpha = \R_\alpha^*$. So $\R_\alpha$ is passive. If $B^{\M|\alpha}\neq \emptyset$, then let $\R_\alpha = (|\R_\alpha^*|; \emptyset, B^{\M|\alpha}, \bigcup_{\beta<\alpha} S^{\R_\beta}, \N ,\param)$. Suppose $E^{\M|\alpha}\neq \emptyset$; let $E^* = E^{\M|\alpha}\cap |\R_\alpha^*|$, then we let $\R_\alpha=(|\R_\alpha^*|; E^*, \emptyset, \bigcup_{\beta<\alpha} S^{\R_\beta}, \N ,\param)$. By the hypothesis, we have $\R_\alpha[\M|\delta]=\M|\alpha$.
\item Suppose we have already constructed $\R_\alpha$ and (by the hypothesis) maintain that $\R_\alpha[\M|\delta]=\M|\alpha$. Then $\R_{\alpha+1} = \Jop_\omega(\R_\alpha)$.
\item $\lambda$ is such that $\R_\lambda[\M|\delta]=\M$. We set $\R_\lambda = \R$.
\end{enumerate} 	
We note that the full constructions from \cite{schindler2009self} do not require that $\delta$ is a cutpoint of $\M$ but we don't need the full power of the $S$-constructions in our paper. Also, the fact that $\M$ is g-organized (or $\Theta$-g-organized) is important for our constructions above because it allows us to get past levels $\M|\alpha$ for which $B^{\M|\alpha}\neq \emptyset$. Because of this fact, in this paper, hod mice are reorganized into the g-organized hierarchy, that is if $\P$ is a hod mouse then $\P(\alpha+1)$ is a g-organized $\Sigma_{\P(\alpha)}$-premouse for all $\alpha<\lambda^\P$. The $S$-constructions are also important in many other contexts. One such context is the local HOD analysis of levels of Lp$^{^\gTheta\Fop}(\mathbb{R},\Fop\rest\mathbb{R})$, which features in the scales analysis of Lp$^{^\gTheta\Fop}(\mathbb{R},\Fop\rest\mathbb{R})$ (cf. \cite{trang2013}).

In the following, a transitive structure $N$ is \emph{closed} under an operator $\Omega$ if whenever $x\in \dom(\Omega)\cap N$, then $\Omega(x)\in N$. We are now in a position to introduce the core model induction operators that we will need in this paper. These are particular kinds of mouse operators (in the sense of \cite[Example 3.41]{schlutzenberg2016fine}) that are constructed during the course of the core model induction. These operators can be shown to satisfy the sort of condensation described in \cite[Section 3]{schlutzenberg2016fine} (e.g. condense coarsely and  condense finely), relativize well,  and determine themselves on generic extensions. 

\begin{definition}[relativizes well]\label{relativizeWell}
Let $\Omega$ be an a $Y$-mouse operator for some operator $Y$.\footnote{$Y$ may be the rud operator, in which case $\Omega$ is just a mouse operator in the usual sense.} We say that $\Omega$ \emph{relativizes well} if there is a formula $\phi(x,y,z)$ such that for any $a,b\in \mathrm{dom}(\Omega)$ such that $a\in L_1(b)$, whenever $N$ is a transitive model of $\mathsf{ZFC}^-$ such that $N$ is closed under $Y$ and $a,b, \Omega(b)\in N$, then $\Omega(a)\in N$ and is the unique $x\in N$ such that $N\vDash \phi[x, a, \Omega(b)]$.
\end{definition}
\begin{definition}[determines itself on generic extensions]\label{detGenExts}
Suppose $\Omega$ is an operator. We say that $\Omega$ \emph{determines itself on generic extensions} if there is a formula $\phi(x,y,z)$ and a parameter $c\in HC$ such that for any countable transitive structure $N$ of $\mathsf{ZFC}^-$ such that $N$ contains $c$ and is closed under $\Omega$, for any generic extension $N[g]$ of $N$ in $V$, $\Omega\cap N[g]\in N[g]$ and is definable over $N[g]$ via $(\phi,c)$, i.e. for any $e\in N[g]\cap \dom(\Omega)$, $\Omega(e)=d$ if and only if $d$ is the unique $d'\in N[g]$ such that $N[g]\vDash \phi[c,d',e]$.
\end{definition}


\begin{definition}\label{dfn:C_Gamma}
Let $\Gamma$ be an inductive-like pointclass. For $x\in\mathbb{R}$, $C_\Gamma(x)$ denotes the set 
of all $y\in\mathbb{R}$ such that for some ordinal $\gamma<\omega_1$, $y$ (as a subset of $\omega$) is 
$\Delta_\Gamma(\{\gamma,x\})$.

Let $x\in\HC$ be transitive and let $f:\omega\to x$ be a surjection. Then $c_f\in\mathbb{R}$ 
denotes the code for $(x,\in)$ determined by $f$.
And $C_\Gamma(x)$ denotes the set of all 
$y\in\HC\cap\powerset(x)$ such that for all surjections $f:\omega\to x$ we have $f^{-1}(y)\in 
C_\Gamma(c_f)$.
\end{definition}

We say that $\vec{A}$ is a self-justifying-system (sjs) if for any $A\in \textrm{rng}(\vec{A})$, $\neg A\in \textrm{rng}(\vec{A})$ and there is a scale $\varphi$ on $A$ such that the set of prewellorderings associated with $\varphi$ is a subset of $\textrm{rng}(\vec{A})$. A set $Y\subseteq \mathbb{R}$ is \textit{self-scaled} if there are scales on $Y$ and $\mathbb{R}\backslash Y$ which are projective in $Y$.

In the following, $\eta$ is a strong cutpoint of $\N$ if there is no extender $E$ on the sequence of $\N$ such that crt$(E) \leq \eta \leq \textrm{lh}(E)$.

\begin{definition}\label{dfn:k-suitable}
  Let $(\Omega,A)$ be as above and let $t\in\HC$ with $\MFsharp\in\J_1(t)$. Let $1\leq k<\om$. A 
premouse $\N$ over $t$ is
\emph{$\Omega$-$\Gamma$-$k$-suitable} (or just \emph{$k$-suitable} if $\Gamma$ and $\Omega$ are clear from the context) iff there is a strictly increasing sequence
$\left<\delta_i\right>_{i<k}$ such that
\begin{enumerate}
 \item $\all\delta\in\N$, $\N\sats$``$\delta$ is Woodin'' if and only if $\ex 
i<k\,(\delta=\delta_i)$.
 \item $\OR(\N)=\sup_{i<\om}(\delta_{k-1}^{+i})^\N$.
\item\label{item:cutpoint} If $\N|\eta$ is a strong 
cutpoint of $\N$ then
$\N|(\eta^+)^\N=\Lp_+^{^\g\Omega,\Gamma}(\N|\eta)$.
 \item\label{item:Qstructure} Let $\xi<\OR(\N)$, where $\N\sats$``$\xi$ 
is 
not Woodin''. Then $C_\Gamma(\N|\xi)\sats$``$\xi$ is not Woodin''.
 \end{enumerate}
We write $\delta^\N_i=\delta_i$; also let $\delta_{-1}^\N=0$ and $\delta_k^\N=\OR(\N)$.\footnote{We could also define a suitable premouse $\N$ as a $\Theta$-g-organized $\Fop$-premouse and all the results that follow in this paper will be unaffected.}

If $\N$ is $1$-suitable, we simply say $\N$ is suitable, and we write $\delta^\N$ for $\delta_0^\N$.

\end{definition}

Let $\N$ be $1$-suitable and let $\xi\in\OR(\N)$ be a limit ordinal such that 
$\N\sats$``$\xi$ isn't Woodin''. Let $Q\pins\N$ be the Q-structure for 
$\xi$. Let $\alpha$ be such that $\xi=\OR(\N|\alpha)$. If $\xi$ is a strong 
cutpoint of $\N$ then $Q\pins\Lp_+^{^\g\Omega,\Gamma}(\N|\xi)$ by clause \ref{item:cutpoint} of the definition.
Assume now that $\N$ is reasonably iterable. If $\xi$ is a strong cutpoint of 
$Q$, our mouse capturing hypothesis combined 
with clause \ref{item:Qstructure} gives that $Q\pins\Lp_+^{^\g\Omega,\Gamma}(\N|\xi)$. If $\xi$ 
is an $\N$-cardinal then indeed $\xi$ is a strong cutpoint of $Q$, since $\N$ 
has only finitely many Woodins. If $\xi$ is not a strong cutpoint of $Q$, then 
by definition, we do not have $Q\pins\Lp_+^{^\g\Omega,\Gamma}(\N|\xi)$. However, 
using 
$*$-translation (see \cite{DMATM}), one can find a level of 
$\Lp_+^{^\g\Omega,\Gamma}(\N|\xi)$ which corresponds to $Q$ (and this level is in $C_\Gamma(\N|\xi)$).

If $\Omega$ is a nice operator (in the sense of \cite{trang2013}, see Footnote \ref{footnote:nice}) and $\Sigma$ is an iteration strategy for a $\Omega$-$\Gamma$-$1$-suitable premouse $\P$ such that $\Sigma$ has branch condensation and is $\Gamma$-fullness preserving (for some pointclass $\Gamma$), then we say that $(\P,\Sigma)$ is a \textit{$\Omega$-$\Gamma$-suitable pair} or just \textit{$\Gamma$-suitable pair} or just \textit{suitable pair} if the pointclass and/or the operator $\Omega$ is clear from the context.


\begin{definition}[Core model induction operators]\label{cmi operator} \index{core model induction operators}Suppose $(\P, \Sigma)$ is a $\G$-$\Omega^*$-suitable pair for some nice operator $\G$ or a hod pair such that $\Sigma$ has branch condensation and is $\Omega^*$-fullness preserving for some inductive-like $\Omega^*$. Let $\Omega = \Sigma$. Assume $\mathrm{Code}(\Omega)$ is self-scaled. We say $J$ is a \emph{$\Sigma$-core model induction operator} or just a \emph{$\Sigma$-cmi operator} if one of the following holds:
\begin{enumerate}
 \item $J$ is a nice $\Omega$-mouse operator (or $\g$-organized $\Omega$-mouse operator) defined on a cone of HC above some $a\in \textrm{HC}$. Furthermore, $J$ condenses finely, relativizes well and determines itself on generic extensions. 
 
\item  For some $\a\in \mathrm{OR}$ such that $\a$ ends either a weak or a strong gap in the sense of \cite{K(R)} and \cite{trang2013}, letting $M=\mathrm{Lp}^{^\gTheta\Omega}(\mathbb{R},\Omega\rest\textrm{HC})|\a$ and $\Gamma = (\Sigma_1)^M$, $M\models \mathsf{AD}^++\mathsf{MC}(\Sigma)$.\footnote{\label{MC}$\textsf{MC}(\Sigma)$ stands for Mouse Capturing relative to $\Sigma$ which says that for $x, y\in \mathbb{R}$, $x$ is $\mathrm{OD}(\Sigma, y)$ (or equivalently $x$ is $\mathrm{OD}(\Omega, y)$) iff $x$ is in some $g$-organized $\Omega$-mouse over $y$. $\mathsf{SMC}$ is the statement that for every hod pair $(\P,\Sigma)$ such that $\Sigma$ is fullness preserving and has branch condensation, $\textsf{MC}(\Sigma)$ holds.} For some transitive $b\in \textrm{HC}$ and some $1$-suitable (or more fully $\Omega$-$\Gamma$-$1$-suitable) $\Omega$-premouse $\Q$ over $b$, $J=\Lambda$, where $\Lambda$ is an $(\omega_1, \omega_1)$-iteration strategy for $\Q$ which is $\Gamma$-fullness preserving, has branch condensation and is guided by some self-justifying-system (sjs) $\vec{A}=(A_i: i<\omega)$ such that for some real $x$, for each $i$, $A_i\in \mathrm{OD}_{b, \Sigma, x}^M$  and $\vec{A}$ seals the gap that ends at $\alpha$. 
\end{enumerate}

When $\Sigma$ is clear from the context or that we don't want to specify $\Sigma$, we simply say $J$ is a cmi operator.
\end{definition}

\begin{remark}
Let $\Gamma, M$ be as in clause 2 above. The (lightface) envelope of $\Gamma$ is defined as: $A\in \textrm{Env}(\Gamma)$ iff for every countable $\sigma\subset\mathbb{R}$ there is some $A'$ such that $A'$ is $\Delta_1$-definable over $M$ from ordinal parameters and $A\cap \sigma = A'\cap \sigma$. For a real $x$, we define Env$(\Gamma(x))$ similarly: here $\Gamma(x) = \Sigma_1(x)^M$ and $A\in \textrm{Env}(\Gamma(x))$ iff for every countable $\sigma\subset \mathbb{R}$ there is some $A'$ that is $\Delta_1(x)$-definable over $M$ from ordinal parameters such that $A\cap \sigma = A'\cap \sigma$. We now let $\bfEnv(\Gamma) = \bigcup_{x\in\mathbb{R}} \textrm{Env}(\Gamma(x))$. Note that $\bfEnv(\Gamma) = \powerset(\mathbb{R})^M$ if $\alpha$ ends a weak gap and $\bfEnv(\Gamma) = \powerset(\mathbb{R})^{\mathrm{Lp}^\Sigma(\mathbb{R})|(\alpha+1)}$ if $\alpha$ ends a strong gap.

In clause 2 above, $\vec{A}$ is Wadge cofinal in $\bfEnv(\Gamma)$ where $\Gamma = \Sigma_1^{M}$.
\end{remark}

The following definitions are obvious generalizations of those defined in \cite{CMI}. For example, see \cite[Definition 3.2.1]{CMI} for the definition of a coarse $(k,U)$-Woodin mouse.

\begin{definition}\label{dfn:coarse_witness}
We say that the coarse mouse witness condition $W^{*,^\g\Omega}_\gamma$ holds if, whenever $U\subseteq \mathbb{R}$ and both $U$ and its complement have scales in $\mathrm{Lp}^{^\gTheta\Omega}(\mathbb{R},\Omega\rest\textrm{HC})|\gamma$, then for all $k< \omega$ and $x \in \mathbb{R}$ there is a coarse $(k,U)$-Woodin mouse $M$ containing $x$ and closed under the strategy $\Lambda$ of $\M_1^{\Omega,\sharp}$ with an $(\omega_1 + 1)$-iteration strategy whose restriction to HC is in $\mathrm{Lp}^{^\gTheta\Omega}(\mathbb{R},\Omega\rest\textrm{HC})|\gamma$.\footnote{We demand the strategy has the property that iterates of $M$ according to the strategy are closed under $\Lambda$.}
\end{definition}
\begin{remark}
By the proof of \cite[Lemma 3.3.5]{CMI}, $W^{*,^\g\Omega}_\gamma$ implies $\mathrm{Lp}^{^\gTheta\Omega}(\mathbb{R},\Omega\rest\textrm{HC})|\gamma \vDash \mathsf{AD}^+$.
\end{remark}
\begin{definition}
An ordinal $\gamma$ is a \emph{critical ordinal} in $\mathrm{Lp}^{^\gTheta\Omega}(\mathbb{R},\Omega\rest\textrm{HC})$ if there is some $U \subseteq \mathbb{R}$
such that $U$ and $\mathbb{R} \backslash U$ have scales in $\mathrm{Lp}^{^\gTheta\Omega}(\mathbb{R},\Omega\rest\textrm{HC})|(\gamma + 1)$ but not in $\mathrm{Lp}^{^\gTheta\Omega}(\mathbb{R},\Omega\rest\textrm{HC})|\gamma$. In other words, $\gamma$
is critical in $\mathrm{Lp}^{^\gTheta\Omega}(\mathbb{R},\Omega\rest\textrm{HC})$ just in case $W^{*,^\g\Omega}_{\gamma+1}$ does not follow trivially from $W^{*,^\g\Omega}_{\gamma}$.
\end{definition}

To any $\Sigma_1$ formula $\theta(v)$ in the language of Lp$^{^\gTheta\Omega}(\mathbb{R},\Omega\rest\textrm{HC})$ we associate formulae $\theta_k(v)$ for $k \in \omega$, such that $\theta_k$ is $\Sigma_k$, and for any $\gamma$ and any real $x$,
\begin{center}
Lp$^{^\gTheta\Omega}(\mathbb{R},\Omega\rest\textrm{HC})|(\gamma+1) \vDash \theta[x] \iff \exists k < \omega$ Lp$^{^\gTheta\Omega}(\mathbb{R},\Omega\rest\textrm{HC})|\gamma \vDash \theta_k[x].$
\end{center}
\begin{definition}\label{prewitness}
Suppose $\theta(v)$ is a $\Sigma_1$ formula (in the language of set theory expanded by a name for $\mathbb{R}$ and a predicate for $^\gTheta\Omega$), and $z$ is a real; then a \emph{$\langle\theta,z\rangle$-prewitness} is an $\omega$-sound $g$-organized $\Omega$-premouse $N$ over $z$ in which there are $\delta_0 < \dots < \delta_9$, $S$, and $T$ such that $N$ satisfies the formulae expressing
\begin{enumerate}[(a)]
\item $\mathsf{ZFC}$,
\item $\delta_0, \ldots, \delta_9$ are Woodin,
\item $S$ and $T$ are trees on some $\omega\times \eta$ which are absolutely complementing in $V^{\mathrm{Col}(\omega,\delta_9)}$, and
\item For some $k <\omega$, $p[T]$ is the $\Sigma_{k+3}$-theory (in the language with names for each real and predicate for $^\gTheta\Omega$) of $\mathrm{Lp}^{^\gTheta\Omega}(\mathbb{R},\Omega\rest\mathrm{HC})|\gamma$, where $\gamma$ is least such that $\mathrm{Lp}^{^\gTheta\Omega}(\mathbb{R},\Omega\rest\mathrm{HC})|\gamma \vDash \theta_k[z]$.
\end{enumerate}
If $N$ is also $(\omega, \omega_1, \omega_1 + 1)$-iterable (as a $\g$-organized $\Omega$-mouse), then we call it a \emph{$\langle\theta, z\rangle$-witness}.
\end{definition}

\begin{definition}\label{def:fine_mouse_witness_cond}
We say that the fine mouse witness condition $W_\gamma^{^g\Omega}$ holds if whenever $\theta(v)$ is a $\Sigma_1$ formula (in the language $\mathcal{L}^+$ of $\g$-organized $\Omega$-premice (cf. \cite{trang2013})), $z$ is a real, and $\mathrm{Lp}^{^\gTheta\Omega}(\mathbb{R},\Omega\rest\textrm{HC})|\gamma\vDash \theta[z]$, then there is a $\langle\theta,z\rangle$-witness $\N$ whose $^{^g}\Omega$-iteration strategy, when restricted to countable trees on $\N$, is in $\mathrm{Lp}^{^\gTheta\Omega}(\mathbb{R},\Omega\rest\textrm{HC})|\gamma$.
\end{definition}

\begin{lemma}
$W_\gamma^{*,^g\Omega}$ implies $W_\gamma^{^g\Omega}$ for limit $\gamma$.
\end{lemma}

The proof of the above lemma is a straightforward adaptation of that of \cite[Lemma 3.5.4]{CMI}. One main point is the use of the $g$-organization: $g$-organized $\Omega$-mice behave well with respect to generic extensions in the sense that if $\P$ is a $g$-organized $\Omega$-mouse and $h$ is set generic over $\P$ then $\P[h]$ can be rearranged to a $g$-organized $\Omega$-mouse over $h$.

\begin{remark}
In light of the discussion above, the core model induction (through $\mathrm{Lp}^{^\gTheta\Omega}(\mathbb{R},\Omega\rest\textrm{HC})$) inductively shows $\mathrm{Lp}^{^\gTheta\Omega}(\mathbb{R},\Omega\rest\textrm{HC})|\gamma\models \sf{AD}^+$ by showing that $W_\gamma^{*,^g\Omega}$ holds for critical ordinals $\gamma$. This, in turn, is done by constructing appropriate $\Omega$-cmi operators ``capturing" the theory of those levels (as specified in Definitions \ref{dfn:coarse_witness} and  \ref{def:fine_mouse_witness_cond}).
\end{remark}

Finally, as in \cite{trangwilson2016Spct}, the maximal model of $\Theta = \theta_\Omega$ is $\mathrm{sLp}^{^\gTheta\Omega}(\mathbb{R},\mathrm{Code}(\Omega))$, an initial segment (possibly strict) of $\mathrm{Lp}^{^\gTheta\Omega}(\mathbb{R},\mathrm{Code}(\Omega))$. 
\begin{definition}\label{def:sLp}
We define $\mathrm{sLp}^{^\gTheta\Omega}(\mathbb{R},\mathrm{Code}(\Omega))$ to be the union of those $\M\lhd \mathrm{Lp}^{^\gTheta\Omega}(\mathbb{R},\mathrm{Code}(\Omega))$ such that whenever $\pi:\M^*\rightarrow \M$ is elementary, $\P\in \pi^{-1}(\mathrm{HC})$, and $\M^*$ is countable and transitive, then $\M^*$ is $X$-$(\omega_1+1)$-iterable with unique strategy $\Lambda$ such that $\Lambda\rest \HC \in \M$.
\end{definition}

In Section \ref{sec:suc}, we will outline the core model induction in the next section, showing that $\mathrm{Lp}^{^\gTheta\Omega}(\mathbb{R},\mathrm{Code}(\Omega)) \models \sf{AD}^+ + \sf{MC}$$(\Omega)$\footnote{$\sf{MC}$$(\Omega)$ states that if $x,y\in\mathbb{R}$ and $x\in OD(y,\Omega)$, then there is a $\Omega$-mouse $\M$ over $y$ such that $\M$ is sound, $\rho_\omega(\M) = \omega$, and $x\in \M$.}  for sufficiently nice $\Omega$. We note that by \cite{sargsyan2014Rmice}, if $M$ is a model of $\sf{AD}^+ + \sf{MC}(\Omega)$ satisfying  $\Theta = \Theta_\Omega$  and $V = L(\powerset(\mathbb{R}))$, then $M$ satisfies that every set of reals $A$ belongs to $\mathrm{sLp}^{^\gTheta\Omega}(\mathbb{R},\mathrm{Code}(\Omega))$. So in fact, in the situation of this paper,  
\begin{center}
$\mathrm{sLp}^{^\gTheta\Omega}(\mathbb{R},\mathrm{Code}(\Omega))=\mathrm{Lp}^{^\gTheta\Omega}(\mathbb{R},\mathrm{Code}(\Omega))$.
\end{center}
For notational simplicity, from now on, we denote $\textrm{Lp}^{^\gTheta\Omega}(\mathbb{R}, \Sigma\rest \textrm{HC})$ by Lp$^\Omega(\mathbb{R})$.

\subsection{Hod mice}\label{sec:hodmice_prelim}


In this paper, a hod premouse $\P$ is one defined as in \cite{ATHM}. The reader is advised to consult \cite{ATHM} for basic results and notations concerning hod premice and mice.

We recall that if $\P$ is a hod premouse and $\gamma$ is an ordinal, then we say $\gamma$ is a \textit{cutpoint} of $\P$ if there is no extender $E$ on the $\P$-sequence such that crt$(E) < \gamma < \textrm{lh}(E)$. We say $\gamma$ is a \textit{strong cutpont} of $\P$ if there is no extender $E$ on the $\P$-sequence such that crt$(E) \leq \gamma < \textrm{lh}(E)$. By $\P|\gamma$, we mean the model $\P$ up to $\gamma$, including the top extender (if one exists); by By $\P||\gamma$, we mean the model $\P$ up to $\gamma$, not including the top extender.

Let us mention some basic first-order properties of a hod premouse $\P$. There are an ordinal $\lambda^\P$ and sequences $\langle(\P(\alpha),\Sigma^\P_\alpha) \ | \ \alpha < \lambda^\P\rangle$ and $\langle \delta^\P_\alpha \ | \ \alpha \leq \lambda^\P  \rangle$ such that 
\begin{enumerate}
\item $\langle \delta^\P_\alpha \ | \ \alpha \leq \lambda^\P  \rangle$ is increasing and continuous and if $\alpha$ is a successor ordinal then $\P \vDash \delta^\P_\alpha$ is Woodin;
\item every Woodin cardinal or limit of Woodin cardinals of $\P$ is of the form $\delta^\P_\alpha$ for some $\alpha$;
\item $\P(0) = Lp_\omega(\P|\delta_0)^\P$; for $\alpha < \lambda^\P$, $\P(\alpha+1) = (Lp_\omega^{\Sigma^\P_\alpha}(\P|\delta_{\alpha+1}))^\P$;\footnote{$\P(\alpha+1)$ is a ($g$-organized) $\Sigma_\alpha$-premouse in the sense defined above.} for limit $\alpha\leq \lambda^\P$, $\P(\alpha) = (Lp_\omega^{\oplus_{\beta<\alpha}\Sigma^\P_\beta}(\P|\delta_\alpha))^\P$;
\item $\P \vDash \Sigma^\P_\alpha$ is a $(\omega,o(\P),o(\P))$\footnote{This just means $\Sigma^\P_\alpha$ acts on all stacks of $\omega$-maximal, normal trees in $\P$.}-strategy for $\P(\alpha)$ with hull condensation;
\item if $\alpha < \beta < \lambda^\P$ then $\Sigma^\P_\beta$ extends $\Sigma^\P_\alpha$.
\end{enumerate}
We will write $\delta^\P$ for $\delta^\P_{\lambda^\P}$ and $\Sigma^\P=\oplus_{\beta<\lambda^\P}\Sigma^\P_{\beta}$. Note that $\P(0)$ is a pure extender model. Suppose $\P$ and $\Q$ are two hod premice. Then $\P\trianglelefteq_{hod}\Q$\index{$\trianglelefteq_{hod}$} if there is $\a\leq\lambda^\Q$ such that $\P=\Q(\a)$. We say then that $\P$ is a \textit{hod initial segment} of $\Q$. We say $(\P,\Sigma)$ is a \textit{hod pair} if $\P$ is a hod premouse and $\Sigma$ is a strategy for $\P$ (acting on countable stacks of countable normal trees) such that $\Sigma^\P \subseteq \Sigma$ and this fact is preserved under $\Sigma$-iterations. Typically, we will construct hod pairs $(\P,\Sigma)$ such that $\Sigma$ has hull condensation, branch condensation, and is $\Gamma$-fullness preserving for some pointclass $\Gamma$.

See \cite{ATHM} for the definition of hulls of an iteration tree/stack and \cite{SteelNormalization} for a more general notion of a pseudo-hull of a stack.

\begin{definition}\label{dfn:condensation}
Let $\P$ be a hod premouse in the sense of \cite{ATHM} and $\Sigma$ be an iteration strategy for $\P$.
\begin{enumerate}[(a)]
\item $\Sigma$ has \textit{branch condensation} if whenever $\vec{\T},\vec{\U}$ are stacks according to $\Sigma$, $b = \Sigma(\vec{\T})$ is a non-dropping branch, and $c$ is a cofinal, nondropping branch of $\vec{\U}$ such that there is an elementary $\sigma: \M^{\vec{\U}}_c \rightarrow \M^{\vec{\T}}_b$ with the property that $\pi^{\vec{\T}}_b = \sigma\circ \pi^{\vec{\U}}_c$, then $c = \Sigma(\vec{\U})$.
\item $\Sigma$ has \textit{strong hull condensation} if whenever $\vec{\T}$ is according to $\Sigma$ and $\vec{\U}$ is a pseudo-hull of $\vec{\T}$ then $\vec{\U}$ is according to $\Sigma$. $\Sigma$ has \textit{hull condensation} if whenever $\vec{\T}$ is according to $\Sigma$ and $\vec{\U}$ is a hull of $\vec{\T}$ then $\vec{\U}$ is according to $\Sigma$. 
\end{enumerate}

\end{definition}

Strong hull condensation easily implies hull condensation because every hull is a pseudo-hull. We note that strategies for hod pairs are assumed to have hull condensation, but it is not clear that hod mouse strategies constructed in \cite{ATHM} can have strong hull condensation.  See \cite{ATHM} for the definition of $\Gamma(\P,\Sigma)$. Roughly,  $\Gamma(\P,\Sigma)$ is the pointclass generated by $\Sigma$. In the case $\lambda^\P$ is a limit ordinal, $\Gamma(\P,\Sigma)$ is the set of $B$ such that there is some $(\Q,\Lambda)\in B(\P,\Sigma)$, $B\leq_w \Lambda$. See \cite{ATHM} for the definition of $\Gamma(\Q,\Sigma)$ in the case $\lambda^\Q$ is a successor ordinal. In Lemma \ref{claim:bcshc}, we show that if $(\P,\Sigma)$ is a hod pair such that $\Sigma$ has branch condensation and $\Gamma(\P,\Sigma)$-fullness preserving then $\Sigma$ has strong hull condensation.  Lemma \ref{claim:bcshc} appears to be a new fact in hod mice theory at the level of ``$\sf{AD}_\mathbb{R} + $$\Theta$ is regular." \footnote{Lemma \ref{claim:bcshc} should also hold for hod mice in a minimal model of $\sf{LSA}$ but we have not checked all details of this claim.} The lemma is used essentially in the proof of Lemma \ref{lem:pullback}, which is a key part in the proof of Theorem \ref{thm:DI}.

The reader should also consult \cite{ATHM} for the definition of $B(\Q,\Sigma)$ and $I(\Q,\Sigma)$. Roughly speaking, $B(\Q,\Sigma)$ is the collection of all hod pairs which are strict hod initial segments of a $\Sigma$-iterate of $\Q$ and $I(\Q,\Sigma)$ is the collection of all $\Sigma$-iterates of $\Q$. In the case $\lambda^\Q$ is limit, the pointclass $\Gamma(\Q,\Sigma)$ is the collection of $A\subseteq \mathbb{R}$ such that $A$ is Wadge reducible to some $\Psi$ for which there is some $\R$ such that $(\R,\Psi)\in B(\Q,\Sigma)$.  If $(\P, \Sigma)$ is a hod pair, and $\vec{\T}$ is according to $\Sigma$ with last model $\Q$, then we write $\Sigma_{\Q,\vec{\T}}$ for the $\vec{\T}$-tail strategy of $\Q$ induced by $\Sigma$, i.e. $\Sigma_{\Q,\vec{\T}}(\vec{\U}) = \Sigma(\vec{\T}^\smallfrown \vec{\U})$.

Suppose $(\Q,\Sigma)$ is a hod pair such that $\Sigma$ has hull condensation. We say $\P$ is a $(\Q,\Sigma)$-\emph{hod premouse} if there are an ordinal $\lambda^\P$ and sequences $\langle(\P(\alpha),\Sigma^\P_\alpha) \ | \ \alpha < \lambda^\P\rangle$ and $\langle \delta^\P_\alpha \ | \ \alpha \leq \lambda^\P  \rangle$ such that 
\begin{enumerate}
\item $\langle \delta^\P_\alpha \ | \ \alpha \leq \lambda^\P  \rangle$ is increasing and continuous and if $\alpha$ is a successor ordinal then $\P \vDash \delta^\P_\alpha$ is Woodin;
\item every Woodin cardinal or limit of Woodin cardinals of $\P$ is of the form $\delta^\P_\alpha$ for some $\alpha$;
\item $\P(0) = Lp^{\Sigma}_\omega(\P|\delta_0)^\P$ (so $\P(0)$ is a $\Sigma$-premouse built over $\Q$); for $\alpha < \lambda^\P$, $\P(\alpha+1) = (Lp_\omega^{\Sigma\oplus\Sigma^\P_\alpha}(\P|\delta_\alpha))^\P$; for limit $\alpha\leq \lambda^\P$, $\P(\alpha) = (Lp_\omega^{\oplus_{\beta<\alpha}\Sigma^\P_\beta}(\P|\delta_\alpha))^\P$;
\item $\P \vDash \Sigma\cap \P$ is a $(\omega,o(\P),o(\P))$-strategy for $\Q$ with hull condensation;
\item $\P \vDash \Sigma^\P_\alpha$ is a $(\omega,o(\P),o(\P))$-strategy for $\P(\alpha)$ with hull condensation;
\item if $\alpha < \beta < \lambda^\P$ then $\Sigma^\P_\beta$ extends $\Sigma^\P_\alpha$.
\end{enumerate}
Inside $\P$, the strategies $\Sigma^\P_\alpha$ act on stacks above $\Q$ and every $\Sigma^P_\alpha$ iterate is a $\Sigma$-premouse. Again, we write $\delta^\P$ for $\delta^\P_{\lambda^\P}$ and $\Sigma^\P=\oplus_{\beta<\lambda^\P}\Sigma^\P_{\beta}$. We say $(\P,\Lambda)$ is a $(\Q,\Sigma)$-\emph{hod pair} if $\P$ is a $(\Q,\Sigma)$-hod premouse and $\Lambda$ is a strategy for $\P$ such that $\Sigma^P\subseteq \Lambda$ and this fact is preserved under $\Lambda$-iterations. The reader should consult \cite{ATHM} for the definition of $B(\Q,\Sigma)$ and $I(\Q,\Sigma)$. Roughly speaking, $B(\Q,\Sigma)$ is the collection of all hod pairs which are strict hod initial segments of a $\Sigma$-iterate of $\Q$ and $I(\Q,\Sigma)$ is the collection of all $\Sigma$-iterates of $\Q$. In the case $\lambda^\Q$ is limit, the pointclass $\Gamma(\Q,\Sigma)$ is the collection of $A\subseteq \mathbb{R}$ such that $A$ is Wadge reducible to some $\Psi$ for which there is some $\R$ such that $(\R,\Psi)\in B(\Q,\Sigma)$. See \cite{ATHM} for the definition of $\Gamma(\Q,\Sigma)$ in the case $\lambda^\Q$ is a successor ordinal. If $(\P, \Sigma)$ is a hod pair, and $\vec{\T}$ is according to $\Sigma$ with last model $\Q$, then we write $\Sigma_{\Q,\vec{\T}}$ for the $\vec{\T}$-tail strategy of $\Q$ induced by $\Sigma$, i.e. $\Sigma_{\Q,\vec{\T}}(\vec{\U}) = \Sigma(\vec{\T}^\smallfrown \vec{\U})$.

Suppose $(\R,\Lambda)$ is a hod pair and $\Gamma$ is a nice pointclass. We say that $\Lambda$ is \textit{$\Gamma$-$\Q$-structure guided} if whenever $\T$ is according to $\Lambda$ and short, then $\Lambda(\T) = b$ is such that $\Q(b,\T)$ exists and the phalanx $\Phi(\T^\smallfrown b)$\footnote{This is the set of models in the tree $\T^\smallfrown b$.} is $(\omega_1,\omega_1)$-iterable with unique strategy in $\Gamma$. We show in essence that the branch $b$ must be unique in Lemma \ref{lem:uniqueQ}. We also note that if $\delta(\T)$ is a cutpoint of $\Q(b,\T)$ then the phalanx iterability condition reduces to the iterability of $\Q(b,\T)$ above $\delta(\T)$.

Suppose $\P$ is $\Sigma$-suitable and $A\subseteq \mathbb{R}$ is $OD_\Sigma$. We say $\P$ \textit{weakly term captures} $A$ if letting $\d=\d^\P$, for each $n<\omega$ there is a term relation $\tau\in \P^{Coll(\omega, (\d^{+n})^\P)}$ such that for comeager many $\P$-generics $g\subseteq Coll(\omega, (\d^{+n})^\P)$, we have $\tau_g=\P[g]\cap A$. We say $\P$ \textit{term captures} $A$ if the equality holds for all generics. Given a $\Sigma$-suitable $\P$  and an $OD_\Sigma$ set of reals $A$, we let $\tau_{A, n}^\P$ be the standard name for a set of reals in $\P^{Coll(\omega, (\d^{+n})^\P)}$ witnessing the fact that $\P$ weakly captures $A$ and let $$\gamma^\P_A = sup(\d^\P \cap Hull^{\P}_1( \{ \tau^\P_{A,n} : n<\omega \})). $$ See \cite{trang2013, ATHM} for all relevant definitions; in particular, discussions on $\Sigma$-suitable premice and term capturing are given in \cite[Section 3]{trang2013}. We let 
\begin{equation}\label{eqn:terms}
f_A(\P)=\la \tau^\P_{A, n} : n<\omega\ra.
\end{equation}

Suppose $(\R,\Lambda)$ is a hod pair and $\lambda^\R = \alpha+1$ for some $\alpha\geq 0$, where $\lambda^\R$ is the order type of the set $\{\delta : \delta \textrm{ is either a Woodin cardinal or a limit of Woodin cardinals in } \R\}$; we will write $\delta^\R_\alpha$ for the $\alpha$-th member of this set. Recall the notations $(\R^-,\Lambda_{\R^-})$, $\mathbb{B}(\R^-,\Lambda_{\R^-})$ from \cite{ATHM}.\footnote{$\R^- = \R(\alpha-1)$ and $\Lambda_{\R^-}$ is just $\Lambda_{\R(\alpha-1)}$. In the case $\alpha=0$, $(\R^-,\Lambda_{\R^-}) = (\emptyset,\emptyset)$.} \cite[Lemma 5.19]{ATHM} gives that $\sf{AD}^+$ implies there is some tail $(\S,\Psi)$ of $(\R,\Lambda)$ and some $\vec{B}= \{B_i : i<\omega\}$ that \textit{strongly guides $\Psi$}. This means that 
\begin{itemize}
\item $\Lambda$ is $\Gamma$-$\Q$-structure guided, where $\Gamma = \Gamma(\R,\Lambda)$.
\item There are terms $(\tau_i^\S = \tau_{B_i,0}^\S  : i < \omega, \tau_i^\S\in \S^{Coll(\omega,\delta^\S)})$ for $B_i$ such that whenever $k: \S \rightarrow \Q$ is an iteration map by $\Psi$ of a maximal tree, then for each $i<\omega, k(\tau_i^\S) = \tau_{B_i,0}^\Q$ is the term that captures $B_i$ over $\Q$, sup$\{\gamma^\Q_{B_i}: i < \omega\}=\delta^\S$, the branch $b$ giving rise to the embedding $k$ is the unique branch whose branch embedding moves the terms for $B_i$'s correctly, and whenever $\vec{\T}$ is according to $\Psi$ with branch embedding $\pi$, $\vec{\U}$ is according to $\Psi$, and suppose $b$ is a cofinal branch of $\vec{\U}$ such that there is an elementary map $\sigma: \M^{\vec{\U}}_b \rightarrow \M^{\vec{\T}}$ such that $\sigma \circ \pi^{\vec{\U}}_b = \pi^{\vec{\T}}$, then for each $i$,
\begin{center}
$\sigma^{-1}(\tau^{\M^{\vec{\T}}}_{B_i,0}) = \tau^{\M^{\vec{\U}}_b}_{B_i,0} = \pi^{\vec{\U}}_b(\tau^\S_{B_i,0})$.
\end{center}
\end{itemize}
When we don't want to specify the $B_i$'s or the particular $B_i$'s are not important to specify, we simply say $\Psi$ is strongly guided. The above notion of strongly guided can be defined in an obvious way for $(\R,\Lambda)$, where $\lambda^\R = \alpha+n$ for some $n<\omega$. We omit details and refer the reader to \cite{ATHM} for a full discussion. The next section will elaborate more on this topic in the context of the HOD analysis.


\begin{definition}[$\Gamma$-Fullness preservation]\label{gamma fullness preservation} Suppose $(\P, \Sigma)$ is a hod pair such that $\P\in HC$ and $\Gamma$ is a nice pointclass. We say $\Sigma$ is $\Gamma$-fullness preserving if $\Sigma$ is $\Gamma$-$\Q$-structure guided and the following holds for all $(\Q, \VT)\in I(\P, \Sigma)$.
\begin{enumerate}

\item For all limit $\alpha < \lambda^\Q$, letting $\R= \Q(\alpha)$, then 
\begin{center}
$\R=Lp^{\Gamma, \oplus_{\beta<\alpha}\Sigma_{\R(\beta), \VT}}_\omega(\R|\d^\R)$.
\end{center}
\item  For all successor $\alpha < \lambda^\Q$, letting $\R= \Q(\alpha)$ and $\beta=\alpha-1$, 
\begin{center}
$\R=Lp^{\Gamma, \Sigma_{\R(\beta), \VT}}_\omega(\R|\d^\R)$.
\end{center}

\item If $\eta$ is a cardinal strong cutpoint of $\Q$, letting $\alpha$ be the largest such that $\Q(\alpha)\lhd \Q|\eta$ and $\R = \Q(\alpha)$, then 
\begin{center}
 $\Q|(\eta^+)^\Q=Lp^{\Gamma, \Sigma_{\R, \VT}}(\Q|\eta)$.
\end{center}

\item Furthermore, letting for $\alpha+1\leq \lambda^\Q$,
\begin{center}
$U_{\Q(\alpha),\Sigma} = \{ (x,y)\in \mathbb{R}^2 : x\in\mathbb{R} \textrm{ codes a countable set } a \textrm{ and } y \textrm{ codes a sound } \Sigma_{\Q(\alpha)}\textrm{-mouse } \M \newline  \textrm{ over }  a \textrm{ whose unique strategy is in } \Gamma \textrm{ such that } \rho(\M) = a \}$,
\end{center}
and 
\begin{center}
$W_{\Q(\alpha),\Sigma} = \{ (x,y,z)\in \mathbb{R}^3 : (x,y)\in U_{\Q(\alpha),\Sigma} \textrm{ and } z \textrm{ codes an iteration tree on the mouse } \M \newline \textrm{coded by } y\}$,
\end{center}
then whenever $(\vec{\U},\R)\in I(\Q(\alpha+1),\Sigma_{\Q(\alpha+1),\vec{\T}})$ such that $\vec{\U}$ only uses extenders with critical points above $\delta_\alpha^\Q$ and its images along branch embeddings of $\vec{\U}$, we have 
\begin{center}
$\pi^{\vec{\U}}(f_A(\Q)) = f_A(\R)$,
\end{center}
where $A = U_{\Q(\alpha),\Sigma}\oplus W_{\Q(\alpha),\Sigma}$ and $f_A$ is defined in \eqref{eqn:terms} below.
\end{enumerate}
\end{definition}

\begin{remark}
In \cite{ATHM}, clauses (1)--(3) comprise the definition of fullness preservation of $\Sigma$; if in addition, clause (4) holds for $\Sigma$, then $\Sigma$ is said to be super fullness preserving (with respect to $\Gamma$). We simplify the terminology by combining these two notions into one definition. 
\end{remark}

Under $\textsf{AD}^+$ and the hypothesis that there are no models of $\textsf{AD}_\mathbb{R}+{}$``$\Theta$ is regular,'' \cite{ATHM} constructs  hod pairs that are fullness preserving and have branch condensation (see \cite{ATHM} for a full discussion of these notions). Such hod pairs are particularly important for our computation as they are points in the direct limit system giving rise to \textrm{HOD} of $\textsf{AD}^+$ models. Under $\sf{AD}^+$, for hod pairs $(\M_\Sigma, \Sigma)$, if $\Sigma$ is a strategy with branch condensation and $\VT$ is a stack on $\M_\Sigma$ with last model $\N$, then $\Sigma_{\N, \VT}$ is independent of $\VT$. Therefore, later on we will omit the subscript $\VT$ from $\Sigma_{N, \VT}$ whenever $\Sigma$ is a strategy with branch condensation and $\M_\Sigma$ is a hod mouse. In a core model induction, at the moment $(\M_\Sigma,\Sigma)$ is constructed we don't quite have an $\textsf{AD}^+$-model $M$ such that $(\M_\Sigma,\Sigma)\in M$, but we do know that every $(\R,\Lambda)\in B(\M_\Sigma,\Sigma)$ belongs to such a model. We then can show (using our hypothesis) that $(\M_\Sigma,\Sigma)$ belongs to an $\textsf{AD}^+$-model.

We briefly review definitions and notations related to the analysis of stacks in \cite[Section 6.2]{ATHM}; see \cite[Section 6.2]{ATHM} for a more detailed discussion. These notions will be useful in Section \ref{sec:lim}. Suppose $\P$ is a hod premouse and $\vec{\T}$ is a stack on $\P$. Let $\S$ be a model that appears in $\vec{\T}$. By $\vec{\T}_{\leq\S}$ we mean the part of $\vec{\T}$ up to and including $\S$ (according to the tree order of $\vec{\T}$), we define $\vec{\T}_{\geq \S}, \vec{\T}_{<\S},\vec{\T}_{>\S}$ similarly. We let $(\M_\alpha, \T_\alpha: \alpha<\eta)$ be the normal components of $\vec{\T}$, i.e. $\M_0 = \P$, $\T_\alpha$ is a normal tree on $\M_\alpha$, and $\M_{\alpha+1} = \M^{\T_\alpha}$. We say $\R$ is a \textit{terminal node} of $\vec{\T}$ if for some $\alpha,\beta$, $\R = \M^{\T_\alpha}_\beta$ and $\pi^{\T_\alpha}_{0,\beta}$ is defined. We say $\R$ is a \textit{non-trivial terminal node} of $\vec{\T}$ if letting $(\alpha,\beta)$ witness that $\R$ is a terminal node of $\vec{\T}$, the extender $E^{\T_\alpha}_\beta$ is applied to $\R$ in the tree $\T_\alpha$ to obtain the model $\M^{\T_\alpha}_{\beta+1}$. We write $tn(\vec{\T})$ for the set of terminal nodes of $\vec{\T}$ and $ntn(\vec{\T})$ for the set of non-trivial terminal nodes of $\vec{\T}$.

For $\Q,\R\in tn(\vec{\T})$, we write $\Q\prec^{\vec{\T}} \R$ if the $\Q$-to-$\R$ iteration embedding in $\vec{\T}$ exists, and we write $\pi^{\vec{\T}}_{\Q,\R}$ for this embedding. We write $\Q\prec^{\vec{\T},s} \R$ if letting $\vec{\U}$ be the part of $\vec{\T}$ between $\Q$ and $\R$, then $\vec{\U}$ is an iteration on $\Q$. We write $\vec{\T}_{\Q,\R}$ for $\vec{\U}$. 

Let $C\subseteq tn(\vec{\T})$. We say $C$ is \emph{linear} (\emph{strongly linear} respectively) if $C$ is linearly ordered by $\prec^{\vec{T}}$ ($\prec^{\vec{T},s}$ respectively). We say $C$ is \emph{closed} if $C$ is strongly linear and whenever $\alpha$ is a limit point of $C$, then letting $\R$ be the direct limit of $C\rest \alpha$ (under the iteration embeddings), we have $\R\in C$. We say $C$ is \emph{cofinal} if for every $\S\in \vec{\T}$, there are $\Q,\R\in C$ such that $\Q\prec^{\vec{\T},s} \R$ and $\S$ is in $\vec{\T}_{\Q,\R}$. Note that if $\vec{\T}$ doesn't have a last model, but there is a strongly closed and cofinal $C\subseteq tn(\vec{\T})$, then $C$ uniquely determines a cofinal branch of $\vec{\T}$. If such a $C$ doesn't exist, then $\eta$ is a successor ordinal, say $\eta = \alpha+1$. Let $\U = \vec{\T}_\alpha$ and $D=\{\S\in tn(\U) : \U_{\geq \S} \textrm{ is a tree on } \S\}$. In this case $D$ has a $\prec^{\vec{\T},s}$-largest element and we write $\S_{\vec{\T}}$ for this element.
Then $\vec{\T}_{\S_{\vec{\T}}}$ is a normal tree based on $\S_{\vec{\T}}(\beta+1)$ and above $\delta_\beta^{\S_{\vec{\T}}}$ for some $\beta<\lambda^{\S_{\vec{\T}}}$.

\subsection{$\H$ and $\H_\Sigma$ under $\sf{AD}^+$}
Suppose $\Sigma$ is an iteration strategy of some hod mouse $\Q$ and suppose $\Sigma$ is fullness preserving (see \cite{ATHM}) and has branch condensation. Assume further that  $V=L(\powerset(\mathbb{R}))$ and $\sf{MC}$$(\Sigma)$ holds and $\Theta=\theta_\Sigma$. 

\begin{definition}[$S(\Gamma, \Sigma)$ and $F(\Gamma, \Sigma)$] 
Suppose $\Gamma$ is a pointclass. Let $S(\Gamma, \Sigma)=\{ \Q: \Q$ is $ \Sigma$-suitable$\}$. Also, we let $F(\Gamma, \Sigma)$ be the set of functions $f$ such that $dom(f)=S(\Gamma, \Sigma)$ and for each $\P\in S(\Gamma, \Sigma)$, $f(\P)\subseteq \P$ and $f(\P)$ is amenable to $\P$, i.e., for every $X\in \P$, $X\cap f(\P)\in \P$.
\end{definition}

We let $\Gamma=\powerset(\mathbb{R})$ and for the duration of this subsection, we drop $\Gamma$ from our notation whenever it is unambiguous to do so. Thus, a $\Sigma$-suitable premouse is a $\Sigma$-$\Gamma$-suitable premouse etc. We remark that by \cite{sargsyan2014Rmice}, 
\begin{center}
$V = L(\textrm{Lp}^\Sigma(\mathbb{R}))$.
\end{center}
Also, we allow for the case $(\P,\Sigma) = (\emptyset, \emptyset)$, in which case $V = L(\textrm{Lp}(\mathbb{R}))$ and $\H_\Sigma = \H$. The following lemma is essentially due to Woodin and the proof for mice can be found in \cite{CMI}.

\begin{lemma} Suppose $\P$ is $\Sigma$-suitable and $A\subseteq \mathbb{R}$ is $OD_\Sigma$. Then $\P$ weakly term captures $A$. Moreover, there is a $\Sigma$-suitable $\Q$ which term captures $A$.
\end{lemma}

%


 The following lemma is one of the most fundamental lemmas used to compute $\H$ and it is originally due to Woodin. Again, the proof can be found in \cite{CMI}. See also \cite{CMI} for detailed discussions of related standard notions like $f$-iterability and $f$-quasi-iterability.

\begin{theorem}\label{existence of quasi-iterable premice} For each $f\in F_{\Sigma, od}$, there is a $\Sigma$-suitable premouse $\P$ which is strongly $f$-iterable.
\end{theorem}

To save some ink, in what follows, we will sometimes say $A$-iterable instead of $f_A$-iterable and similarly for other notions. Also, we will use $A$ in our subscripts instead of $f_A$.

 Given $\P\in S(\Gamma, \Sigma)$ and $f\in F_{\Sigma,od}$ we let $f_n(\P)=f(\P)\cap \P|((\d^\P)^{+n})^\P$. Then $f(\P)=\bigcup_{n<\omega}f_n(\P)$. We also let 
\begin{center}
$\gg^\P_{f}=sup(\d^\P \cap Hull^{\P}_1( \{ f_n(\P) : n<\omega \})) $.
\end{center}
Notice that
\begin{center}
$\gg^\P_{f}=\d^\P \cap Hull^{\P}_1(\gg_{f}^\P\cup\{ f_n(\P) : n<\omega \} )$.
\end{center}
We then let 
\begin{center}
$H_{f}^\P =Hull^{\P}_1(\gg^\P_{f}\cup \{ f_n(\P) : n<\omega \} )$.
\end{center}
If $\P\in S(\Gamma, \Sigma)$, $f\in F_{\Sigma,od}$, and $i: \P\rightarrow \Q$ is an embedding, then we let $i(f(\P))=\bigcup_{n<\omega}i(f_n(\P))$.

The following are the next block of definitions that routinely generalize into our context: (1) $(f, \Sigma)$-iterability,
(2) $\vec{b}=\la b_k: k< m\ra$ witnesses $(f, \Sigma)$-iterability for $\VT=\la \T_k, \P_k : k< m\ra$, and
(3) strong $(f, \Sigma)$-iterability.

If $\P$ is strongly $(f, \Sigma)$-iterable and $\VT$ is a $(\Gamma, \Sigma)$-correctly guided finite stack on $\P$ with last
model $\R$ then we let
\begin{center}
$\pi^{\Sigma}_{\P, \R, f}:H_f^\P\rightarrow H_f^\R$
\end{center}
be the embedding given by any $\vec{b}$ which witnesses the $(f, \Sigma)$-iterability of $\VT$, i.e., fixing
$\vec{b}$ which witnesses $f$-iterability for $\VT$,
\begin{center}
$\pi^\Sigma_{\P, \R, f} =\pi_{\VT, \vec{b}}\rest H_f^\P$.
\end{center}
Clearly, $\pi^{\Sigma}_{\P, \R, f}$ is independent of $\VT$ and $\vec{b}$. Here we keep $\Sigma$ in our notation for $\pi^{\Sigma}_{\P, \R, f}$ because it depends on a $(\Gamma, \Sigma)$-correct iteration. It is conceivable that $\R$ might also be a $(\Gamma, \Lambda)$-correct iterate of $\P$ for another $\Lambda$, in which case $\pi^{\Sigma}_{\P, \R, f}$ might be different from $\pi^{\Lambda}_{\P, \R, f}$. However, the point is that these embeddings agree on $H_f^\P$. 

Given a finite sequence of functions $\vec{f}=\la f_i : i<n\ra$ in $F_{\Sigma,od}$, we let $\oplus_{i<n}f_i\in F_{\Sigma,od}$ be the function given by $(\oplus_{i<n}f_i)(\P)=\la f_i(\P): i<n\ra$. We set $\oplus\vec{f}= \oplus_{i<n}f_i$. 

We let $F = F_{\Sigma,od}$ and 
\begin{center}
$\mathcal{I}_{F, \Sigma}=\{ (\P, \vec{f}): \P\in S(\Gamma, \Sigma)$, $\vec{f}\in F^{<\omega}$ and $\P$ is strongly $\oplus\vec{f}$-iterable$\}$
\end{center}
and 
\begin{center}
$\mathcal{F}_{F, \Sigma}=\{ H^\P_f: (\P, f)\in \mathcal{I}_{F, \Sigma}\}$.
\end{center}
We then define $\preceq_{F, \Sigma}$ on $\mathcal{I}_{ F, \Sigma}$ by letting $(\P, \vec{f})\preceq_{F, \Sigma} (\Q, \vec{g})$ iff $\Q$ is a $\Sigma$-correct iterate of $\P$ and $\vec{f}\subseteq \vec{g}$. Given $(\P, \vec{f})\preceq_{F, \Sigma} (\Q, \vec{g})$, we have 
\begin{center}
$\pi^\Sigma_{\P, \Q, \vec{f}}: H^\P_{\oplus\vec{f}}\rightarrow H^\Q_{\oplus\vec{f}}.$
\end{center}
Notice that $\preceq_{F, \Sigma}$ is directed. Let then $\M_{\infty,F, \Sigma}$ be the direct limit of $(\mathcal{F}_{F, \Sigma}, \preceq_{F, \Sigma})$ under the maps $\pi^\Sigma_{\P, \Q, \vec{f}}$. Given $(\P, \vec{f})\in \mathcal{I}_{F, \Sigma}$, we let $\pi^\Sigma_{\P, \vec{f}, \infty}: H^\P_{\oplus\vec{f}}\rightarrow \M_{\infty, F, \Sigma}$ be the direct limit embedding. Let $$\M_\infty=\M_{\infty, F, \Sigma}.$$ 

\begin{theorem}[Woodin, \cite{CMI}]\label{hod theorem} $\d^{\M_\infty}=\Theta$, $\M_\infty\in \H_\Sigma$, and 
\begin{center}
$\M_\infty|\Theta=(V_\Theta^{\H_\Sigma}, \vec{E}^{\M_\infty|\Theta}, S^{\M_\infty}, \in)$,
\end{center}
 where $S^{\M_\infty}$ is the predicate of $\M_\infty$ describing $\Sigma$.
\end{theorem}

\begin{remark}\label{rmk:one_cardinal}
In some of the arguments below, for convenience, we actually use the ``one cardinal" version of suitability. More precisely, for $(\P, f)\in \mathcal{I}_{F,\Sigma}$ we consider direct limits of $(\hat{\P}, \hat{\vec{f}})$ where $\delta = \delta^\P$, $\hat{\P}= \P|(\delta^+)^\P$, and $\hat{\vec{f}} = \vec{f}(\P)\cap \P|(\delta^+)^\P$. We define $\gamma^{\hat{\P}}_{\hat{f}} = sup(\d^\P \cap Hull^{\P}_1( \{ f_0(\P) \}))$ etc. We let $\hat{\M}_\infty$ be the direct limit of such pairs $(\hat{\P}, \hat{\vec{f}})$. Then it is easy to see also that $\hat{\M}_\infty|\Theta = (V_\Theta^{\H_\Sigma}, \vec{E}^{\M_\infty|\Theta}, S^{\M_\infty}, \in)$.
\end{remark}

Finally, if $a\in H_{\omega_1}$ is self-wellordered then we could define $\M_\infty(a)$ by working with $\Sigma$-suitable premice over $a$. Everything we have said about $\Sigma$-suitable premice can also be said about $\Sigma$-suitable premice over $a$, and in particular the equivalent of \rthm{hod theorem} can be proven using $\H_{(\Sigma, a)\cup\{a\}}$ instead of $\H_\Sigma$ and $\M_\infty(a)$ instead of $\M_\infty$.

\cite{ATHM} computes HOD (up to $\Theta$) in models of $(V=L(\powerset(\mathbb{R})))+\textsf{SMC} + \textsf{AD}_\mathbb{R}$ below $\textsf{AD}_\mathbb{R} + {}$``$\Theta$ is regular'' by exhibiting a hod premouse $\M_\infty$ satisfying 
\begin{enumerate}
\item $\M_\infty \in \H$.
\item $\M_\infty$ is a hod premouse.
\item $\M_\infty|\Theta = (V_\Theta^{\H}, \vec{E}^{\M_\infty|\Theta}, S^{\M_\infty}, \in)$, where $S^{\M_\infty|\Theta}$ is the predicate for strategies of hod initial segments of $\M_\infty|\Theta$.
\end{enumerate}
Here $\sf{SMC}$ is Strong Mouse Capturing, which is the statement that for any $x,y\in\mathbb{R}$, if $x\in OD_{y,\Sigma}$ where $(\P,\Sigma)$ is a hod pair such that $\Sigma$ has branch condensation and is fullness preserving, then $x$ is in a $\Sigma$-mouse $\M$ over $y$. We call $\M_\infty$ the \textit{hod limit}. Here $\M_\infty = \bigcup_{(\Q,\Lambda)} \M_\infty(\Q,\Lambda)$, where $(\Q,\Lambda)$ is a hod pair with branch condensation and is fullness preserving and $\M_\infty(\Q,\Lambda)$ is the direct limit of all (non-dropping) $\Lambda$-iterates of $\Q$.

\subsection{Strategies with strong hull condensation pulls back}\label{sec:pullback}

\begin{definition}\label{def:reasonable}
We say a hod pair $(\P,\Sigma)$ \textit{reasonable} if it has the following additional properties:
\begin{itemize}
\item $\Sigma$ has branch condensation.
\item $\Sigma$ is $\Gamma(\P,\Sigma)$-fullness preserving.
\end{itemize}
\end{definition}

We will show that properties listed above for $\Sigma$ hold for hold mice constructed in this paper. For the next several proofs, the reader is advised to review \cite{ATHM} for basic properties and terminologies of hod pair strategies. See also \cite[Lemma 3.18]{trang2013} for a similar argument.

\begin{lemma}\label{lem:uniqueQ}
Suppose $(\R,\Lambda)$ is a reasonable hod pair. Let $\Gamma = \Gamma(\R,\Lambda)$. Suppose $\vec{\U}$ is according to $\Lambda$ with the following properties:
\begin{itemize}
\item $\vec{\U} = \vec{\U}_0^\smallfrown \vec{\U}_1$, where $\vec{\U}_0 = \W^\smallfrown d$, where $d = \Lambda(\W)$, 
\item letting $\S = \M^{\W}_d$, there is $\beta < \lambda^\S$ such that the set of generators used in $\vec{\U}_0$ $\alpha(\vec{\U}_0)\subset (\delta^\S_\beta)^{<\omega}$,
\item $\vec{\U}_1$ is based on $\S(\beta+1)$ and is above $\delta^\S_\beta$, 
\item suppose $b$ is a cofinal well-founded branch such that $\Q(b, \vec{\U}_1)$ exists and the phalanx $\Phi(\vec{\U}_1^\smallfrown b)$ is iterable in $\Gamma$.
\end{itemize}
Then $b = \Lambda_{\vec{\U}_0,\S}(\vec{\U}_1)$.
\end{lemma}
\begin{proof}
Let $\Sigma = \Lambda_{\vec{\U}_0,\S}\rest \S(\beta)$. Let $c =  \Lambda_{\vec{\U}_0,\S}(\vec{\U}_1)$. We want to show $b=c$. There are two cases. 

Suppose $\delta(\vec{\U}_1) =_{def}\delta$ is a cutpoint of $\Q(b, \vec{\U}_1)$.\footnote{Technically, this is the $\Q$-structure for the last normal component of $\vec{\U}_1$, but we abuse notation here.} This means that 
\begin{center}
$\Q(b,\vec{\U}_1)\lhd Lp^{\Sigma,\Gamma}(\M(\vec{\U}_1))$. 
\end{center}
This follows from the fact that $\Q(b,\vec{\U}_1)$ must be iterable in $\Gamma$ for trees above $\delta$. But by $\Gamma$-fullness preservation of $\Lambda$, $\Q(c,\vec{\U}_1)$ exists and  $\Q(c,\vec{\U}_1) =  \Q(b,\vec{\U}_1)$. So $b=c$.

Suppose now $\delta$ is not a cutpoint of $\Q(b,\vec{\U}_1)$. Let $E$ be the least extender on the $\Q(b,\vec{\U}_1)$-sequence with the property that crt$(E) < \delta(\vec{\U}_1) < \textrm{lh}(E)$. Let $\U' = \vec{\U}_1^\smallfrown \langle E \rangle$. Let $\lambda = lh(\vec{\U}_1)$, $\kappa = \textrm{crt}(E)$, $\xi = \U'-pred(\lambda+1)$. Then it is easy to see that there is a $\gamma < o(\M^{\U'}_\xi)$ such that $$\M^{\U'}_\infty = \textrm{Ult}_n(\M^{\U'}_\xi|\gamma, E),$$ where $n$ is least such that $\rho_{n+1}(\M^{\U'}_\xi)\leq \kappa$.\footnote{This situation is what Sargsyan calls a ``fatal drop" in \cite{ATHM}. See also \cite{trang2013} for an alternative treatment and more details of such a situation.} By the minimality of $E$, we also have $$\M^{\U'}_\xi|\gamma \models ``\kappa \textrm{ is a limit of cutpoints}"$$ and $$\M^{\U'}_\infty \models ``\delta(\vec{\U}_1) \textrm{ is a cutpoint}''.$$ This implies $\M^{\U'}_\infty \lhd (Lp^\Sigma(\M(\U_1)))^{\Gamma(\R,\Lambda)}$ and since $\rho_{n+1}(\M^{\U'}_\infty) < \delta(\vec{\U}_1)$, we must have that $c$ drops; so $\Q(\vec{\U}_1,c)$ exists. Suppose $\Q(\vec{\U}_1,b)\neq \Q(\vec{\U}_1,c)$. By the argument in \cite[Claim 3.20]{trang2013}, letting $\mathcal{Y}, \mathcal{Z}$ be the results of comparing the phalanxes $\Phi(\vec{\U}_1^\smallfrown b), \Phi(\vec{\U}_1^\smallfrown c)$, then for every $\alpha\geq \lambda$, $[0,\alpha]_{\mathcal{Y}}, [0,\alpha]_{\mathcal{Z}}$ both drop. This gives a standard contradiction.\footnote{The fact that the last branches of $\mathcal{Y},\mathcal{Z}$ drop give that some pairs of extenders in $\mathcal{Y},\mathcal{Z}$ must be compatible. This contradicts the fact that $\mathcal{Y},\mathcal{Z}$ are comparison trees.} We note that the phalanx $\Phi(\vec{\U}_1^\smallfrown c)$ is iterable by the strategy induced by $\Lambda$.

We give the argument in \cite[Claim 3.20]{trang2013} here for the reader's convenience. Suppose not. Let $\alpha\geq \lambda$ be least such that either $F=E^{\mathcal{Y}}_\alpha$ or $F=E^{\mathcal{Z}}_\alpha$ overlaps $\delta$, i.e. crt$(F)<\delta < \textrm{lh}(F)$. Then $[0,\alpha']_{\mathcal{Y}}$ and $[0,\alpha']_{\mathcal{Z}}$ both drop for $\alpha'\in [\lambda,\alpha]$. Note that $\delta$ is Woodin in $M^{\mathcal{Y}}||\textrm{lh}(F)$ and if there is any $F'$ on the sequence of $M^{\mathcal{Y}}||\textrm{lh}(F)$ that overlaps $\delta$, then $[0,\beta]_{\mathcal{Y}}, [0,\beta]_{\mathcal{Y}}$ both drop for all $\beta > \alpha$. This is because Woodin cardinals are cutpoints of hod mice we consider (i.e. below ``$\sf{AD}_\mathbb{R} + $$\Theta$ is measurable").

Now we consider the case $F$ being the least extender overlapping $\delta$, and so  $\alpha = \lambda$. Let $\kappa'=\textrm{crt}(F)$ and $\epsilon$ be the least such that $F$ is applied to some $\Q\unlhd \M^{\mathcal{Y}}_\epsilon$ or $\Q\unlhd \M^{\mathcal{Z}}_\epsilon$ according to the rules of normal trees. Then $\mathcal{Y}\rest [\epsilon, \textrm{lh}(\mathcal{Y}))$ and $\mathcal{Z}\rest [\epsilon,\textrm{lh}(\mathcal{Z}))$ are equivalent to above-$\kappa'$, normal trees on $\Q$. If $\Q\lhd \M^{\vec{\U}_1}_\epsilon$, we are done. Otherwise, $[0,\epsilon]_{\vec{\U}_1}$ must drop because our hod mice are below ``$\sf{AD}_\mathbb{R} + $$\Theta$ is measurable" and $\kappa'$ is an inaccessible limit of Woodin cardinals. 

So $\Q(\vec{\U}_1,b) = \Q(\vec{\U}_1,c)$ and hence $b=c$.

\end{proof}

\begin{lemma}\label{claim:bcshc}
Suppose $(\R,\Lambda)$ is a reasonable hod pair, then $\Lambda$ has strong hull condensation. 
\end{lemma}
\begin{proof}
%
Suppose $\vec{\T}$ is according to $\Lambda$ and $\vec{\U}$ is a pseudo-hull of $\vec{\T}$. We assume for ease of notations in the following argument that $\vec{\U} = \vec{\U}_0^\smallfrown \vec{\U}_1$ and letting $\Q$ be the last model of $\vec{\U_0}$, then there is an ordinal $\beta$ such that:
\begin{itemize}
\item $\alpha(\vec{\U}_0)$, the set of generators used in $\vec{\U}_0$, is contained in $(\delta^\Q_\beta)^{<\omega}$.
\item $\vec{\U}_1$ is based on $\Q(\beta+1)$ and is above $\delta^\Q_\beta$.
\item $\vec{\U}_0$ is non-dropping and is according to $\Lambda$.
\end{itemize}
This is indeed the main case; the proof of other cases is similar and we will leave that to the reader.

In this case, we also have that $\vec{\T} = \vec{\T}_0^\smallfrown \vec{\T}_1$, where $\vec{\T}_0$ has last model $\S$ and the embedding $\varphi: \Q\rightarrow \S$ is the natural map. Hence, we have $\varphi\circ i_0 = j_0$ where $i_0$ is the iteration map given by $\vec{\U}_0$ and $j_0$ is the iteration map given by $\vec{\T}_0$. We also have that $\S$ is such that $\vec{\T}_1$ is above $\S(\varphi(\beta))$. Let $$b^*= \Lambda_{\vec{\T}_0,\S}(\vec{\T}_1)$$ and $$\S^* = \M^{\vec{\T}_1}_{b^*} .$$

Suppose the following holds.
\begin{equation}\label{eqn:pullbackequal}
\Lambda_{\vec{\T}_0,\S}^\varphi\rest \Q(\beta) = \Lambda_{\vec{\U}_0,\Q}\rest \Q(\beta).
\end{equation}
Call the strategy in \ref{eqn:pullbackequal} $\Sigma$. Let $\Psi = \Lambda_{\vec{\T}_0,\S}^\varphi$, $b = \Psi(\vec{\U}_1)$ and $c = \Lambda_{\vec{\U}_0,\Q}(\vec{\U}_1)$, we then show that 
\begin{center}
$b = c$.
\end{center}

\begin{figure}
\centering
\begin{tikzpicture}[node distance=2.5cm, auto]
  \node (A) {$\R$};
  \node (B) [right of=A] {$\Q$};
  \node (D) [node distance = 3cm, right of=B] {$\M^{\vec{\U_1}}_b$};
  \node (E) [node distance=2cm, below of=D] {$\M^{\vec{\U_1}}_c$};
  \node (G) [node distance=2cm,above of=A] {$\R$};
  \node (H) [right of=G] {$\S$};
  \node (I) [node distance=2cm, above of=D]{$\S^*$};
  \draw[->] (A) to node {$\vec{\U_0}, i_0$}(B);
  \draw[->] (A) to node {$id$}(G);
  \draw[->] (B) to node {$\varphi$}(H); 
  \draw[->] (D) to node {$\varphi^*$}(I);
  \draw[->] (B) to node {$\vec{\U_1},c$}(E);
  \draw[->] (B) to node {$\vec{\U_1},b$}(D);
  \draw[->] (G) to node {$\vec{\T}_0, j_0$}(H);
  \draw[->] (H) to node {$\vec{\T}_1,b^*$}(I);
  \end{tikzpicture}
\caption{Strong hull condensation.}
\label{fig:pullback}
\end{figure}
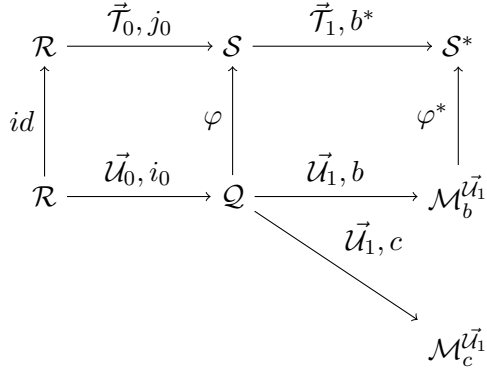

There are two cases. Suppose $b$ does not drop. Then there is a map $\varphi^*: \M^{\vec{\U}_1}_b \rightarrow \S^*$ given by the fact that $\vec{\U}$ is a pseudo-hull of $\vec{\T}$. We then have that $\pi^{\vec{\T}} = \varphi^*\circ \pi^{\vec{\U}_1}_b\circ i_0$. Applying branch condensation, we have that $b = c$. See Figure \ref{fig:pullback}.

Suppose $b$ drops, then $\Q(b,\vec{\U}_1)$ exists and the phalanx $\Phi(\vec{\U}_1^\smallfrown b)$ is iterable (above $\delta^\Q_\beta$) in $\Gamma$ because it is embeddable into the phalanx $\Phi(\vec{\T}_1^\smallfrown b^*)$ and by \cite[Lemma 4.20]{schlutzenberg2021iterability}.\footnote{In the case $\delta(\vec{\U}_1)$ is not a cutpoint of $\Q(b,\vec{\U}_1)$, as mentioned above, we simply have that  $\Q(b,\vec{\U}_1)$ is iterable above $\delta(\vec{\U}_1)$ as a $\Sigma$-mouse via a unique iteration strategy in $\Gamma$.} Lemma \ref{lem:uniqueQ} then implies that $b = c$. 


Now we prove equation \ref{eqn:pullbackequal}. Suppose not. Let $\vec{\W} = \vec{\W}_0^\smallfrown \vec{\W}_1$ be a minimal disagreement stack on $\Q(\beta)$. More precisely, $\vec{\W}_0$ is a nondropping stack on $\Q(\beta)$ according to both $\Lambda_{\vec{\T}_0,\S}^\varphi\rest \Q(\beta) =_{def} \Lambda^1$ and $\Lambda_{\vec{\U}_0,\Q}\rest \Q(\beta) =_{def} \Lambda^2$ with last model $\P^*$, $\alpha(\vec{\W}_0)\subseteq (\delta^{\P^*}_\gamma)^{<\omega}$ and $\vec{\W}_1$ is on $[\delta^{\P^*}_\gamma, \delta^{\P^*}_{\gamma+1})$ for some $\gamma$, and $\Lambda^1_{\vec{\W}_0,\P^*}(\vec{\W}_1)\neq \Lambda^2_{\vec{\W}_0,\P^*}(\vec{\W}_1)$.  Let $\psi:\P^* \rightarrow \R^*$ be the copy map from $\vec{\W}_0$ to $\varphi \vec{\W}_0$. Let $b_i = \Lambda^i_{\vec{\W}_0,\P^*}(\vec{\W}_1)$ for $i=1,2$. 

Again there are two cases just like above. If $b_1$ drops, then the same argument as above gives us $b_1 = b_2$. Now suppose $b_1$ does not drop. Let $\tau: \M^{\vec{\W}_1}_{b_1} \rightarrow \R^{**}$ obtained by copying $\vec{\W}_1^\smallfrown b_1$ to $\psi \vec{\W}_1^\smallfrown b_1$. Note that 
\begin{center}
$\tau\circ \pi^{\vec{\W}_1}_{b_1}\circ \pi^{\vec{\W}_0}\circ \pi^{\vec{\U}_0} = \pi^{\psi\vec{\W}_1}_{b_1}\circ \pi^{\varphi \vec{\W}_0} \circ \pi^{\vec{\T}_1}$.
\end{center}
By branch condensation, again, $b_1 = \Lambda_{\vec{\U}_0^\smallfrown \vec{\W}_0}(\vec{\W}_1)$. Therefore, $b_1 = b_2$. Contradiction. This shows Equation \ref{eqn:pullbackequal} holds and hence completes the proof of the lemma.

\end{proof}

The following lemma will be used in Lemma \ref{lem:pullback}. Lemma \ref{lem:pullback} also uses Lemma \ref{claim:bcshc} in an essential way. Lemma \ref{lem:pullback} may also be of independent interest and is used in an essential way in the proof of several theorems, including \ref{thm:sjs}, \ref{thm:omega1UB}.

\begin{lemma}\label{lem:omega2inaccessible}
Assume $\sf{CH} + $ there is an $\omega_1$-dense ideal $\mathcal{I}$ on $\omega_1$. Let $g\subseteq \mathbb{P}_{\mathcal{I}}$ be $V$-generic and $j = j_g : V\rightarrow M$ be the corresponding generic embedding. Suppose $(\R,\Lambda)$ is a reasonable hod pair where $\Lambda$ is an $(\omega_2,\omega_2)$-strategy. Suppose $A\subseteq \omega_1^V$ and $A$ codes $H_{\omega_1}^V$. Then in $X = L^{\Lambda}_{\omega_2^V}[A][g]$, there is no largest cardinal.
\end{lemma}
\begin{proof}
First, suppose $\pi: \P \rightarrow \R$ is elementary and $\P$ is countable. Let $\Psi = \Lambda^\pi$. $\Lambda$ has hull condensation, and hence $\Psi$ has hull condensation (see \cite{ATHM} for a proof that hull condensation ``pulls back"). We first claim that for any $x\in HC$ containing $\P$, 
\begin{center}
$L^\Psi_{\omega_2^V}[x] = j(L^\Psi_{\omega_1^V}[x]). \ \ \ \ \ \ \ \ \ \ \ \ \ \ \ \ \ \ \ \ \ \ \ \ \ \ \ (*)$
\end{center}
Suppose not. Then let $\T$ be a tree in $L^\Psi_{\omega_2^V}[x] \cap j(L^\Psi_{\omega_1^V}[x])$ such that $\Psi(\T) \neq j(\Psi)(\T)$. Let $\T$ be the least such (in the constructibility order of the models). Now the following are easy to see:
\begin{enumerate}[(a)]
\item $\T\in V$.
\item $j(\pi) = j\circ \pi$.
\item $j(\Psi) = j(\Lambda^\pi) = j(\Lambda)^{j\circ \pi}$.
\end{enumerate}
So 
\begin{equation}\label{eqn:pullback_equal}
j(\Psi)(\T) = j(\Lambda)^{j\circ \pi}(\T) = j(\Lambda)(j\circ \pi \T) = \Lambda(\pi \T) = \Psi(\T).
\end{equation}
The first equality follows from (c). The second and last equalities follow from definitions. To see the third equality, first note that by (a), $\pi\T\in V$ and therefore, $j\circ \pi \T$ is a hull of $j(\pi\T)$. Since $j(\pi\T)$ is according to $j(\Lambda)$, so is $j\circ \pi \T$ by hull condensation of $j(\Lambda)$.\footnote{We note that $j\circ \pi \T$ is countable in $V[g]$ and therefore is in $M$.} Now let $b = \Lambda(\pi\T)$, then $j(b) = j(\Lambda)(j(\pi\T))$ and $j\circ \pi\T^\smallfrown b$ is a hull of $j(\pi\T)^\smallfrown j(b)$. By hull condensation of $j(\Lambda)$, $b = j(\Lambda)(j\circ \pi \T)$ as desired. This is a contradiction. So $(*)$ holds.

$(*)$ implies that there is no $\alpha < \omega^V_1$ such that 
\begin{center}
$L^\Psi_{\omega_2^V}[x] \models \alpha^+ = \omega_1^V$.
\end{center}
This is because otherwise, in $j(L^\Psi_{\omega_2^V}[x]) \models \alpha^+ = j(\omega_1^V) = \omega_2^V$. This implies then that 
\begin{center}
$j(L^\Psi_{\omega_1^V}[x]) \models ``\omega_1^V$ is not a cardinal".
\end{center}
On the other hand, 
\begin{center}
$L^\Psi_{\omega_2^V}[x] \models ``\omega_1^V$ is a cardinal."
\end{center}
$(*)$ then immediately gives a contradiction. 

Now let $A\subseteq \omega_1^V$ and $A$ codes $H_{\omega_1}^V$. To see that there is no largest cardinal in $X$. It is enough to show there is no largest cardinal in $L^\Lambda_{\omega_2^V}[A]$.\footnote{Again, we use that $\mathbb{P}_{\mathcal{I}}$ is forcing equivalent to $Coll(\omega,\omega_1)$.} The argument above (showing Equations \ref{eqn:pullback_equal} hold) shows that 
\begin{center}
$\Lambda = j(\Lambda)^j\rest V$,
\end{center} 
and 
\begin{equation}\label{eqn:jLambdaj}
L^\Lambda_{\omega_2^V} = L^{j(\Lambda)^j}_{\omega^V_2}[A].
\end{equation}
Now, $\R$ is countable in $M$, $A\in HC^M$, and $j\rest \R: \R\rightarrow j(\R)$ is elementary in $M$, so the argument above, applied in $M$, shows that
\begin{center}
$\omega_1^M = \omega_2^V$ is not a successor cardinal in $L^{j(\Lambda)^j}_{\omega_2^V}[A] \ \ \ \ \ \ \ \ \ \ \ \ \ \ \ \ \  (**)$
\end{center}
$(**)$ and Equation \ref{eqn:jLambdaj} imply that there is no largest cardinal in $L^\Lambda_{\omega_2^V}[A]$.
\end{proof}

\begin{lemma}\label{lem:pullback}
Assume $\sf{CH} + $ there is an $\omega_1$-dense ideal $\mathcal{I}$ on $\omega_1$. Suppose $(\R,\Lambda)$ is a reasonable hod pair such that $|\R|^V\leq \omega_1$ and $\Lambda$ is an $\omega_2$-iteration strategy for $\R$. Let $g\subseteq \mathbb{P}_{\mathcal{I}}$ be $V$-generic and $j = j_g : V\rightarrow M$ be the corresponding generic embedding. Then $\Lambda = j(\Lambda)^j$.
\end{lemma}
\begin{proof}
By Lemma \ref{claim:bcshc}, $\Lambda$ has strong hull condensation. By strong hull condensation and \cite[Theorem 7.3]{schlutzenberg2021iterability}, there is a unique extension of $\Lambda$ in $V[g]$. Hence we identify $\Lambda$ with its canonical extension in $V[g]$.  First let $\T\in V$ be according to $\Lambda$. Then $j\T$ is a hull of $j(\T)$ and $j(\T)$ is according to $j(\Lambda)$, so $j\T$ is according to $j(\Lambda)$ by strong hull condensation of $j(\Lambda)$. But then $\T$ is by $j(\Lambda)^j$.

Suppose $\T\in M$ is according to $\Lambda$. Then there is a $\U\in V$ according to $\Lambda$ such that $\T$ is a pseudo-hull of $\U$ (see \cite[Theorem 7.3]{schlutzenberg2021iterability}); we note that to apply \cite[Theorem 7.3]{schlutzenberg2021iterability} to get the existence of $\U$, we need to work inside $X = L^\Lambda_{\omega_2^V}[tr.cl.(\{\dot{\T}\}\cup H_{\omega_1}^V)][g]$, where $\dot{T}\in H^V_{\omega_2}$ is a $Coll(\omega,\omega_1)$-name of $\T$. For \cite[Theorem 7.3]{schlutzenberg2021iterability} to apply, we need that $\omega_2^V > (\omega ^+)^X$. This follows from Lemma \ref{lem:omega2inaccessible}.

This means $j\T$ is a pseudo-hull of $j\U$\footnote{This fact can be easily verified, by chasing through the definition of pseudo-hull. See \cite{SteelNormalization}. Furthermore, \cite[Theorem 7.3]{schlutzenberg2021iterability} gives that if $\T$ is nondropping, then so is $\U$.} and $j\U$ is by $j(\Lambda)$ by the argument above. By strong hull condensation of $j(\Lambda)$, $j\T$ is by $j(\Lambda)$. Therefore, $\T$ is by $j(\Lambda)^j$. 

\end{proof}

\subsection{Boolean-valued comparison and $\sf{ZFC}$ comparison of hod pairs}\label{sec:booleanComp}

Suppose $(\P,\Sigma)$ is a reasonable hod pair such that $\Sigma$ is $\omega_1$-UB. Suppose $p\in Coll(\omega,\omega_1^V)$ and $G\subset Coll(\omega,\omega_1^V)$ is $V$-generic and $p\in G$; let let $g\subseteq \mathbb{P}_{\mathcal{I}}$ be the corresponding induced by $\G,\pi$ and $j_{g}:V\rightarrow M$ be the corresponding generic embedding.  Suppose $\Gamma\in M$ is an inductive-like pointclass. In cases of interest, $\Gamma$ is typically the largest Suslin pointclass in an $\sf{AD}^+$ model. For each $q\leq p$, let $G_p = G - G\rest dom(q) \cup q$ be the ``finite variation" of $G$ induced by $q$. Note that $V[G] = V[G_q]$ for all $q\leq p$; for each $q$, let $g_q\subseteq \mathbb{P}_{\mathcal{I}}$ be the corresponding induced by $\G_q,\pi$ and $j_{g_q}:V\rightarrow M_q$ be the corresponding generic embedding. Suppose $(\P_q, \Lambda_q)$ is a (countable) $\Sigma$-$\Gamma$-suitable mouse with $\Lambda_q$ being a $(\omega_1,\omega_1+1)$, $\Gamma$-fullness preserving strategy for $\P_q$ and $\Lambda_q$ is strongly guided by a sjs $\mathcal{A}_q$ that seals $\undertilde{Env(\Gamma)}$ (see Section \ref{sec:suc}). Then Woodin's Boolean comparison theorem (\cite{CMI}) gives us that we can compare $\{(\P_q,\Lambda_q) : q\leq p\}$ in $V[G]$ and the comparison results in a pair $(\R,\Lambda)$ such that $\R\in V$, $|\R|^V\leq \omega_1$, $\Lambda \rest H_{\omega_2}^V\in V$. Furthermore, $\Lambda$ is the tail of all the $\Lambda_q$'s via the iteration trees that appear in the comparison.

In our present context,\footnote{Another context, where the conditions for $\Lambda_q$'s below may not satisfy, occurs in the proof of Claim \ref{claim:independent}. We will show in that case the Boolean comparison still succeeds.} we only know $\Lambda_q$ is an $(\omega_1,\omega_1)$-iteration strategy in $V[G]$ for each $q$. However, we can still conclude the comparison above terminates in less than $\omega_1^{V[G]} = \omega_2^V$ many steps. This is because by $\Sigma_1$-reflection (inside the pointclass $j_g(\Gamma)$), we have that for every $q\leq p$, there is a countable tree $\T_q$ such that:
\begin{itemize}
\item $\T_q$ is correctly guided, i.e. whenever $\alpha < lh(\T_q)$ is limit, then $\Q(\T_q\rest \alpha)$ exists and $\Q(\T_q\rest \alpha)\lhd (Lp^\Sigma(\M(\T_q\rest\alpha)))^{j_g(\Gamma)}$.\footnote{We note that the fatal drop cases can be ruled out in the boolean comparison.}
\item $\T_q$ is maximal and has last model $M_q= (Lp^\Sigma(\M(\T_q)))^{j(\Gamma)}$.
\item $\{\T_q : q\leq p\}$ are obtained by the least-extender disagreement process.
\item For $q \neq r$, $M_q = M_r$.
\end{itemize}

The tree $\T_q$'s above are precisely the trees occurred during the Boolean comparison process. Notice we never referred to the strategies $\Lambda_q$ in the above process. $\Lambda_q$ is used to define $\Lambda_q(\T_q)$ at the end (i.e. picking the last, maximal branch of $\T_q$). This is possible because $\T_q$ is countable; that $\T_q$ is countable is a consequence of the fact that $\omega_1$ is measurable in $j_g(\Gamma)$. Therefore, the comparison process succeeds and results in $(\R,\Lambda)$ above.


We now introduce concepts needed for the proof of Claim \ref{claim:terminate}. In essence, the proof of Claim  \ref{claim:terminate} is a proof that a Boolean comparison between hod pairs $\{(\P_q,\Lambda_q) : q\in Coll(\omega,\omega_1^V)\}$\footnote{More generally, we compare pairs $(\P_q,\Lambda_q)$ for $q\leq p$, for some fixed condition $p$.} terminates in $V[G]$ (in less than $\omega_1$ many steps), where for each $q$, $\P_q$ is a hod mouse such that $\lambda^{\P_q}$ is a limit ordinal and $\Lambda_q$ is an $(\omega_1,\omega_1)$-strategy with branch condensation and for each $(\Q,\Psi) \in B(\P_q, \Lambda_q)$, $\Psi$ is a $(\omega_1,\omega_1+1)$-strategy and $\Psi\rest HC$ belongs to an $\sf{AD}^+$ model \footnote{See the definition of $\Gamma$ in the next section. In this paper, we will have that $(\Q,\Psi)\in j_g(\Gamma)$, even though $\Lambda_q$ need not belong to $j_g(\Gamma)$ a priori.}. Furthermore, we assume that for $p\neq q$, $(\P_q,\Lambda_q)$, $(\P_p,\Lambda_p)$ are hod pairs of the ``same kind" in that whenever $(\Q_1,\Psi_1)\in B(\P_q,\Lambda_q)$ and $(\Q_2,\Psi_2)\in B(\P_p,\Lambda_p)$, and suppose there is $\alpha < min(\lambda^{\Q_1},\lambda^{\Q_2})$ such that $(\Q_1(\alpha),(\Psi_1)_{\Q_1(\alpha)}) = (\Q_2(\alpha),(\Psi_2)_{\Q_2(\alpha)})$, then there are normal trees $\T_i$ according to $\Psi_i$ on the window $(\delta^{\Q_i}_\alpha,\delta^{\Q_i}_{\alpha+1})$ such that letting $\R_i$ be the end model of $\T_i$ and $\Lambda_i = (\Psi_i)_{\T_i,\R_i}$, then $(\R_1(\alpha+1),(\Lambda_1)_{\R_1(\alpha+1)})=(\R_2(\alpha+1),(\Lambda_2)_{\R_2(\alpha+1)})$.

Typically, $\Psi_1,\Psi_2$ are Suslin coSuslin in an $\sf{AD}^+$ model $M$ (e.g. $M$ is of the form $L(A,\mathbb{R})$ for $A\in \Gamma$). We let $(N,\delta,\Sigma)$ be a coarse $\Omega$-Woodin mouse for some inductive-like pointclass $\Omega\in M$ that contains all projective sets in $(\Psi_1,\Psi_2)$ and $(N,\delta,\Sigma)$ Suslin captures $\Psi_1,\Psi_2$.\footnote{See \cite{DMATM} for more details on coarse Woodin mice.} More precisely, $(N,\delta,\Sigma)$ has the following properties:
\begin{itemize}
\item $N\models \sf{ZFC}$.
\item $\delta$ is the unique Woodin cardinal of $N$.
\item $\Sigma$ is an iteration strategy for $N$.
\item $\Q_1,\Q_2\in N$.
\item For each $i\in\{1,2\}$, there are trees $(T_i, U_i) \in N$ that witnesses $(N,\delta,\Sigma)$ Suslin captures $\Psi_i$ at $\delta$, i.e. for any countable $\Sigma$-iterate $N'$ of $N$ such that there is an iteration map $i: N \rightarrow N'$, for any $h\subset Coll(\omega,i(\delta))$ such that $h\in V$ is $N'$-generic, $p[i(T_i)]\cap N'[h] = \Psi_i\cap N'[h]$ and $p[i(U_i)]\cap N'[h] = \mathbb{R}^{N'[h]} - \Psi_i$.\footnote{Here we fix a canonical coding of elements of $HC$ by reals identify $\Psi_i$ with its code.}
\end{itemize}

The existence of $\T_i$ is then easy to see. Let $\Lambda = (\Psi_1)_{\Q_1(\alpha)}=(\Psi_2)_{\Q_2(\alpha)}$. In $N$, iterate $(\Psi_1)_{\Q_1(\alpha+1)}$ and $(\Psi_2)_{\Q_2(\alpha+1)}$ into the $\Lambda$-hod mouse construction of $V^N_\delta$. Since these two strategies have branch condensation, there are normal trees $\T_i$ (as specified above) and iteration maps $k_i: \Q_i \rightarrow \R_i$ according to $\Psi_i$ such that 
\begin{enumerate}[(a)]
\item $(\R_1(\alpha+1),(\Lambda_1)_{\R_1(\alpha+1)})=(\R_2(\alpha+1),(\Lambda_2)_{\R_2(\alpha+1)})$.
\item $\R_1(\alpha+1)$ is model in the  $\Lambda$-hod mouse construction of $V^N_\delta$ and $(\Lambda_1)_{\R_i(\alpha+1)}$ is the background induced strategy.
\end{enumerate}
See \cite{ATHM} for more details. The above argument generalizes easily to countably many hod pairs (as in the proof of Claim \ref{claim:terminate}).

The comparisons described above are the building blocks of the ``diamond comparison" described in Claim \ref{claim:terminate}. The ``diamond comparison" of all pairs of the form $(\P_q,\Lambda_q)$ for $q\in Coll(\omega,\omega_1)$ must end in $<\omega_1$ steps in $V[G]$; see the proof of Claim \ref{claim:terminate} for more details.


\section{OUTLINE OF THE PROOF OF THEOREM \ref{thm:DI}}\label{sec:outline}
We outline the proof of Theorem \ref{thm:DI}.  In $V$, define the maximal pointclass 
\begin{center}
$\Gamma= \{A\subseteq \mathbb{R} : L(A,\mathbb{R})\models \sf{AD}^+\}$.
\end{center}
The goal is to show that $\Gamma$ is sufficiently rich in that there is a $\Omega\subseteq \Gamma$ such that $L(\Omega,\mathbb{R})\models \sf{AD}_\mathbb{R} + $$ \Theta$ is regular. So suppose not. We assume:
\begin{center}
$(\ddag): \ $ No $\sf{AD}^+$ models satisfy ``$\sf{AD}_\mathbb{R} + $$\Theta$ is regular."
\end{center}

As part of the induction, we maintain:

\begin{center}
$(\dag)$: All cmi operators $J$ are $\omega_1$-UB.
\end{center}

We will analyze the complexity of $\Gamma$, ultimately showing that there is some Wadge initial segment $\Omega$ of $\Gamma$ (possibly $\Omega = \Gamma$) such that $L(\Omega, \mathbb{R})\models \sf{AD}_\mathbb{R} + {}$``$\Theta$ is regular.'' There are two major cases. We summarize the key points of each case below before jumping into the details.
\begin{enumerate}[(i)]
\item The successor case (Section \ref{sec:suc}): we first show that if $(\P,\Sigma)\in \Gamma$ ($\Sigma$ may be $\emptyset$) is a hod pair such that $\Sigma$ is $\Gamma$-fullness preserving and has branch condensation, then Lp$^\Sigma(\mathbb{R})\models \sf{AD}^+$, and therefore $\powerset(\mathbb{R})\cap \textrm{Lp}^\Sigma(\mathbb{R})\subset \Gamma$. This is via a standard core model induction argument similar to that showing $\sf{AD}$ holds in $L(\mathbb{R})$ (\cite{CMI, wilson2012contributions}). One wrinkle that appears in the case that $\Sigma\neq \emptyset$ is that one needs to show $\M_1^{\Sigma,\sharp}$ exists before being able to define Lp$^\Sigma(\mathbb{R})$ as done in \cite{trang2013}. The argument showing that $\M_1^{\Sigma,\sharp}$ exists is given in Theorem \ref{thm:m1sharp}. 

As part of the induction, we maintain $(\dag)$, the hypothesis that for every $\Sigma$-cmi operator $J$ (including the operator induced by $\Sigma$), $J$ is $\omega_1$-UB. This is what we need to carry out the proof of Theorem \ref{thm:m1sharp}. This then allows us to adapt the standard arguments in \cite{CMI, wilson2012contributions} to show Lp$^\Sigma(\mathbb{R})\models \sf{AD}^+$. 

In Section \ref{sec:suc} (see in particular Theorem \ref{thm:sjs}), we adapt the argument in \cite{wilson2012contributions} to show that there is a self-justifying system $\mathcal{A}$ consisting of sets Wadge cofinal in Lp$^\Sigma(\mathbb{R})$, and a $\Sigma$-suitable pair $(\Q,\Lambda)$ where $\Lambda$ is the strategy guided by $\mathcal{A}$.\footnote{This argument allows us to construct $(\Q,\Lambda)$ without the technical hypothesis $\sf{HI}(c)$ in Ketchersid's thesis. See \cite{CMI,ketchersid2000toward} for an alternative argument constructing $(\Q,\Lambda)$ that uses a seemingly stronger hypothesis.} Therefore, $\Lambda$ is $\Gamma$-fullness preserving and has branch condensation and $\Lambda\notin$ Lp$^\Sigma(\mathbb{R})$. 

We can then show Lp$^\Lambda(\mathbb{R})\models \sf{AD}^+$ and therefore $\Lambda\in \Gamma$. To do this, we first need to show some such $\Lambda$ can be extended to an $\omega_2$-strategy in $V$ and is $\omega_1$-UB (Theorem \ref{thm:omega1UB}). Crucially, we use Lemma \ref{lem:pullback} in this argument. 

\item The limit case (Section \ref{sec:lim}): assuming $(\ddag)$ and letting $\mathcal{H},\mathcal{H}^+$ and $\Sigma$ be defined as in Section \ref{sec:lim}, we use the generic embedding $j:V\rightarrow M$ induced by a $V$-generic $G\subset \textrm{Coll}(\omega,\omega_1)$ to derive a nice strategy $\Lambda$ for $\mathcal{H}^+$ in $M$. The strategy $\Lambda$ is $j(\Gamma)$-fullness preserving, has branch condensation, and most importantly, if $\Gamma(\mathcal{H}^+,\Lambda)\subsetneq j(\Gamma)$, then letting $\M_\infty(\mathcal{H}^+,\Lambda)$ be the direct limit of non-dropping iterates of $(\mathcal{H}^+,\Lambda)$ in $j(\Gamma)$, we have $\M_\infty(\mathcal{H}^+,\Lambda) = \mathcal{H}(\delta)$ where $\delta =  \delta^{\M_\infty(\mathcal{H}^+,\Lambda)}$, and there is a factor map $\sigma: \M_\infty(\mathcal{H}^+,\Lambda)\rightarrow j(\mathcal{H}^+)$ such that crt$(\sigma) = \delta$. This property is a consequence of the \textit{$j$-condensation lemma}, Theorem \ref{thm:cond}. This result is crucial here and its variations are important in many other arguments (cf. \cite{sargsyanCovering2013, LSA, Trang2015PFA}).

Again, Lemma \ref{lem:pullback} will be useful in proving Theorem \ref{thm:cond} and Lemma \ref{fullness}. Part of the proof of Lemma \ref{fullness} is to show that $j\rest o(\mathcal{H}^+)$ is continuous. This continuity property is also important in the proof of Theorem \ref{thm:cond}.

Now there are two cases. Suppose first that $\Gamma(\mathcal{H}^+,\Lambda) = j(\Gamma)$. Then by elementarity, in $V$ there is a hod pair $(\P,\Sigma)$ such that $\Gamma(\P,\Sigma) = \Gamma$; in particular, $\Sigma\notin \Gamma$. By a core model induction as in the successor case, Lp$^\Sigma(\mathbb{R})\models \sf{AD}^+$. To show this, we again have to show we can extend $\Sigma$ to $H_{\omega_2}^V$ and that $\Sigma$ is $\omega_1$-UB (see Lemma \ref{lem:limitUB}). This implies $\Sigma\in \Gamma$. Contradiction. Otherwise, $\Gamma(\mathcal{H}^+,\Lambda) \subsetneq j(\Gamma)$. Therefore $\sigma$ exists and $\delta$ is a regular cardinal which is a limit of Woodin cardinals in $\M_\infty(\mathcal{H}^+,\Lambda)$. By standard arguments, $L(j(\Gamma)\rest \delta, \mathbb{R}^M)\models \sf{AD}_\mathbb{R}+ {}$``$\Theta$ is regular.'' This is again a contradiction, so $(\ddag)$ fails. This completes the outline of the proof.
\end{enumerate}

\section{SUCCESSOR STEP}\label{sec:suc}

Suppose $(\P,\Sigma)\in \Omega$ is a reasonable hod pair such that $\Sigma$ is $\Omega$-fullness preserving, has branch condensation, and $\Sigma$ is $\omega_1$-UB (i.e. we assume the hypothesis $(\dag)$ holds for $\Sigma$). This includes the case $(\P,\Sigma) = (\emptyset, \emptyset)$. We show that $Lp^\Sigma(\mathbb{R})\models \sf{AD}^+$. See Sections \ref{sec:strat_premice}, \ref{sec:cmioperators} for a summary of mouse operators and $\Sigma$-cmi operators and related concepts (like the definition of $\mathcal{F}_{\Sigma,\varphi}$).


\begin{theorem}\label{thm:m1sharp}
Suppose $F$ is a nice mouse operator (or a $\Sigma$-cmi operator) on $H_{\omega_1}^V$ that is $\omega_1$-UB, then $\M_1^{F,\sharp}$ is a nice operator (or a $\Sigma$-cmi operator) and is $\omega_1$-UB.
\end{theorem}
\begin{proof}
We assume that $F$ is a nice $\Sigma$-cmi operator where $\Sigma\in \Omega$ has branch condensation, is $\Omega$-fullness preserving, and is $\omega_1$-UB. Without loss of generality, we assume $\mathcal{F} = \mathcal{F}_{\Sigma,\varphi}$ be the operator induced by $\Sigma$ and with $\varphi = \varphi_{\textrm{all}}$ and $\mathcal{F}^+$ the canonical extension of $\mathcal{F}$ in $V[g]$. The case $F\neq\mathcal{F}_{\Sigma,\varphi}$ is similar. The operator $\mathcal{F}$ codes up the same information as $\Sigma$ does; the reader will lose little by pretending $\mathcal{F}=\Sigma$.

The proof that $\mathcal{F}^\sharp$ exists and is $\omega_1$-UB is standard. Details have been given in \cite{wilson2012contributions, CMI}. We only mention some key points here. The operator $\mathcal{F}^+$, the unique extension of $\mathcal{F}$ in $V$, is simply $j(\mathcal{F})\rest V$. Since $\mathcal{F}$ satisfies ($\dag$), the hypothesis $\sf{DI}$ will imply that $j(\mathcal{F})\rest V$ is in $V$ and doesn't depend on $G$; this follows from Lemma \ref{lem:UBprojection} and homogeneity of the forcing $Coll(\omega,\omega_1)$. We will write $\mathcal{F}$ for $\mathcal{F}^+$ for brevity.

To see $\mathcal{F}^\sharp(x)$ is defined for each $x\in dom(\mathcal{F})$, note that from $j$, one can define an ultrafilter $\mu$ over $L^\mathcal{F}[x]$\footnote{This is the model $L_{\mathfrak{c}^+}^{\mathcal{F}}[x]$.} as follows: for each $A\in \powerset(\omega_1^V)\cap L^\mathcal{F}[x]$, 
\begin{center}
$A\in \mu \iff \omega_1\in j(A)$.
\end{center}
By a standard argument, $\mu$ is a countably complete, normal measure over $L^\mathcal{F}[x]$ that is amenable to $L^\mathcal{F|}[x]$ in the sense that for any $Y$ of size $\omega_1^V$ in $L^\mathcal{F}[x]$, we have $\mu\cap Y \in L^\mathcal{F}[x]$. Furthermore, by condensation properties of $\mathcal{F}$ we have Ult$(L^\mathcal{F}[x], \mu) = L^\mathcal{F}[x]$ as it embeds into $j(L^\mathcal{F}[x])$. By standard arguments due to Kunen, the amenable structure $(L^\mathcal{F}[x],\mu)$ is iterable. This implies $\mathcal{F}^\sharp(x)$ exists.

To prove $\M_1^{\mathcal{F},\sharp}$ exists, we need to build the $K^{c,\mathcal{F}}$-construction inside $N = L^{\mathcal{F}^\sharp}(\mathbb{R})$ and run the proof of \cite[Theorem 2.10.2]{CMI}. For contradiction, we get for some $x\in\mathbb{R}$, the core model (relative to $\mathcal{F}$) $K = K^\mathcal{F}(x)$ exists (and iterable) in $N$. We need that $j(K)\in V$. To show this, we need to show $j(N)$ is definable in $V[g]$ from parameters in $V$. Here are some details that execute this plan.

We define the following model $W$ by induction on $\alpha < \omega_2^V$: $W_0 = (HC^V,\in)$, 
\begin{center}
$W_{\alpha+1} = J_\omega(tr.cl. (W_\alpha\cup \{(\T,b): b = \Sigma(\T) \wedge \T\in W_\alpha\wedge \T \textrm{ is according to } \Sigma\}))$, \footnote{Equivalently, $W_{\alpha+1} = J_\omega(tr.cl.(W_\alpha\cup \{(x,\mathcal{F}(x)): x\in W_\alpha\}))$.}
\end{center}
and for $\alpha$ limit, $W_\alpha = \bigcup_{\beta<\alpha} W_\beta$. Finally, $W = \bigcup_{\alpha<\omega_2^V} W_\alpha$. Note that $W\in V$ and $\Sigma\rest W_\alpha\in W$ for all $\alpha<\omega_2^V$.

By the proof of \cite[Lemma 3.35]{ATHM}, we have the following.

\begin{fact}\label{claim:det_gen}
For any poset $\mathbb{P}\in W$ and any $W$-generic $g\subset \mathbb{P}$ such that $g\in V$ (or $g\in M$), $W[g]$ is closed under $\Sigma$ (respectively $j(\Sigma)$).\footnote{\cite[Lemma 3.35]{ATHM} indeed implies that $\mathcal{F}$ determines itself on generic extensions. It is also easy to see that $\mathcal{F}$ relativizes well.}
\end{fact}

Let $\Sigma^+$  be the canonical extension of $\Sigma$ to $V[g]$. We fix trees $T,U\in V$ witnessing $\Sigma$ is $\omega_1$-UB. So in $V[g]$, $Code(\Sigma^+) = p[T] = \mathbb{R}\backslash p[U]$. Note also that $j(\Sigma)\rest V = \Sigma^+\rest V$. Suppose $h\in V$ (or in $M$) is a generic enumeration of $\mathbb{R}^V$ in order type $\omega^V_1$, let $X_h = \bigcup_{\alpha<\omega_2^V} X_\alpha$, where $X_0 = tr.cl.(h\cup\{h\})$, $X_1 = \mathcal{F}_0^+(X_0)$, and for $\alpha\geq 1$, $X_{\alpha+1}=\mathcal{F}_1^+(X_\alpha)$, and $X_\lambda = \bigcup_{\alpha<\lambda} X_\alpha$ for $\lambda$ a limit ordinal; here $\mathcal{F}^+$ codes $\Sigma^+$ the same way $\mathcal{F}$ codes $\Sigma$ and see \cite{trang2013} for the precise definition of $\mathcal{F}_0^+, \mathcal{F}_1^+$.\footnote{The reader will lose little by pretending $\mathcal{F}^+ = \mathcal{F}^+_0 = \mathcal{F}^+_1$.} We note that $X_h$ contains $\mathbb{R}^V$ and is closed under $\Sigma^+$. Now, if we let $W^{X_h}$ be the structure $W$ defined as above, but the definition is carried out inside $X_h$, then 
\begin{center}
$W^{X_h} = W$.
\end{center}
This means that the model $W$ is independent of $h$.

In a similar manner, letting $\mathcal{G} = (\mathcal{F}^+)^\sharp$, we define $X_h = \bigcup_{\alpha<\omega_2^V} X_\alpha$, where $X_0 = tr.cl.(h\cup\{h\})$, $X_1 = \mathcal{G}_0(X_0)$, and for $\alpha\geq 1$, $X_{\alpha+1}=\mathcal{G}_1(X_\alpha)$, and $X_\lambda = \bigcup_{\alpha<\lambda} X_\alpha$ for $\lambda$ a limit ordinal\footnote{$X_g$ is a potential $\mathcal{G}$-premouse over $g$ and it is closed under $\mathcal{G}$ because $\mathcal{G}$ relativizes well.} and let $W$ be the model defined in $X_h$ as above, but using $\mathcal{G}$ instead of $\mathcal{F}$. In particular, it is easy to verify that $W$ now has the following properties:
\begin{itemize}
\item $o(W) = \omega_2$ and $W$ is a transitive model over $\mathbb{R}^V$;
\item for any $a\in H_{\omega_2}\cap W$, we have $\mathcal{G}(a)\in W$; in particular, $W$ is closed under $\Sigma^+$ and if $h\in V$ (or in $M$) is $W$-generic, then $W[h]$ is closed under $\Sigma^+$;
\item $W$ is independent of $h$; in other words, suppose $h_1,h_2\in V$ (or in $M$) are two enumerations of $\mathbb{R}$ in order type $\omega_1$, then $W^{X_{h_1}} = W^{X_{h_2}}$.
\item If $h\in V$ (or in $M$) is Coll$(\omega_1,\mathbb{R}^V)$-generic over $W$, then the universe of $X_h$ is just the universe of $W[h]$.
\end{itemize}

Suppose that on a cone of $x\in HC$,  $\M_1^{\mathcal{F},\sharp}(x)$ does not exist.  Then in $W[h]$ where $h\in V$ is Coll$(\omega_1,\mathbb{R}^V)$-generic over $W$, the core model $K =_{\textrm{def}} K^{\mathcal{F}^+}(x)$ exists\footnote{Here the core model relative to $\Sigma^+$ is defined in the sense of \cite{jensen2013k} and $o(K) =\xi < o(W)$ and $\omega_1^V < \xi$ is a sufficiently large indiscernible relative to $\mathcal{G}$.}. Here $K$ is a $\mathcal{F}$-mouse and is in $W$. 

\begin{claim}\label{claim:j(K)inV}
$j(K)\in V$. 
\end{claim}
\begin{proof}[Proof of Claim \ref{claim:j(K)inV}]
To see that $j(K)\in V$, it suffices to show that $j(W)$ is definable in $V[G]$ from parameters in $V$. To see this, first note that $j(j(\Sigma)\rest V)$ is definable in Ult$(V,g)$ as the unique extension of $j(\Sigma)$ to $H_{\omega_2}$ that condenses well. Let $\Lambda = j(j(\Sigma)\rest V)$. Note that $\Lambda$ may not be definable in $V[G]$; the main wrinkle is that $H_{\omega_2}^{V[G]}$ may differ from $H_{\omega_2}^{\textrm{Ult}(V,g)}$. But in $V[G]$ we can define $\Psi$, the maximal (partial) strategy on $H_{\omega_2}$ that extends $j(\Sigma)$ with the property that whenever $\T$ is according to $\Psi$, the branch $\Psi(\T)$ (if defined) is the unique well-founded branch $b$ such that whenever $(\U,c)$ is a countable hull of $(\T,b)$, then $\U$ is according to $j(\Sigma)$ and $c = j(\Sigma)(\U)$. Note that if $\Lambda_1$ and $\Lambda_2$ are partial strategies extending $j(\Sigma)$ satisfying the above properties, then whenever $\T \in dom(\Lambda_1)\cap dom(\Lambda_2)$, we have $\Lambda_1(\T) = \Lambda_2(\T)$. As a result, $\Psi$ is simply the union of all such partial strategies, and since $\Lambda$ is one such partial strategy, 
\begin{center}
$\Lambda \subseteq \Psi$.
\end{center}
This easily implies that $j(W)$ is definable in $V[G]$ from $\Psi$ as $\Psi\rest W = \Lambda \rest W$. Hence $j(W)$ is definable in $V[G]$ from $j(\Sigma)$, but $j(\Sigma) = p[j(T)]\cap V^{Ult(V,g)} =  p[T]\cap V[g]$ (by Lemma \ref{lem:UBprojection}); so $j(W)$ is definable in $V[G]$ from $(T,U)$.\footnote{This is the crucial point and is the reason we maintain that operators we construct in this core model induction are $\omega_1$-UB.} By homogeneity, $j(W)\in V$.
\end{proof}

Given this claim, the rest of the proof proceeds as in  \cite[Theorem 2.10.2]{CMI} by showing that for the $(\omega_1^V,\omega_2^V)$-extender $E$ derived from $j$, we have $E\rest \alpha\in j(K)$ for all $\alpha<\omega_2^V$. This implies that $\omega_1^V$ is Shelah in $j(K)$, contradiction. Fixing $\alpha < \omega_2^V$, we give a sketch of $E\rest \alpha\in j(K)$. We note again that $W$ is closed under $\mathcal{G}$. We need to see that the phalanx $(j(K), \textrm{Ult}(j(K), E\rest \alpha), \alpha)$ is iterable in $j(W)$.\footnote{Iterability here is with respect to trees of length $< j(\xi)$ in $j(W)$.} Otherwise in $j(W)$ there is a countable $\mathcal{F}$-premouse $\bar{K}$ and a map $\sigma: \bar{K}\rightarrow \textrm{Ult}(K,E\rest \alpha)$ with crt$(\sigma)=\alpha$ and
\begin{center}
$j(W) \vDash (j(K), \bar{K},\alpha)$ is not $\omega_1$-iterable.
\end{center}
We have a factor map $k: \textrm{Ult}(K,E\rest \alpha)\rightarrow j(j(K))$ with $k\rest \alpha = id$ and 
\begin{center}
$k\circ \sigma: \bar{K}\rightarrow j(j(K))$
\end{center}
such that $k\circ \sigma\rest \alpha = id$. Note that $j(j(K))$ makes sense by the claim above.

Let $\psi = k\circ \sigma$ and $\psi = [\beta \mapsto \psi_\beta]_G$. Let $\bar{K} = [\beta\mapsto K_\beta]_G$ and $\alpha = [\beta \mapsto \alpha_\beta]_G$. We need to see that for $G$-almost all $\beta$, 
\begin{center}
$W\vDash (K, K_\beta,\alpha_\beta)$ is $\omega_1$-iterable.
\end{center}
By absoluteness, in $j(W)$ there is some $\psi'_\beta: K_\beta\rightarrow j(K)$ such that $\psi'_\beta\rest \alpha_\beta = id$. Then in $W$ there is some $\bar{\psi}: K_\beta\rightarrow K$ such that $\bar{\psi}\rest \alpha_\beta = id$. But this means $(K, K_\beta,\alpha_\beta)$ is iterable in $W$. We have reached a contradiction.

Finally, the operator $H: x\mapsto \M_1^{\mathcal{F},\sharp}(x)$ is definable from $\Sigma$. Since $j(\Sigma)\rest V\in V$, we have $j(H)\rest V\in V$ also. It is then standard to show $H$ is $\omega_1$-UB. One shows that for club many countable $X \prec (H_{\omega_2},\in, (T,U))$, $X$ is generically correct about $H$, namely letting $\pi_X: M_X\rightarrow X$ be the uncollapse map, for any forcing $\mathbb{P}\in M_X$ such that $$M_X \models ``|\mathbb{P}|\leq \omega_1",$$ for any $M_X$-generic $g\subset\mathbb{P}$ such that $g\in V$, 
then for any $x\in HC\cap M_X[g]$,

\begin{center}
$V\models \varphi[x, (T,U)] \Leftrightarrow M_X[g] \models  \varphi[x, \pi_X^{-1}(T,U)]$,
\end{center}
here $\varphi(x, (T,U))$ is the natural formula that defines $H(x)$ from $\mathcal{F}$. We give an informal definition of $\varphi(x, y)$ here. $\varphi(x, y)$ is the statement: there is a unique $z$ such that
\begin{enumerate}[(a)]
\item $z$ has the first order properties of $\M_1^{\mathcal{F},\sharp}(x)$, where $Code(\mathcal{F}) = p[(y)_0] = \mathbb{R} - p[(y)_1]$.
\item $z$ has a unique $(\omega_1,\omega_1+1)$-iteration strategy $\Lambda$ with the property that whenever $\T$ is according to $\Lambda$ with limit length (or $\T$ is a stack with last normal component with limit length), $\Lambda(\T)$ is the unique $b$ such that $\M^\T_b \unlhd \mathcal{F}^\sharp(\M(\T))$.
\end{enumerate}

\end{proof}

The induction through Lp$^\Sigma(\mathbb{R})$ proceeds as usual and is organized by the scales pattern in Lp$^\Sigma(\mathbb{R})$ (see \cite{trang2013}). The above theorem takes care of the successor steps in the induction in Lp$^\Sigma(\mathbb{R})$.  The limit step is non-trivial and requires the use of our hypothesis when we reach an inductive-like $\Gamma$; recall here that a pointclass $\Gamma$ is \textit{inductive-like} if it is $\omega$-parametrized, closed under $\forall^\mathbb{R}, \exists^\mathbb{R}$, recursive substitution, and has the scale property. We need to construct an operator that is beyond $\undertilde{Env(\Gamma)}$ to continue the induction.\footnote{In fact, we need the hypothesis in the construction of the ``next" operator when $\Gamma$ is the last scaled pointclass in Lp$^\Sigma(\mathbb{R})$.} We start with a useful lemma.

\begin{lemma}\label{lem:UB}
Suppose $(\P,\Sigma)$ is a reasonable hod pair such that $\Sigma$ is $\omega_1$-UB. Suppose $(\P',\Sigma')$ is a pair such that $\P'$ is a countable $\Sigma$-premouse that is $\Gamma$-suitable for some inductive-like pointclass $\Gamma\subset Lp^\Sigma(\mathbb{R})$ and $\Sigma'$ has branch condensation and is a $\Gamma$-fullness preserving strategy for $\P'$ (as a $\Sigma$-mouse) that can be uniquely extended to an $(\omega_2,\omega_2)$-strategy, then $\Sigma'$ is $\omega_1$-UB.
\end{lemma}
\begin{proof}
We identify $\Sigma'$ with its unique extension to stacks in $H_{\omega_2}^V$. Let $i: \P'\rightarrow \M_\infty$ be the direct limit map of all non-dropping $\Sigma'$-iterates via stacks in $H^V_{\omega_2}$. For a club of countable $Y\prec (H_{\omega_3},\in, (\P',\Sigma'), i, \M_\infty)$, let $\pi_Y: M_Y\rightarrow Y$ be the uncollapse, let $\kappa_Y = crt(\pi_Y)$, and let $a^Y =\pi_Y^{-1}(a)$ for any $a\in Y$. Let $h\subseteq Coll(\omega,\kappa_Y)$ be a $M_Y$-generic in $V$. Let $\T,b\in M_Y[h]$, let $a=(i, \M_\infty)$ where $\T$ is a normal tree, and  let $\varphi(\T, b, a^Y)$ say: 
\begin{itemize}
\item $\T$ is correctly guided i.e. all strict initial segments of $\T$ are given by the $\Q$-structures in $C_{j(\Gamma)}$.\footnote{At this point, we know $C_{j(\Gamma)}$ is independent of generics $g$. To see this, suppose $g_1,g_2$ are such that leting $j_i: V\rightarrow M_i$ be the corresponding generic embeddings, and $Lp^{\Sigma,j_1(\Gamma)}(\M(\T))\lhd Lp^{\Sigma,j_2(\Gamma)}(\M(\T))$. Let $\M$ be the least in $Lp^{\Sigma,j_2(\Gamma)}(\M(\T))- Lp^{\Sigma,j_1(\Gamma)}(\M(\T))$ and $\Lambda_\M$ be its unique strategy. Note that $\Lambda_\M\in j_2(\Gamma)$; therefore, $\Lambda_\M\rest V\in \Gamma$ as $j_2(\Lambda_\M\rest V) = \Lambda_\M$. This means $j_1(\Lambda_\M)\in j_1(\Gamma)$. This contradicts the choice of $\M$.}
\item If $\T$ is short then $b$ is the unique cofinal branch such that $\Q(b,\T)$ exists and the phalanx $\Phi(\T^\smallfrown  \Q(b,\T))$ is iterable with unique strategy in $C_{j(\Gamma)}(\M(\T))$.
\item If $\T$ is maximal then $b$ is the unique non-dropping branch such that there is a map $\sigma: \M^\T_b\rightarrow\M_\infty^Y$ such that $i^Y = \sigma\circ i^\T_b$.
\end{itemize}

We need to see that $M_Y[h] \models \varphi(\T,b,a^Y)$ if and only if $\Sigma'(\T) = b$. Suppose first $\T$ is short. Note that $H_{\omega_2}^{M_Y}$ is closed under $\Sigma'$, so we let $W\in H_{\omega_2}^{M_Y}$ be transitive such that $\T\in W[h]$. Let $\pi: \P\rightarrow \Q$ be the iteration map given by the generic genericity iteration according to $\Sigma'\rest H_{\omega_2}^{M_Y}$ that makes $W$ generically generic. So $\T\in \Q[W,h]$ and $C_{j(\Gamma)}(\M(\T))\in \Q[W,h]$ by $j(\Gamma)$-fullness of $\Q$ and the fact that the operator $C_{j(\Gamma)}$ relativizes well. Therefore, $\Q(b,\T) = \Q(\T) \in M_Y[h]$. This shows that $M_Y[h]$ is correct about the shortness of $\T$ and can compute the correct $\Q$-structure and hence the branch $\Sigma'(\T)$.

Suppose $\T$ is maximal. The above calculation shows that this is equivalent to $\T$ being maximal in $\M_Y[h]$. If $M_Y[h] \models \varphi(\T,b,a^Y)$ then there is a $\sigma: \M^\T_b \rightarrow \M^Y_\infty$ such that $i^Y = \sigma\circ i^\T_b$. In $V$, let $\vec{\T}$ be according to $\Sigma'$ with last model $\M_\infty^Y$ such that $i^Y = i^\T$. Then by branch condensation of $\Sigma'$, $b = \Sigma'(\T)$. Conversely, suppose $\Sigma'(\T)= b$. Let $c = j_h(\Sigma' \rest M_Y)(\T)$. Then by boolean comparisons, it is easy to see there is a $\sigma: \M^\T_c\rightarrow \M^Y_\infty$ such that $\sigma\circ i^\T_c = i^Y$. But $i^Y$ is an iteration map according to $\Sigma'$ (in $V$), by branch condensation of $\Sigma'$, $c = \Sigma'(\T)$. So $b= c$.

The argument for stacks is similar. We leave the details to the reader. This completes the proof of the lemma.
 
\end{proof}

\begin{theorem}\label{thm:sjs}
Suppose $\Gamma^*\subset \textrm{Lp}^\Sigma(\mathbb{R})$ is an inductive-like pointclass such that:
\begin{itemize}
\item $\Gamma^* \models \sf{AD}^+$, and
\item $\Gamma^*$-$\sf{MC}$$(\Sigma)$ holds.
\end{itemize}
Then 
\begin{enumerate}[(a)]
\item for any $A\in \undertilde{Env(\Gamma^*)}$, there is a scale on $A$ whose norms are in $\undertilde{Env(\Gamma^*)}$;
\item there is a self-justifying system (sjs) $(A_i : i<\omega)$ sealing $\undertilde{Env(\Gamma^*)}$.
\end{enumerate}
\end{theorem}
\begin{proof}

We assume for simplicity that $\Sigma=\emptyset$, so $Lp^\Sigma(\mathbb{R}) = Lp(\mathbb{R})$; the general case is just more notationally complicated. We assume $\Gamma^* = \Sigma_1^{Lp(\mathbb{R})}$, i.e. the largest scaled pointclass in $Lp(\mathbb{R})$. The other cases are taken care of by the scales analysis in $Lp(\mathbb{R})$ (see \cite{CMI, schlutzenberg2016scales, Scalesendgap}). Let $T$ be the tree of a $\Gamma$-scale on a universal $\Gamma^*$ set; $T$ is a tree on $\omega\times \kappa$, where $\kappa$ is the largest Suslin cardinal of $P = Lp(\mathbb{R})$.\footnote{The following argument works for $P = Lp^\Sigma(\mathbb{R})$. One just needs to put the trees $T,U$ witnessing $\Sigma$ is $\omega_1$-UB into the parameters that define all the relevant objects below.} Let $A = p[T]$ be the universal $\Gamma^*$-set induced by $T$. We note that at this point, we know that $P \models \sf{AD}^+$; this is because by essentially the Kechris-Woodin transfer theorem (see \cite{wilson2012contributions} for more discussions in this particular context),  $\sf{AD}^+$ holds for sets in $\undertilde{Env(\Gamma^*)}$ and $\undertilde{Env(\Gamma^*)} = \powerset(\mathbb{R})\cap P$ by arguments in \cite[Lemma 4.5.1]{wilson2012contributions}. We assume for contradiction that (a) (and hence (b)) fails.

\begin{claim}\label{claim:independent}
For any $V$-generic $g\subset \mathbb{P}_\mathcal{I}$, suppose $j_g: V\rightarrow \textrm{Ult}(V,g)= M$ is the associated ultrapower map and $G\subseteq Coll(\omega,\omega_1^V)$ is the $V$-generic filter associated with $g$, then 
\begin{enumerate}[(i)]
\item $j_g\rest \kappa$ is independent of $g$
\item $j_g(\kappa)$ is the largest Suslin cardinal of $(Lp(\mathbb{R}))^{V[G]}$ and hence is independent of $g$.
\item $j_g\rest \powerset^{\undertilde{\Gamma^*}}(\kappa)$ is independent of $g$.
\end{enumerate}
\end{claim}
\begin{proof}

To see (i), let $\gamma < \kappa$ be arbitrary and let $A\in Lp(\mathbb{R})$ be of Wadge rank $\gamma$. Note that since $\gamma < \kappa$, by our induction hypothesis, $A$ is $\omega_1$-UB as witnessed by $(S,W)$. Notice then that by Lemma \ref{lem:UBprojection}, 
\begin{center}
$j_g(A) = p[S]\cap V[g]$.
\end{center}
If (i) fails at $\gamma$, let $p \Vdash j_g(\gamma) = \gamma_0$ and $q \Vdash j_g(\gamma) = \gamma_1$ with $\gamma_0 \neq \gamma_1$. Let $g_0,g_1\subset \mathbb{P}_{\mathcal{I}}$ be $V$-generic such that $p\in g_0, q\in g_1$, and $V[g_1] = V[g_0]$; such $g_0, g_1$ can be easily obtained using the homogeneity of the forcing \footnote{Given $p\in g_0$, we can find an automorphism $\tau:\mathbb{P}_{\mathcal{I}}\rightarrow \mathbb{P}_{\mathcal{I}}$ such that $\tau(p_0)\leq q$. Then let $g_1 = \tau[g_0]$. $g_0, g_1$ are as desired.}. Let $M_0 = j_{g_0}(Lp(\mathbb{R}))$, $M_1 = j_{g_1}(Lp(\mathbb{R}))$. Note that $M_0\unlhd M_1$ or $M_1\unlhd M_0$. Write $j_i$ for $j_{g_i}$ and note that 
\begin{equation}\label{eqn:ord_diff}
j_0(\gamma) = \gamma_0 \neq \gamma_1 = j_1(\gamma).
\end{equation}
Note also by the fact that $V[g_0] = V[g_1]$, 
\begin{equation}\label{eqn:set_equal}
j_0(A) = j_1(A) = p[S]\cap V[g_1] = p[S]\cap V[g_0].
\end{equation}
The fact that the Wadge hierarchies of $M_0, M_1$ are compatible gives us
\begin{equation}\label{eqn:intersection}
j_0(A)= j_1(A) \in M_0\cap M_1.
\end{equation}
\ref{eqn:set_equal} and \ref{eqn:intersection} give us the Wadge rank of $j_0(A) = j_1(A)$ is $j_0(\gamma) = j_1(\gamma)$, which clearly contradicts \ref{eqn:ord_diff}. So (i) holds.

We now show (ii). We do not claim here that $(Lp(\mathbb{R}))^{V[G]}\models \sf{AD}^+$.  Suppose the statement of the claim is false. Fix  $g$ as above witnessing the failure of the claim. Then there is $\M\lhd (Lp(\mathbb{R}))^{V[G]}$ such that $j(\Gamma^*)$ is Suslin co-Suslin in $\M$ and $\M\vDash \sf{AD}^+$. 

By the scales analysis and $\sf{MC}$ in $\M$ ($\sf{MC}$ holds in $\M$ by our smallness assumption $(\ddag)$ and results in \cite{ATHM}), there is a sjs $\vec{A}$ sealing $\undertilde{Env(j(\Gamma^*))}$ in $V[G]$. Let $(\P,\Sigma)\in \M$ be guided by $\vec{A}$. By Boolean-valued comparisons (described in the previous section), there is an iterate $(\R,\Lambda)$ of $(\P,\Sigma)$ such that $\R\in V$ and $\Lambda\rest H^V_{\omega_2} \in V$. Now, $\Lambda$ has branch condensation and is $j(\Gamma^*)$-fullness preserving and hence by Lemma \ref{claim:bcshc} has strong hull condensation. By Lemma \ref{lem:pullback}, $\Lambda = j(\Lambda)^j$. Therefore, $\Lambda\in M$ and is $j(\Gamma^*)$-fullness preserving. 

Now note that $\Lambda$ is $\omega_1$-UB in $M$ by Lemma \ref{lem:UB}, and so by the core model induction similar to the above, $Lp^{\Lambda}(\mathbb{R}^M)\models \sf{AD}^+$ (here by density, $\mathbb{R}^M =\mathbb{R}^{V[G]}$). This implies that $L(\Lambda,\mathbb{R}^M)\models \Theta > \theta_0$ since $\Lambda\notin Lp(\mathbb{R})^M$. This in particular implies, via standard results (cf. \cite{wilson2012contributions}), that conclusion (a) and (b) holds for  $\undertilde{Env(j(\Gamma^*))}$ in $M$. By elementarity, (a) and (b) hold for  $\undertilde{Env(\Gamma^*)}$. This contradicts our assumption that (a), (b) fail.

To see that $j_g\rest \powerset^{\undertilde{\Gamma^*}}(\kappa)$ is independent of the choice of $g$ in (iii), fix a $\undertilde{\Gamma^*}$-prewellorder $\preceq$ of $\mathbb{R}$ of length $\kappa$; by choosing a minimal definition, we can assume $\preceq$ is definable from a real $y$ and $\kappa$ in $Lp(\mathbb{R})$. More precisely, we choose the least $\xi$ such that $Lp(\mathbb{R})|\xi$ ordinal defines such a $\preceq$ from a real $y$. By minimizing the ordinal parameters, we can then get that $\preceq$ is definable over $Lp(\mathbb{R})|\xi$ from $\{y,\kappa\}$, say  by formula $\varphi$. Note that any $X\in \powerset^{\undertilde{\Gamma^*}}(\kappa)$ is $\Sigma^1_1(\preceq,z)$ for some real $z$ by the Coding Lemma. Suppose $X$ witnesses the failure of (c) and $X$ is  $\Sigma^1_1(\preceq,z)$ for some real $z$. Let $g_0, g_1$ be such that $V[g_0] = V[g_1]$ and $j_i = j_{g_i}$ be the associated generic embeddings with the property that $j_1(X) \neq j_0(X)$. Let $\kappa^* = j_0(\kappa) = j_1(\kappa)$. By the choice of $\preceq$ and part (i), $j_0(\preceq) = j_1 (\preceq)$; this is because $j_0(\preceq), j_1(\preceq)$ are both definable from $\{\kappa^*, y\}$ via formula $\varphi$ over the least $\M\lhd (Lp(\mathbb{R}))^{V[g_0]}$ that ordinal defines a prewellorder of $\mathbb{R}^{V[g_0]}$ of length $\kappa^*$. Since $j_0(X), j_1(X)$ are $\Sigma_1^1$-definable from $j_0(\preceq)$ from $y$ via the same formula, $ j_0(X)= j_1(X)$. Contradiction.


\end{proof}

\begin{remark}\label{rem:clause_ii_revised}
In the proof of Claim \ref{claim:independent}(ii), it appears that we need to assume the failure of Theorem \ref{thm:sjs}(a). However, one can show

\medskip

\noindent (ii') $j_g(\kappa)$ is independent of $g$

\medskip

\noindent without assuming the failure of Theorem \ref{thm:sjs}(a). Suppose (ii') fails. We can then find $g_1, g_2$ such that $V[g_1] = V[g_2]$ and $j_{g_1}(\kappa)<j_{g_2}(\kappa)$. Let $j_i = j_{g_i}$ for $i\in 2$ and $j_i: V\rightarrow M_i$. We can run the argument in the proof of Claim \ref{claim:independent}(ii) to get $(\R,\Lambda)$ as there, where $\Lambda$ is $\omega_1$-UB in $M_1$. By elementarity, there is such a pair $(\R,\Lambda)\in V$ such that $\Lambda$ is $\omega_1$-UB, Lp$^\Lambda(\mathbb{R})\models \sf{AD}^+$, and $\Lambda\notin \textrm{Lp}(\mathbb{R})$. Since $j_0(\kappa)\neq j_1(\kappa)$, it is easy to see that $j_0(\Lambda)\neq j_1(\Lambda)$. But since $\Lambda$ is $\omega_1$-UB as witnessed by trees $(T,U)$ and $V[g_0] = V[g_1]$, $j_0(\Lambda) = p[T]\cap V[g_0] = p[j_1(T)]\cap M_1 = j_1(\Lambda)$. Contradiction.
\end{remark}

From the claim above and homogeneity, we easily see that the value of $j_g(\kappa), j_g(T)$ is independent of $g$; from now on, we will write $j(\kappa)$ for $j_g(\kappa)$ etc. Let 
\begin{center}
$\sigma = j'' \textrm{meas}^{\undertilde{\Gamma}}(\kappa^{<\omega})$. 
\end{center}
Note also that $\sigma$ is independent of $g$. Let $\lambda$ be the length of the well-ordering of Env$(\Gamma)$. We have $\lambda < j(\omega_1^V) = \omega_2^V$.  It follows that $j''\lambda$ (and hence also $\sigma$) is in $\textrm{Ult}(V,g)$ and is countable there. This then implies that $\sigma \in M$. 

Let $\mu\in\sigma$. Suppose $\mu$ concentrates on $j(\kappa)^n$ and let $\langle\mu_i \ | \ i\leq n \rangle$ be the projections of $\mu$, meaning $A\in \mu_i \iff \{ s\in j(\kappa)^n \mid s\rest i \in A \} \in \mu$. Note that $\mu_0$ is the trivial measure.

In Ult$(V,g)$, we define the following putative scale $\{\varphi_\mu : \mu\in\sigma\}$ on $\mathbb{R}\backslash p[j(T)]$ as follows. For each $\mu\in\sigma$, and for each $x\in \mathbb{R}\backslash p[j(T)]$ (so $j(T)_x$ is well-founded),
\begin{center}
$\varphi_\mu(x) = [\textrm{rank}_{j(T)_x}]_\mu$.\footnote{rank$_{j(T)_x}(t)$ denotes the rank of the node $t$ in the tree $j(T)_x$, and is considered to be zero if $t \notin j(T)_x$ and undefined if $j(T)_x$ is illfounded below $t$.}
\end{center}
We now define the following closed game $G^{\sigma,\mu}_{j(T)}$ in Ult$(V,g)$ (equivalently in $V[G]$, recalling that $\mathbb{R}^{V[G]}=\mathbb{R}^{\textrm{Ult}(V,g)}$ and the pointclass $j(\Gamma^*)$ is ordinal definable in $V[G]$): player I starts by playing $m_0, \ldots, m_n$ and $s_n$, $h_n$, and player II responds by playing a measure $\mu_{n+1}$. In each subsequent move (numbered $i > n$,) player I plays $m_{i}$, $s_i$, $h_i$, and player II plays a measure $\mu_{i+1}$.

\noindent \textbf{Rules for player I:}
\begin{itemize}
\item $m_k < \omega$ for all $k<\omega$
\item $j(T)_{(m_0,\ldots, m_{n-1})} \in \mu = \mu_n$
\item $s_i \in j_{\mu_i}(j(T)_{(m_0,\ldots, m_i)})$, and in particular $s_i \in j_{\mu_i}(j(\kappa))^{i+1}$ for all $i\geq n$
\item $s_n \supsetneq [\rm{id}]_{\mu_n}$
 \item $j_{\mu_i,\mu_{i+1}}(s_i) \subsetneq s_{i+1}$ for all $i\geq n$
\item $h_i \in \rm{OR}$ for all $i\geq n$
\item $j_{\mu_i,\mu_{i+1}}(h_i)>h_{i+1}$ for all $i\geq n$
\end{itemize}
\noindent \textbf{Rules for player II:}
\begin{itemize}
 \item $\mu_{i+1}\in \sigma$ is a measure on $j(\kappa)^{i+1}$ projecting to $\mu_i$
\item $\mu_{i+1}$ concentrates on the set $j(T)_{(m_0,\ldots,m_i)}\subset j(\kappa)^{i+1}$.
\end{itemize}
The first player that violates one of these rules loses, and if both players follow the rules for all $\omega$ moves, then player I wins.

The game is closed, hence determined by the Gale--Stewart theorem. Intuitively, player I is building a real $x = (m_0,m_1, \dots)$, player II is trying to build a tower $\vec{\mu}$ of measures in $\sigma$ concentrating on $j(T)_x$, and player I is trying to build a continuous witness
$\vec{h}$ to the illfoundedness of $\vec{\mu}$ as well as  a special kind of branch $(j_{i,\infty}(s_i) : i\geq n)$ through the
direct limit $j_{0,\infty}(j(T)_x)$ of $j(T)_x$ along $\vec{\mu}$. The following is the main lemma.
 
\begin{lemma}\label{lem:II_wins}
Player II has a winning strategy in the game $G^{\sigma,\mu}_{j(T)}$ for each $\mu\in \sigma$.
\end{lemma}
\begin{proof}
First note that $j(T)\in V$; this is because $T$ is ordinal definable in $V$.The parameter defining $j(T)$ in $V[G]$ has the form $j(s)$ for some finite sequence of ordinals $s\in V$. Therefore, $j(s)\in V$ and $j(T)\in V$ by homogeneity.\footnote{In the case $P = Lp^\Sigma(\mathbb{R})$, $T$ is ordinal definable from $\Sigma$ and there are trees $(W,S)$ witnessing $\Sigma$ is $\omega_1$-UB. Then $j(T)$ is ordinal definable in $V[G]$ from $(W,S)$ by the fact that $p[W] = p[j(W)]$ and $p[S] = p[j(S)]$ (see a similar calculation in the proof of Claim \ref{claim:j(K)inV}). Therefore, $j(T)\in V$ by homogeneity.}  In fact, by Claim \ref{claim:independent} and the remark after, $j(\kappa), j(T), j\rest \powerset^{\undertilde{\Gamma^*}}(\kappa^{<\omega})$ are independent of $g$.

Fix $\mu\in\sigma$. We define a winning strategy for player II in $G^{\sigma,\mu}_{j(T)}$ in Ult$(V,g)$. Let $\mu_0, \dots ,\mu_n$ be the projections of $\mu$ in order (here $\mu_n = \mu$). Let $j(\bar{\mu_i}) = \mu_i$ for $i = 0, \dots, n$. Note that for all $i$,
\begin{center}
$j_{\mu_i}\circ j = j\circ j_{\bar{\mu_i}}$.
\end{center}

Suppose player I starts the game by playing integers $m_0, \dots, m_n$, a finite sequence of ordinals $s_n\in j_{\mu_n}(j(T_{m_0,\dots,m_n}))\cap j_{\mu_n}(j(\kappa)^{n+1})$, and an ordinal $h_n$. Define the measure $\bar{\mu}_{n+1}\in \textrm{meas}^{\undertilde{\Gamma^*}}(\kappa^{<\omega})$ as follows.
\begin{center}
$X\in \bar{\mu}_{n+1} \iff s_n \in j_{\mu_n}(j(X)) $.
\end{center}
$\bar{\mu}_{n+1}$ is $OD^{V[g]}$ from a finite sequence of ordinals, some real $x\in\mathbb{R}^V$\footnote{The real $x$ can be taken to be the real that appears in the definition of $j^{-1}(\mu)$.} and $j_g\rest  \powerset^{\undertilde{\Gamma^*}}(\kappa)$. Since $j_g\rest  \powerset^{\undertilde{\Gamma^*}}(\kappa)$ is independent of $g$,  $\bar{\mu}_{n+1}\in V$.\footnote{In the general case $P=Lp^\Sigma(\mathbb{R})$, we reach the same conclusion because $\bar{\mu}_{n+1}$ is OD$^{V[g]}$ from a real, a finite sequence of ordinals, $j_g\rest \powerset^{\undertilde{\Gamma^*}}(\kappa^{<\omega})$, and $(W,S)$, where $(W,S)$ witnesses $\Sigma$ is $\omega_1$-UB.}


For $i>n$, suppose player I has played an integer $m_i$, a finite sequence of ordinals $s_i\in j_{\mu_i}(j(T_{m_0,\dots,m_i}))\cap j_{\mu_i}(j(\kappa)^{i+1})$, and an ordinal $h_i$. Define the measure $\bar{\mu}_{i+1}\in \textrm{meas}^{\undertilde{\Gamma^*}}(\kappa^{<\omega})$ as follows.
\begin{center}
$X\in \bar{\mu}_{i+1} \iff s_i \in j_{\mu_i}(j(X))$.
\end{center}
As before, the measure $\bar{\mu}_{i+1}$ is in $V$, concentrates on $T_{m_0,\dots,m_i}$, and projects to $\bar{\mu}_i$. Let player II play the measure $\mu_{i+1} = j(\bar{\mu}_{i+1})$.

Assume for contradiction that player I is able to play $\omega$ many moves, following all the rules of the game. We get a real $x = (m_0, m_1, \dots)$, a tower of measures $(\mu_i : i<\omega)$ in $\sigma$, and a countable sequence of ordinals $(h_i : i < \omega)$ witnessing the illfoundedness of this tower. By elementarity, the tower $(\bar{\mu}_i : i<\omega)$ is also illfounded. 

Take a wellfounded tree $W \in \bigcup_{x\in \mathbb{R}} L[T,x]$ on $\kappa$ on which each measure $\bar{\mu}_i$ in this tower concentrates,
and such that the function $\bar{h}:\omega\rightarrow \textrm{Ord}$ defined by $\bar{h}(i) = [\textrm{rank}_W]_{\mu_i}$ is a pointwise
minimal witness to the illfoundedness of the tower $(\bar{\mu}_i : i<\omega)$ (see \cite[Lemma 3.5.9]{wilson2012contributions}). Then by
the elementarity of $j$, the function $h = j(\bar{h})$ is a pointwise minimal witness to the illfoundedness
of the tower $(\mu_i : i<\omega)$.\footnote{Actually we only need the minimality of $h(n)$.}  Because $\bar{\mu}_i$ concentrates on $W$ we have $s_i \in j_{\mu_i}(j(W))$ for all $i <\omega$. Define a function $h': \omega\rightarrow \textrm{Ord}$ by $h'(i) = \textrm{rank}_{j_{\mu_i}(j(W))}(s_i)$. Then from the rules for player I concerning the finite sequences $s_i$ we have $j_{\mu_i,\mu_{i+1}}(h'(i)) > h'(i+1)$ and also $h'(n) < \textrm{rank}_{j_{\mu_n}(j(W))}([\textrm{id}]_{\mu_n}) = h(n)$, contradicting the minimality of $h(n)$.
\end{proof}

\begin{remark}
In the above proof, we use $\sf{CH}$ in a crucial way. $\sf{CH}$ implies that $\lambda < \omega_2^V$ and we in turns get that $\sigma\in M$ and is countable there. These two facts are key for the proof. As mentioned in the introduction, without $\sf{CH}$ the existence of an $\omega_1$-dense ideal on $\omega_1$ is equiconsistent with $\sf{AD}$. 
\end{remark}

The proof of Claim \ref{claim:independent} and the argument in the following remark give us the following useful corollary.
\begin{corollary}\label{cor:independent}
Suppose $A\in \Gamma^*$ is $\omega_1$-UB and let $\gamma = w(A)$ in $\Gamma^*$. Then $j_g(\gamma)$ is independent of $g$.
\end{corollary}


\begin{lemma}\label{lem:scale}
In Ult$(V,g)$, the set of norms $\{\varphi_\mu : \mu\in \sigma\}$ defined by $\varphi_\mu(x) = [\textrm{rank}_{j(T)_x}]_\mu$ (or more precisely, any enumeration of this countable set of norms in order type $\omega$) is a scale on the complement of $p[j(T)]$.
\end{lemma}
\begin{proof}
Work in Ult$(V,g)$. Let $\mu\in \sigma$. We say that $\sigma$ \textit{stabilizes}\footnote{The idea of this definition comes from a similar notion of stability used in unpublished work of S.~Jackson.} $\mu$ if, whenever $(x_k : k  < \omega)$ is a sequence of
reals in $\mathbb{R}\backslash p[j(T)]$ converging to a limit $x$ and such that for each $\mu' \in \sigma$, the ordinals $\varphi_{\mu'}(x_k)$
are eventually constant, we have $\varphi_\mu(x)\leq \textrm{lim}_{k\rightarrow \omega}\varphi_\mu(x_k)$. (In particular, $\varphi_\mu(x) < \infty$.)  

It is clear from the definition that if $\sigma$ stabilizes every $\mu\in \sigma$, then $\{\varphi_\mu : \mu\in \sigma\}$ is a scale. So fix a measure $\mu\in\sigma$. We want to show $\sigma$ stabilizes $\mu$. Suppose not. We describe a winning strategy for player I in $G^{\mu,\sigma}_{j(T)}$. Let  $(x_k : k < \omega)$ witness that $\sigma$ does not stabilize $\mu$. That is, $x_k\in \mathbb{R}\backslash p[j(T)]$
for each $k <\omega$, and the sequence of ordinals $(\varphi_\nu(x_k): k<\omega)$ has an eventually constant
value $h(\nu)$ for each measure $\nu\in\sigma$ but the limit $x$ of the sequence $(x_k : k < \omega)$ satisfies $\varphi_\mu(x) > \textrm{lim}_{k\rightarrow \omega}\varphi_\mu(x_k)$. (This includes the possibility that $\varphi_\mu(x) =\infty$.)

Define $m_i = x(i)$ and $h(\nu) = \textrm{lim}_{k\rightarrow \omega} \varphi_\nu(x_k)$. Let $n$ be the unique integer such that $\mu$ concentrates on $j(\kappa)^n$ and let $\mu_i$ be the projection of $\mu$ onto $j(\kappa)^i$ for all $i\leq n$. In particular, $\mu_n = \mu$. By definition,
\begin{center}
$\varphi_{\mu_n}(x) = [s\mapsto \textrm{rank}_{j(T)_x}(s)]_{\mu_n} = \textrm{rank}_{j_{\mu_n}(j(T)_x)}([\textrm{id}]_{\mu_n}) > h(\mu_n)$.
\end{center}
So there is a finite sequence $s_n \supsetneq [id]_{\mu_n}$ with rank $\geq h(\mu_n)$ in the tree $j_{\mu_n}(j(T)_x)$. Let player I play as his first move the integers $m_0, \dots, m_n$, the ordinal $h_n = h(\mu_n)$, and $s_n$, where $s_n$ is the least such sequence. For $i\geq n$, we will show inductively that player I can maintain the inequality
\begin{equation}\label{eqn:rank-inequality}
\operatorname{rank}_{j_{\mu_i}(j(T)_x)}(s_i) \geq h(\mu_i).
\end{equation}
Whenever player II plays a measure $\mu_{i+1}$ according to the rules of the game, we have
\begin{center}
rank$_{j_{\mu_{i+1}}(j(T)_x)}(j_{\mu_i,\mu_{i+1}}(s_i)) = j_{\mu_i,\mu_{i+1}}(\textrm{rank}_{j_{\mu_i}(j(T)_x)}(s_i)) \geq j_{\mu_i,\mu_{i+1}}(h_i)> h_{i+1}$.
\end{center}
To show the last step $j_{\mu_i,\mu_{i+1}}(h_i)> h_{i+1}$, we argue as follows. Recall that for each $l$ we have $h_l = h(\mu_l) = \textrm{lim}_{k\rightarrow \omega}\varphi_{\mu_l}(x_k)$. Since the measure $\mu_{i+1}$ concentrates on $j(T)_{x\rest (i+1)}$ and projects to $\mu_i$, for each $k$ we have
\begin{center}
$j_{\mu_i,\mu_{i+1}}(\varphi_{\mu_{i}}(x_k)) = j_{\mu_i,\mu_{i+1}}([\textrm{rank}_{j(T)_{x_k}}]_{\mu_{i}}) = [\textrm{ext}_{i,i+1}\textrm{rank}_{j(T)_{x_k}}]_{\mu_{i+1}}$,
\end{center}
where the ``extension'' of a function $F: j(\kappa)^i \rightarrow \textrm{Ord}$ to $j(\kappa)^{i+1}$ is defined by ext$_{i,i+1} F(s) = F(s\rest i)$ for all $s\in j(\kappa)^{i+1}$. Note that
\begin{center}
$ [\textrm{ext}_{i,i+1}\textrm{rank}_{j(T)_{x_k}}]_{\mu_{i+1}} > [\textrm{rank}_{j(T)_{x_k}}]_{\mu_{i+1}} = \varphi_{\mu_{i+1}}(x_k)$.
\end{center}
Finally, since for each $l$ the ordinal $h_l$ is the eventual value of $\varphi_{\mu_l}(x_k)$ as $k \to \omega$, consideration of sufficiently large $k$ gives $j_{\mu_i,\mu_{i+1}}(h_i)> h_{i+1}$.

This shows that player I can choose a successor $s_{i+1} \supsetneq j_{\mu_i,\mu_{i+1}}(s_i)$ of rank at least $h(\mu_{i+1})$ in the tree
$j_{\mu_{i+1}}(j(T)_x)$, thereby maintaining the desired inequality \eqref{eqn:rank-inequality} for one more step. Then player I can play the integer $m_{i+1} = x(i+1)$, the
least such finite sequence $s_{i+1}$, and the ordinal $h_{i+1} = h(\mu_{i+1})$.
By playing in this way, player I can follow the rules forever. This contradicts the previous lemma, which showed that player II has a winning strategy.
\end{proof}

The previous claims and elementarity establishes (a) for $A$ being the universal $\check{\Gamma^*}$-set. By standard arguments, see \cite[Section 4.3]{wilson2012contributions}, the rest of (a) and (b) follow. This contradicts our assumption. Therefore, (a) and (b) hold after all.

\end{proof}

\begin{theorem}\label{thm:omega1UB}
There is a hod pair $(\P',\Sigma')$ in $V$ such that 
\begin{enumerate}
\item $\Sigma'$ is Lp$^\Sigma(\mathbb{R})$-fullness preserving and $\Sigma'\notin $ Lp$^\Sigma(\mathbb{R})$.
\item $\Sigma'$ has branch condensation.
\item $\Sigma'$ is $\omega_1$-UB.
\end{enumerate}
\end{theorem}
\begin{proof}

Let $\Gamma^*$ be the largest Suslin pointclass of $Lp^\Sigma(\mathbb{R})$. Let $\vec{A} = (A_i : i<\omega)$ be the sjs sealing $\undertilde{Env(\Gamma^*)}$ as in the previous theorem. Let $(\P'',\Sigma'')$ be a pair such that $\Sigma''$ is guided by $\vec{A}$. $\Sigma''$ has properties (1) and (2), but (3) may fail for $\Sigma''$. Here one can regard $\P''$ as a $\Sigma$-suitable mouse with one Woodin cardinal or a hod mouse. We take the first viewpoint and hence we regard $\Sigma''$ as an iteration strategy for $\P''$ as a $\Sigma$-mouse (so all $\vec{\T}$ according to $\Sigma''$ are above $\P$ and iterates of $\P''$ according to $\Sigma''$ are $\Sigma$-premice).

For each $p\in Coll(\omega,\omega_1)$, let $G_p$ be the ``finite variation" of $G$ induced by $p$ and let $g_p$ be the corresponding $\mathbb{P}_\mathcal{I}$-generic induced by $\pi$ and $G_p$. We let $\vec{A}^p = (A^p_i:i<\omega)$ be $j_{g_p}(\vec{A})$. Let $(\Q_p,\Sigma_p)$ be a hod pair in $V[g] = V[g_p]$ guided by $\vec{A}^p$\footnote{We can take $(\Q_p,\Sigma_p)$ to be $(\P'', j_{g_p}(\Sigma''))$.} and $(N,\Lambda)$ be obtained by Boolean comparing all $(\Q_p,\Sigma_p)$. So $N\in V$ and $\Lambda\rest V\in V$ is a strategy acting on stacks in $H^V_{\omega_2}$ such that $\Lambda$ is has branch condensation (and is guided by $\mathcal{B} = \bigcup_{p} rng(\vec{A}^p)$), strong hull condensation, and is $j(\Gamma^*)$-fullness preserving. Note that $\R$ is countable in $M$ and $\Lambda \notin j(Lp^\Sigma(\mathbb{R}))$.

Applying Lemma \ref{lem:pullback}, we get that $\Lambda = j(\Lambda)^j$. By elementarity, in $V$, there is a pair $(\P',\Sigma')$ and an elementary embedding $\pi: \P'\rightarrow \R$ such that
\begin{enumerate}[(a)]
\item $\Sigma' = \Lambda^\pi$.
\item $\Sigma'\rest HC$ is $\Gamma^*$-fullness preserving and has branch condensation.
\item $\P'$ is a countable $\Sigma$-mouse (i.e. $\pi\rest \P = id$) that is $\Gamma^*$-suitable.
\end{enumerate}



$(\P',\Sigma')$ satisfies (1) and (2). We note that property (a) above gives that $\Sigma'$ is an $(\omega_2,\omega_2)$-strategy.   Now Lemma \ref{lem:UB} implies that $\Sigma'$ is $\omega_1$-UB. This completes the proof of the theorem. 
\end{proof}

\section{The LIMIT CASE}\label{sec:lim}

Recall we let $g\subseteq \mathbb{P}_\mathcal{I}$ be $V$-generic and $j = j_g: V\rightarrow M = \textrm{Ult}(V,g)$ be the corresponding ultrapower map; by our hypothesis, $g$ corresponds to a $V$-generic $G\subset \textrm{Coll}(\omega,\omega_1)$. We also let $k:M \rightarrow N$ be the generic ultrapower map induced by a generic $h\subset j(\mathbb{P}_\mathcal{I})$. We remind the reader that $\sf{CH}$ holds, so the continuum $\mathfrak{c}$ is $\omega_1$.

Let $\langle\theta_\alpha : \gamma < \gamma \rangle$ be the Solovay sequence computed in $\Gamma$ (our maximal model) and $\Theta = \textrm{sup}_{\gamma} \theta_\gamma$. By the previous section, $\gamma$ is a limit ordinal and $\Theta$ is the Wadge ordinal of $\Gamma$. For $\alpha \leq \Theta$, by $\Gamma\rest \alpha$, we mean the set of $B\in\Gamma$ such that the Wadge rank of $B$ is less than $\alpha$. We also remind the reader that our inductive hypothesis implies that every $B\in \Gamma$ is $\omega_1$-UB; in particular, because $\Theta$ is a limit of Suslin cardinals in $\Gamma$, by Corollary \ref{cor:independent}, $j\rest \Theta$ is independent of $G$. First we claim 
\begin{center}
$|\Gamma|\leq \mathfrak{c}$. 
\end{center}
\begin{lemma}\label{lem:theta_small}
Suppose $|\Gamma| = \mathfrak{c}^+$. Then $\Gamma = \powerset(\mathbb{R})\cap L(\Gamma,\mathbb{R})$.
\end{lemma}
\begin{proof}
\indent Suppose not. Let $\alpha$ be the least such that $\rho_\omega(J_\alpha(\Gamma,\mathbb{R})) = \mathbb{R}$, i.e. $J_\alpha(\Gamma,\mathbb{R})$ defines a set of reals $A$ such that $A\notin \Gamma$. Hence $\alpha \geq \mathfrak{c}^+$ by our assumption. Let $f: \alpha \times \Gamma \twoheadrightarrow J_\alpha(\Gamma,\mathbb{R})$ be a surjection that is definable over $J_\alpha(\Gamma,\mathbb{R})$ (from parameters). 

We first define a sequence $\langle H_i \ | \ i<\omega\rangle$ as follows. Let $H_0 = \mathbb{R}$. By induction, suppose $H_n$ is defined and there is a surjection from $\mathbb{R}$ onto $H_n$. Suppose $(\psi, a)$ is such that $a\in H_n$ and $J_\alpha(\Gamma,\mathbb{R}) \vDash \exists x\psi[x,a]$. Let $(\gamma_{a,\psi},\beta_{a,\psi})$ be the $<_{lex}$-least pair such that there is a $B\in \Gamma$ with Wadge rank $\beta_{a,\psi}$ such that
\begin{center}
$J_\alpha(\Gamma,\mathbb{R}) \vDash \psi[f(\gamma_{a,\psi},B),a]$.
\end{center}
Let then $H_{n+1} = H_n\cup \{f(\gamma_{a,\psi},B) \ | \ J_\alpha(\Gamma,\mathbb{R}) \vDash \exists x\psi[x,a]\wedge w(B) = \beta_{a,\psi}\wedge a\in H_n\}$. It's easy to see that there is a surjection from $\mathbb{R}$ onto $H_{n+1}$. This uses the fact that $\Theta^\Gamma = \mathfrak{c}^+$ is regular, which implies sup$\{\beta_{a,\psi} \ | \ a\in H_n \wedge L_\alpha(\Gamma,\mathbb{R})\vDash \exists x\psi[x,a]\} < \Theta = \mathfrak{c}^+$. Let $H = \bigcup_n H_n$. By construction, $H\prec J_\alpha(\Gamma,\mathbb{R})$. Finally, let $M$ be the transitive collapse of $H$. 
\\
\indent Say $M = J_\beta(\Gamma^*,\mathbb{R})$. By construction, it is easy to see that $\Gamma^* = \Gamma\rest \theta_\gamma$ for some $\gamma$ such that $\theta_\gamma < \Theta$. But then $\rho_\omega(J_\beta(\Gamma^*,\mathbb{R})) = \mathbb{R}$.\footnote{For instance, to see that $\Gamma\rest \theta_0 \subset \Gamma^*$, let $A\in\Gamma$ be $OD$ in $J_\alpha(\Gamma,\mathbb{R})$ from a real $x$. Suppose $A\notin M$. By minimizing the Wadge rank of $A$ and minimizing the ordinal parameters defining $A$, we may assume $A$ is definable in $J_\alpha(\Gamma,\mathbb{R})$ from $x$. By elementarity, $A$ is definable in $M$ from $x$, so $A\in \Gamma^*$. Contradiction.} This contradicts that $\Gamma^*$ is constructibly closed. 
\end{proof}

The lemma gives $\Gamma = \powerset(\mathbb{R})\cap L(\Gamma,\mathbb{R})$ and in fact, $L(\Gamma,\mathbb{R})\vDash ``\textsf{AD}_\mathbb{R} + \Theta$ is regular". This is because $\Theta = \mathfrak{c}^+$ in this case. This contradicts $(\ddag)$. Therefore, $|\Gamma|\leq \mathfrak{c}$ as desired.

Let $\mathcal{H}$ be the direct limit of hod pairs $(\P,\Sigma)\in \Gamma$ such that $\Sigma$ has branch condensation and is fullness preserving under iteration embeddings by $\Sigma$. So $\lambda^\mathcal{H}$ is a limit ordinal. For each $\alpha < \lambda^{\mathcal{H}}$, let $\Sigma_\alpha$ be the strategy of $\mathcal{H}(\alpha)$ in $j(\Gamma)$ obtained as a tail of some (any) $j(\Sigma)$, where $(\P,\Sigma)$ is a hod pair in $\Gamma$ with branch condensation and is fullness preserving such that $\M(\P,\Sigma) = \mathcal{H}(\alpha)$.  Let 
\begin{center}
$\Sigma = \oplus_{\alpha<\lambda^\mathcal{H}} \Sigma_\alpha$. 
\end{center}

Now note that
\begin{center}
$j$ is continuous at $\lambda^\mathcal{H}$ if and only if cof$^V(\lambda^\mathcal{H}) = \omega$.
\end{center}


First note that $j\rest \omega_1^V \in M$. If $j$ is continuous at $\lambda^{\mathcal{H}}$ and cof$^V(\lambda^\mathcal{H})=\omega_1$, then $j(\omega_1^V)$ is singular in $M$. This contradicts the fact that $j(\omega_1^V)$ is a successor cardinal, hence regular, in $M$. This implies cof$^V(\lambda^{\mathcal{H}})\neq \omega_1^V$ and hence   cof$^V(\lambda^{\mathcal{H}}) = \omega$. 



\begin{lemma}\label{lem:in_V}
\begin{itemize}
\item $\Sigma\rest V \in V$ and $\Sigma$ does not depend on $G$.
\item $j\rest \Theta^\Gamma$ is independent of $G$.
\end{itemize}
\end{lemma}

\begin{proof}
This follows from our induction hypothesis, i.e. for each $\alpha$, $\Sigma_\alpha$ is $\omega_1$-UB by the inductive hypothesis, and hence $\Sigma_\alpha\rest V\in V$ and does not depend on $G$. This gives the first item. The argument for the second item is given at the beginning of the section.

\end{proof}

Let
\begin{equation}\label{eqn:cases}
\mathcal{H}^+ = \begin{cases*} 
Lp^{\Sigma,j(\Gamma)}(\mathcal{H}) & if $\forall \M\lhd Lp^{\Sigma,j(\Gamma)}(\mathcal{H}) \ \rho_\omega(\M)\geq \Theta$ \\
\P                                                       & where $\P\lhd Lp^{\Sigma,j(\Gamma)}(\mathcal{H})$ is least $\N$ such that $\rho_\omega(\N)<\Theta$.
\end{cases*} 
\end{equation}

To be technically correct, by $Lp^{\Sigma,j(\Gamma)}(\mathcal{H})$ we mean $Lp^{\Sigma}(\mathcal{H})$ defined inside $L(j(\mathbb{R}), C)$ for some $C\in j(\Gamma)$. This makes sense as $\Sigma\in j(\Gamma)$ and the Solovay sequence of $j(\Gamma)$ has limit length. By Lemma \ref{lem:in_V}, we get that $$\mathcal{H}^+\in V.$$ This is because $\mathcal{H}^+$ is definable in $V[G]$ from $\mathcal{H}, \Sigma\rest V$ and by Lemma \ref{lem:in_V}, $\Sigma\rest V \in V$ and does not depend on $G$.


\begin{proposition}\label{cor:H_small}
$|\mathcal{H}^+| \leq \mathfrak{c}$. Therefore, $j\rest \mathcal{H}^+\in M$
\end{proposition}
\begin{proof}
Suppose we have $\mathcal{H}^+ = \textrm{Lp}^{\Sigma,j(\Gamma)}(\mathcal{H})$. If $|\mathcal{H}^+|=\mathfrak{c}^+$, we would get an $\omega_1$-sequence of distinct reals in $j(\Gamma)$, noting that $(\mathfrak{c}^+)^V=\omega_2^V$ is $\omega_1$ in $M$ by the density of $\mathcal{I}$. Contradiction. Therefore, $|\mathcal{H}^+| = \mathfrak{c}$, and hence $\mathcal{H}^+$ is countable in $M$. Again, by density of $\mathcal{I}$, $j\rest \mathcal{H}^+\in M$. A similar argument also works for the second case of (\ref{eqn:cases}).
\end{proof}

Using the embedding $j$, the fact that $j\rest \mathcal{H}^+\in M$, and the construction in \cite[Section 11]{sargsyanCovering2013}, we obtain a strategy $\Lambda$ for $\mathcal{H}^+$ such that
\begin{enumerate}
\item $\Lambda$ extends $\Sigma$;
\item for any $\Lambda$-iterate $\P$ of $\mathcal{H}^+$ via a stack $\vec{\mathcal{T}}$ such that $i^{\vec{\mathcal{T}}}$ exists, there is an embedding $\sigma:\P \rightarrow j(\mathcal{H}^+)$ such that $j\rest\mathcal{H}^+ = \sigma \circ i^{\vec{\mathcal{T}}}$. Furthermore, letting $\Lambda_\P$ be the $\vec{\mathcal{T}}$-tail of $\Lambda$, for all $\alpha< \lambda^\P$, $\Lambda_{\P(\alpha)}\in j(\Gamma)$ has branch condensation.
\item $\Lambda$ is $\Gamma(\mathcal{H}^+,\Lambda)$-fullness preserving.
\end{enumerate}

We outline the construction here. We first briefly review definitions and notations related to the analysis of stacks in \cite[Section 6.2]{ATHM} summarized in Section \ref{sec:hodmice_prelim}; see \cite[Section 6.2]{ATHM} for a more detailed discussion.

\begin{definition}[$j$-realizable iterations]\label{dfn:realizable_iterations}
Let $\vec{\T}\in HC^M$ be a stack on $\mathcal{H}^+$. We say $\vec{\T}$ is \textbf{$j$-realizable} if there is a sequence $\langle \sigma_\R : \R\in tn(\vec{\T}) \rangle$ such that
\begin{enumerate}
\item $\sigma_{\mathcal{H}^+} = j\rest\mathcal{H}^+$; for all $\R\in tn(\vec{\T})$, $\sigma_\R: \R\rightarrow j(\mathcal{H}^+)$.
\item For $\R,\Q\in tn(\vec{\T})$ such that $\R\prec^{\vec{\T},s} \Q$, $\sigma_\R = \sigma_\Q \circ \pi^{\vec{\T}}_{\R,\Q}$.
\item For every $\R\in ntn(\vec{\T})$, there is a reasonable hod pair $(\S_\R,\Lambda_\R)\in j(\Gamma)$ that is $j(\Gamma)$-fullness preserving and has branch condensation such that $\sigma_\R[\R(\xi^{\vec{\T},\R}+1)]\subset \textrm{rng}(\pi^{\Lambda_\R}_{\S_\R,\infty})$.

\item For every $\R\in ntn(\vec{\T})$, letting $(\S_\R,\Lambda_\R)$ be as above, and letting $k_\R: \R(\xi^{\vec{\T},\R}+1)\rightarrow \S_\R$ be given by: $k_\R(x) = y$ if and only if $\sigma_\R(x) = \pi^{\Lambda_\R}_{\S_\R,\infty}(y)$ and $k_\R \vec{\T}_\R$ is according to $\Lambda_\R$.

\item For every $\R\in ntn(\vec{\T})$, let $\S^*_\R$ be the last model of $k_\R \vec{\T}_\R$ and let $\Q_\R$ be the last model of $\vec{\T}_\R$ (considered as a stack on all of $\R$. Suppose $\pi^{\vec{\T}_\R}$ is defined (hence, $\Q_\R\in tn(\vec{\T})$ and $\R\prec^{\vec{\T},s} \Q_\R$). Let $k^*_\R: \Q_\R(\zeta)\rightarrow \S^*_\R$ be the natural map that comes from the copying construction, where $\Q_\R(\zeta)$ is the image of $\R(\xi^{\vec{\T},\R}+1)$ under the iteration embedding of $\vec{\T}_\R$. Then we define $\sigma_{\Q_\R}:\Q_\R\rightarrow j(\mathcal{H}^+)$ as follows: for all $x\in \Q_\R$, 
\begin{center}
$\sigma_{\Q_\R}(x) = \sigma_\R(f)(\pi^{\Lambda}_{\S^*_\R,\infty}(k^*_\R(a))$,
\end{center}
where $f\in \R$, and $a\in [\Q(\pi^{\vec{\T}}_{\R,\Q_\R}(\xi^{\vec{\T},\R}+1))]^{<\omega}$ are such that $x = \pi^{\vec{\T}}_{\R,\Q_\R}(f)(a)$; here $\Lambda = (\Lambda_\R)_{k_\R\vec{\T}_\R,\S^*_\R}$.

\item For every trivial terminal node $\R$, for every $\xi < \lambda^\R$, there is a reasonable hod pair $(\S_\R,\Lambda_\R)\in j(\Gamma)$ where $\Lambda$ is $j(\Gamma)$-fullness preserving, and has branch condensation and $\sigma_\R(\xi+1)\subset \textrm{rng}(\pi^{\Lambda_\R}_{\S_\R,\infty})$.
\end{enumerate}

The maps $(\sigma_\R: \R\in tn(\vec{\T}))$ are the $j$-realizable embeddings of $\vec{\T}$. In the above, we may also choose $(\S_\R,\Lambda_\R)$ such that letting $j(\mathcal{H})(\alpha) = \M_\infty(\S_\R,\Lambda_\R)$, then $\alpha$ is minimal.
\end{definition}

Now we define the domain of the strategy $\Lambda$. Basically, it consists of $j$-realizable stacks. See \cite[Definition 11.5]{sargsyanCovering2013}.

\begin{definition}\label{dfn:dom}
Let $\vec{\T}\in HC^M$ be a stack of on $\mathcal{H}^+$.\footnote{$\vec{\T}$ either has a strongly linear, closed and cofinal set $C\subseteq tn(\vec{\T})$ or $\vec{\T}_{\S_{\vec{\T}}}$ is of limit length.} We let $\vec{\T}\in dom(\Lambda)$ iff $\vec{\T}$ is $j$-realizable. Define $\Lambda(\vec{\T}) = b$ iff $\vec{\T}^\smallfrown b$ is $j$-realizable.
\end{definition}

\begin{lemma} \label{lem:total}
Whenever $\vec{\T}\in dom(\vec{\T})$, then $\Lambda(\vec{\T})$ is defined.
\end{lemma} 

See \cite[Lemma 11.6]{sargsyanCovering2013} for a similar argument. In other words, the lemma states that if $\vec{\T}$ is $j$-realizable and has no last model, then we can find a cofinal branch $b$ of $\vec{\T}$ so that $\vec{\T}^\smallfrown b$ is $j$-realizable. We sketch the argument here. 
\begin{proof}
Suppose there is a strongly closed, cofinal $C\subset tn(\vec{\T})$. In this case $\vec{\T}$ has a unique, cofinal, non-dropping branch $b$ determined by $C$. Let $\Q = \M^{\vec{\T}}_b$ and $\sigma_\Q: \Q\rightarrow j(\mathcal{H}^+)$ be the direct limit of the maps $\{\sigma_\R: \R\in C\}$; more precisely, let $\sigma_\Q(x) = y$ if and only if there is some $x^*\in \R$ for some $\R\in C$ such that $\pi^{\vec{\T}}_{\R,\Q}(x^*) = x$ and $\sigma_\R(x^*) = y$. It is easy to see that $\sigma_\Q$ is well-defined and satisfies the clauses of Definition \ref{dfn:realizable_iterations} (note that in this case, $\Q$ is a trivial terminal node). 

Otherwise, we are looking for a branch of $\vec{\T}_{\S_{\vec{\T}}}$. Let $\R = \S_{\vec{\T}}$ and $\U =  \vec{\T}_{\S_{\vec{\T}}}$. By our hypothesis, objects like $\sigma_\R, k_\R, (\S_\R,\Lambda_\R)$ as in (3) and (4) can be defined. Let then $b = \Lambda_\R(k_\R\U)$,  $\Q = \M^{\U}_b$, $\S^* = \M^{k_\R\U}_b$,  $k:\Q(\zeta) \rightarrow \S^*$, $\sigma_\Q: \Q\rightarrow j(\mathcal{H}^+)$ be the objects as described in (5) above. So $b$ is the branch of $\vec{\T}_{\S_{\vec{\T}}}$ we are looking for.

In the following, we assume $\Q$ is a terminal node; otherwise, we're done. We need to verify clause (6) in the case $\Q$ is a trivial terminal node. The case for non-trivial terminal nodes has been dealt with as above. Without loss of generality, we assume $\Q\neq \R$ and there is a $\U$ on $\R$ with last model $\Q$ such that $\pi^\U_{\R,\Q}$ exists. We let $\sigma_\R,k_\R, (\S_\R,\Lambda_\R), \S^*, k$ be the objects associated with $\R, \U, \Q$ as before. We let $\Lambda = (\Lambda_\R)_{k_\R\U,\S^*}$ and $\sigma_\Q = \pi^{\Lambda}_{\S^*,\infty}\circ k$.  Fix $\xi < \lambda^\Q$. Let $(\W,\Psi)\in j(\Gamma)$ be a reasonable hod pair such that $\Psi$ is $j(\Gamma)$-fullness preserving, and such that $\M_\infty(\W,\Psi) = j(\mathcal{H}^+)(\sigma_\Q(\xi+1))$. We can then find $(\S,\Psi_\S)\in I(\W,\Psi)$ such that $\sigma_\Q[\Q(\xi+1)]\subset \textrm{rng}(\pi^{\Psi_\S}_{\S,\infty})$. We are done.
\end{proof}

\begin{remark}\label{rmk:independent}
Suppose $\vec{\T}\in dom(\Lambda)$, then there is at most one $b$ such that $\vec{\T}^\smallfrown b$ is $j$-realizable. In the proof of Lemma \ref{lem:total}, the only case to verify is when $\S_{\vec{\T}}$ exists. Let $\R, \U, \sigma_\R,k_\R, (\S_\R,\Lambda_\R)$ be as there. Suppose $(\S^*_\R,\Lambda^*_\R)$ and $l_\R:\R(\xi^{\vec{\T},\R}+1)\rightarrow \S^*_\R$ are such that $l_\R(x) = y$ if and only if $\sigma_\R(x) = \pi^{\Lambda^*_\R}_{\S^*_\R,\infty}(y)$, $l_\R \vec{\T}_\R$ is according to $\Lambda^*_\R$ and $c= \Lambda^*_\R(l_\R\U)$. To see $b =c$, we let $(\S,\Psi)$ be the common iterate of $(\S_\R,\Lambda_\R)$ and $(\S^*_\R,\Lambda^*_\R)$. Let $\sigma_0: \S_\R\rightarrow \S$ and $\sigma_1: \S^*_\R\rightarrow \S$ be the iteration maps. So $\Lambda_\R = (\Psi)^{\sigma_0}$ and $\Lambda^*_\R=(\Psi)^{\sigma_1}$ because these strategies are pullback consistent. It is also easy to verify that $$\sigma_0\circ k_\R = \sigma_1\circ l_\R;$$ this is because letting $\tau: \S\rightarrow \sigma_\R(\R(\xi^{\vec{\T},R}+1))$ be the direct limit embedding according to $\Psi$, then $$\sigma_\R = \tau\circ \sigma_0\circ k_\R = \tau\circ \sigma_1\circ l_\R.$$ So $\sigma_0\circ k_\R = \sigma_1\circ l_\R$ as desired. Therefore, 
\begin{center}
$b = \Psi^{\sigma_0\circ k_\R}(\U) = \Psi^{\sigma_1\circ l_\R}(\U) = c$.
\end{center} 
\end{remark}

Clearly, if $\Lambda$ is a $j$-realizable strategy, then $\Lambda$ satisfies (1) and the first clause of (2); by basic hod mice theory (cf. \cite{ATHM}), $\Lambda$ also satisfies the ``Furthermore" clause. By the proof of \cite[Lemma 11.8]{sargsyanCovering2013}, we can choose $\Lambda$ so that $\Gamma(\mathcal{H}^+,\Lambda)$ is Wadge minimal (amongst all strategies $\Lambda$ constructed this way) and this particular choice of $\Lambda$ satisfies (3) as well.

\begin{lemma}\label{lem:no_proj_across}
$\mathcal{H}^+ = Lp^{\Sigma,j(\Gamma)}(\mathcal{H})$ and if $j$ is discontinuous at $\lambda^\mathcal{H}$, then $\mathcal{H}^+ \models \textrm{cof}(\lambda^{\mathcal{H}})$ is measurable.
\end{lemma}
\begin{proof}
The second clause follows from the first clause and the case assumption that $j$ is discontinuous at $\lambda^\mathcal{H}$. To see this, assume the first clause. If $\mathcal{H}^+\models ``\lambda^{\mathcal{H}}$ is regular", then by standard results on Vopenka forcing (cf. \cite{trangThesis2013}) $L[\mathcal{H}^+](\Gamma)\cap \powerset(\mathbb{R}) = \Gamma$ and therefore, $L(\Gamma,\mathbb{R})\models ``\sf{AD}_{\mathbb{R}} + $$\Theta$ is regular", contradicting our smallness assumption $(\ddag)$. If $\mathcal{H}^+\models ``\lambda^{\mathcal{H}}$ is singular", then letting $\kappa = \textrm{cof}^{\mathcal{H}^+}(\lambda^\mathcal{H})$, then $\kappa$ must be measurable in $\mathcal{H}^+$. This is because $j\rest (\kappa+1)$ is the iteration embedding of $\mathcal{H}(\alpha)$ according to $\Psi =_{def} \Sigma_{\mathcal{H}(\alpha)}$ in $M$ for some (equivalently any) $\alpha$ such that $\kappa\in \mathcal{H}(\alpha)$; therefore, $i^\Psi_{\mathcal{H}(\alpha),\infty}$ is discontinuous at $\kappa$,\footnote{If $j$ is continuous at $\kappa$, we show that $j$ is continuous at $\lambda^\mathcal{H}$. Suppose $f:\kappa\rightarrow \lambda^{\mathcal{H}}$ is cofinal and increasing and $f\in \mathcal{H}^+$. Then $j(f)\in j(\mathcal{H}^+)$, and $j(f): j(\kappa)\rightarrow j(\lambda^\mathcal{H})$ is cofinal and increasing. But $j(\kappa) = \textrm{sup} \ j''\kappa$, therefore, $j(\lambda^\mathcal{H}) = \textrm{sup} \ j''\lambda^\mathcal{H}$.} implying $\kappa$ is measurable in $\mathcal{H}(\alpha)$, hence in $\mathcal{H}^+$. 

Now, suppose for contradiction that there is a $\P\lhd \mathcal{H}^+$ such that $\rho_\omega(\P)<\Theta$. Let $\P$ be the least such. Let $\beta<\lambda^{\mathcal{H}}$ be least such that $\rho_\omega(\P)\leq \delta_\beta^\P$ and $\delta^\P_\beta > \textrm{cof}^\P(\lambda^\P)$, here $\lambda^\P = \lambda^{\mathcal{H}}$ and $\delta^\P_\alpha = \delta^\mathcal{H}_\alpha$ for all $\alpha <\lambda^\P$. $\P$ can be considered a hod premouse over $(\mathcal{H}(\beta),\Sigma_{\beta})$. Using $j$ and the construction in \cite[Section 11]{sargsyanCovering2013} discussed above, we can define a strategy $\Lambda$ for $\P$ such that $\Lambda$ acts on stacks above $\delta_\beta^\P$ and extends $\oplus_{\alpha<\lambda^\P}\Sigma_{\alpha}$ (the strategy is simply $\oplus_{\alpha<\lambda^\P}\Sigma_{\alpha}$ for stacks based on $\mathcal{H}$ (above $\delta^\P_\beta$), but the point is that it also acts on all of $\P$ because of $j$). This is because given a stack $\vec{\T}$ according to $\Lambda$, there is a map $\sigma: \M^{\vec{\T}} \rightarrow j(\P)$ such that $\sigma \circ i^{\vec{\T}} = j\rest \P$, where for any $f\in \P$, any generator $a$ used along the main branch of $\vec{\T}$, say $a\in \M^{\vec{\T}}(\gamma)$ and $\M^{\vec{\T}}(\gamma)$ is the image of $\P(\gamma^*)$, then letting $\Psi = \Sigma_{\gamma^*}$,
\begin{center}
$\sigma(i^{\vec{\T}}(f)(a)) = j(f)(i^{\Psi_{\vec{\T},\M^{\vec{\T}}(\gamma)}}(a))$.
\end{center}
In the above, we note that $i^{\vec{\T}}$ is continuous at $\lambda^\P$, so we can find $\gamma, \gamma^*$.

Note that $\Lambda$ has branch condensation. By a core model induction as in the successor case, we get that $\Lambda \in j(\Gamma)$.\footnote{$\Lambda$ is essentially $\Sigma$, acting on stacks above $\delta^\P_\beta$, so it has branch condensation. The core model induction (in $M$) as done so far works for $\Lambda$, showing that in $M$ we can uniquely extend $\Lambda$ to an $(\omega_2,\omega_2)$-strategy and  and $\Lambda$ is $\omega_1$-UB. We then proceed to show $\M_1^{\Lambda,\sharp}$ exists, and $Lp^\Lambda(\mathbb{R})\models \sf{AD}^+$ just like before.} In $j(\Gamma)$, let $\mathcal{F}$ be the direct limit system of $\Sigma_{\beta}$-hod pairs $(\Q,\Psi)$ Dodd-Jensen equivalent to $(\P,\Lambda)$.\footnote{$(\P,\Lambda)$ is an anomalous hod pair in the terminology of \cite{ATHM}. $(\Q,\Psi)$ is Dodd-Jensen equivalent to $(\P,\Lambda)$ means that there are non-dropping iterates $(\Q^*,\Psi^*)$ of $(\Q,\Psi)$ and $(\P^*,\Lambda^*)$ of $(\P,\Lambda)$ such that $(\Q^*,\Psi^*) = (\P^*,\Lambda^*)$.} $\mathcal{F}$ can be characterized as the direct limit system of $\Sigma_{\beta}$-hod pairs $(\Q,\Psi)$ in $j(\Gamma)$ such that $\Psi$ is $\Gamma(\P,\Lambda)$-fullness preserving and has branch condensation and $\Gamma(\Q,\Psi)=\Gamma(\P,\Lambda)$. $\mathcal{F}$ only depends on $\Sigma_{\beta}$ and the Wadge rank of $\Gamma(\P,\Lambda)$ and hence is $OD^{L(j(\mathbb{R}),C)}_{\Sigma_{\beta}}$ for some $C\in j(\Gamma)$. 

Fix such a $C$ and note that $L(j(\mathbb{R}),C)\vDash \sf{AD}^+ + \sf{SMC}$. Let $A\subseteq \delta^\P_\beta$ witness $\rho_\omega(\P)\leq \delta_\beta^\P$, that is, $A\notin\P$ and there is a formula $\phi$ such that for all $\alpha\in \delta^\P_\beta$,
\begin{center}
$\alpha\in A \Leftrightarrow \P \vDash \phi[\alpha,p]$,
\end{center}
where $p$ is the standard parameter of $\P$. Now $A$ is $OD_{\Sigma_{\beta}}$ in $L(j(\mathbb{R}),C)$; this is because letting $\M_\infty$ be the direct limit of $\mathcal{F}$ under iteration maps, then in $L(j(\mathbb{R}),C)$, $\M_\infty \in \textrm{HOD}_{\Sigma_{\beta}}$ and $A$ witnesses that $\rho_\omega(\M_\infty)\leq \delta_\beta^\P$. By $\sf{SMC}$ in $L(j(\mathbb{R}),C)$ and the fact that $\mathcal{H}(\beta+1)$ is $j(\Gamma)$-full, we get that $A\in \P$. This is a contradiction.

\end{proof}

\begin{remark} \label{rmk:non_trivial}
The construction of $\Lambda$ is nontrivial in the case that $\mathcal{H}^+ \vDash \textrm{cof}(\Theta)$ is measurable; otherwise, $\Lambda$ is simply $\Sigma$ but because of $j$, it acts on all of $\mathcal{H}^+$ by an argument as in the proof of Lemma \ref{lem:no_proj_across}. So from this point on, we assume \textit{$j$ is not continuous at $\lambda^{\mathcal{H}^+}$}.
\end{remark}

\begin{definition}[Nice strategies]\label{NiceStrategy}
Suppose $\pi_{\mathcal{H}^+,\R}:\mathcal{H}^+\rightarrow \R$, $\sigma: \R\rightarrow j(\mathcal{H}^+)$ are elementary and $\R$ is countable in $M$. Suppose $j\rest \mathcal{H}^+ = \sigma\circ \pi_{\mathcal{H}^+,\R}$. Let $\alpha<\lambda^\R$. We say that an iteration strategy $\Lambda_{\R(\alpha)}$ for $\R(\alpha)$ is \textbf{nice} if and only if
\begin{enumerate}[(i)]
\item $\Lambda_{\R(\alpha)}$ is a $j(\Gamma)$-fullness preserving strategy for $\R(\alpha)$ with branch condensation. $\Lambda_{\R(\alpha)}$ is also positional and commuting.
\item $\pi^{\Lambda_{\R(\alpha)}}_{\R(\alpha),\infty} = \sigma'\rest \R(\alpha)$ for some elementary map $\sigma': \R\rightarrow j(\mathcal{H}^+)$ such that $j\rest\mathcal{H}^+ = \sigma' \circ \pi_{\mathcal{H}^+,\R}$ (so $\Lambda_{\R(\alpha)}$ acts on all of $\R$). 
\item If $\pi_{\mathcal{H}^+,\R}\in M$, then $\Sigma_\alpha\rest M \in M$.
\end{enumerate}
\end{definition}

We want to show some $j$-realizable strategies are nice. This will be accomplished through the next several lemmas.

\begin{lemma}\label{fullness}
Let $\vec{\T},\R,\sigma_\R$ be as above. Then $\R$ is full in $j(\Gamma)$. In fact, letting $\pi: \mathcal{H}^+\rightarrow \R$ and $\sigma: \R \rightarrow j(\mathcal{H}^+)$ be arbitrary elementary embeddings such that $j\rest \mathcal{H}^+ = \sigma\circ \pi$, then $\R$ is full in $j(\Gamma)$. Furthermore, $j\rest \mathcal{H}^+$, and hence $\pi$, must be continuous at $o(\mathcal{H}^+)$.

\end{lemma}
\begin{proof}

We show the last statement of the lemma. The argument is very similar for all the other statements; we briefly indicate the changes at the end of the proof. Suppose $j$ is not continuous at $o(\mathcal{H}^+)$. Suppose without loss of generality that $\pi$ is not continuous at $o(\mathcal{H}^+)$.  Indeed the general case can be reduced to this case. Suppose $j\rest \mathcal{H}^+ = \sigma' \circ \pi'$, where $\sigma': \R' \rightarrow j(\mathcal{H^+})$ is discontinuous at $o(\R')$ and $\pi': \mathcal{H}^+\rightarrow \R'$ is continuous at $o(\mathcal{H}^+)$. In $M$, let $\sigma: \R\rightarrow j(\mathcal{H}^+)$ be elementary such that $\R$ is countable, transitive and rng$(\sigma')\subseteq \textrm{rng}(\sigma)$; such a $\sigma$ can easily be found in $M$ by considering a countable hull $X\prec H^{M}_{\omega_2}$ that contains all relevant objects, then $\sigma$ can be taken to be the restriction of the uncollapse map associated with $X$.  Let $\pi = \sigma^{-1}\circ \sigma'$. It is easy to see then that $\pi$ is not continuous at $o(\mathcal{H}^+)$.

This means there is a mouse $\M\lhd \textrm{Lp}^{\oplus_{\beta<\lambda^\R} \Lambda_{\R(\beta)}}(\R|\delta^\R) = \R$ such that $\M\notin \R|\gamma$ where $\gamma = sup \ \pi[o(\mathcal{H}^+)]$. We take $\M$ to be the least such and let $\Sigma_\M$ be the unique strategy for $\M$ in $j(\Gamma)$ (acting on trees on $\M$ above $\delta^\R$).



\begin{claim}\label{claim:witness}
There is a $\Sigma$-hod pair $(\P,\Phi)$ such that 
\begin{enumerate}[(a)]
\item $\P\in V$, $\Phi\rest V\in V$,\footnote{By $\Phi\rest V$, we mean $\Phi\rest H^V_{\mathfrak{c}^+}$.} and $\Phi\in j(\Gamma)$ is fullness preserving and has branch condensation.
\item $\P$ is countable in $M$, $\lambda^{\P}$ is limit and cof$^{\P}(\lambda^{\P})$ is not measurable in $\P$.
\item $\Phi = j(\Phi)^j$.
\item in $j(\Gamma)$, $\exists \beta$ such that $\Gamma(\P,\Phi)=j(\Gamma)|\theta_{\beta+\omega}$ and $\Gamma(\P,\Phi)|\theta_\beta \models ``(\M,\Sigma_\M)$ witnesses $\pi$ is not continuous at $o(\mathcal{H}^+)$" . 
\item $o(\mathcal{H}^+)$ is a cardinal of $\P$, i.e. $\P \models ``\mathcal{H}^+$ is full."
\end{enumerate}
\end{claim}
\begin{proof}

First note that in $M$, there is some $\alpha$ such that $\Sigma_\M$, the canonical strategy of $\M$, is in $j(\Gamma)|\delta^{\P^*}_\beta$, where $\P^* = \textrm{HOD}^{j(\Gamma)}_{\Sigma}(\alpha)$ \footnote{We identify $\textrm{HOD}^{j(\Gamma)}_{\Sigma}$ with the direct limit of $\Sigma$-hod pairs $(\R,\Psi)$ and $\Psi$ is fullness preserving and has branch condensation in $j(\Gamma)$.}  and $\P^*\models \alpha = \beta+\omega$.  Such $\P^*$ and $\alpha$ exists by our assumptions on $\Gamma$. $\P^*\in V$ follows from homogeneity. Let $\Psi^*$ be the strategy of $\P^*$ which is the tail of some (equivalently, all) $\Sigma$-hod pair $(\R^*,\Psi)\in j(\Gamma)$ where $\Psi$ is fullness preserving and has branch condensation in $j(\Gamma)$ and $\M_\infty(\R,\Psi) = \P^*$. $\Psi^*$ is fullness preserving and has branch condensation in $k(j(\Gamma))$. It follows that $\Psi^*\rest V \in V$: we can ordinal define $\Psi^*\rest V$ in $V[G]$ from $\Sigma$ and $\P$ with the prescription above, using the fact that $j(\Gamma)$ is $OD$ in $V[G]$; so by homogeneity, $\Psi^*\rest V\in V$. 

We want to find a countable-in-$M$ version of $\P^*$ in $V$. Let $(\R,\Psi)$ be a $\Sigma$-hod pair in $j(\Gamma)$ such that $\M_\infty(\R,\Psi)=\P^*$ and $\Psi$ has strong hull condensation, branch condensation, and is $j(\Gamma)$-fullness preserving. By boolean comparisons, there is a $\Psi$-iterate $(\P,\Phi)$ such that $(\P,\Phi)$ satisfies (a). (b) is clear from the choice of $\P^*$. (c) follows from Lemma \ref{lem:pullback}. (d) follows from the choice of $\P^*$ and the fact that $\Gamma(\P,\Phi) = j(\Gamma)|\delta^{\P^*}_\alpha$.

To see (e), suppose not and for simplicity, let $\mathcal{H}^+ \unlhd \N\lhd \P$ be least such that $\rho_1(\N) = \Theta$. Let $f: \kappa^* \rightarrow \Theta$ be an increasing and cofinal map in $\mathcal{H}^+$, where $\kappa^* = \textrm{cof}^{\mathcal{H}^+}(\Theta)$. $\N$ is intercomputable with the sequence $g = \langle \N_\alpha \ | \ \alpha < \kappa^*\rangle$, where $\N_\alpha = Th_{\Sigma_1}^\N( \delta^{\mathcal{H}^+}_{f(\alpha)}\cup \{p_\N\})$. Note that $\N_\alpha\in \mathcal{H}^+$ for each $\alpha<\kappa^*$. Now let $\R_0 = \textrm{Ult}_{0}(\mathcal{H}^+,\mu)$, $\R_1= \textrm{Ult}_{1}(\N,\mu)$, where $\mu\in \mathcal{H}^+$ is the (extender on the sequence of $\mathcal{H}^+$ coding a) measure on $\kappa^*$ with Mitchell order $0$. Let $i_0:\mathcal{H}^+ \rightarrow \R_0$, $i_1:\N\rightarrow \R_1$ be the ultrapower maps. Letting $\delta = \delta_{\lambda^{\mathcal{H}^+}} = \Theta$, it's easy to see that $i_0\rest (\kappa^*+1)=i_1\rest(\kappa^*+1)$ and $\powerset(\delta)^{\R_0} = \powerset(\delta)^{\R_1}$. The second equality follows from the fact that $\R_0$ is full in $j(\Gamma)$ (and hence in $k(j(\Gamma))$).

\end{proof}

Let $(\P,\Phi)$ be as in the claim. Let $\pi^+: \P\rightarrow \S$ be the ultrapower map derived from the $\pi$-extender of length $\delta^\R$. We note that $\pi^+$ is continuous at $o(\mathcal{H}^+)$ and by elementarity, $\S \models ``\R|\gamma$ is full". Therefore, $\M\notin \S$. 

Let $\tau^+:\S\rightarrow j(\P)$ be the factor map, so $j\rest \P = \tau^+\circ \sigma^+$. Let $\Psi = j(\Phi)^{\tau^+}$. By (c) of the claim, $\Phi = \Psi^{\sigma^+}$. Therefore, 
\begin{center}
$\Gamma(\P,\Phi) \subset \Gamma(\R,\Psi)$.
\end{center}
Hence, $ \Gamma(\R,\Psi)|\theta_\beta \models ``(\M,\Sigma_\M)$ witnesses $\pi$ is not continuous at $o(\mathcal{H}^+)$". Now, we iterate $\S$ using $\Psi$ at the top $\omega$ Woodin cardinals of $\S$ to make $\mathbb{R}^M$ generic.\footnote{More precisely, we write $(\delta_i^{\mathcal{S}}: i<\omega)$ for the top $\omega$ Woodin cardinals of $\mathcal{S}$ and a similar notation applies to iterates of $\mathcal{S}$. We work in $M[L]$ where $L\subseteq Coll(\omega,\mathbb{R}^M)$. We have a generic enumeration $(x_n : n<\omega)$ of $\mathbb{R}^M$ and we have a sequence of normal trees and models $(\T_n , \mathcal{S}_n: n<\omega)$ according to $\Psi$, where $\T_0$ is on $\mathcal{S}=\mathcal{S}_0$, $\T_n$ is a $x_n$-genericity iteration tree on $\mathcal{S}_n$ on the window $(\delta^{\mathcal{S}_n}_{n-1},\delta^{\mathcal{S}_n}_n)$ according to the $\T_{n-1}$-tail of $\Psi$, here $\delta_{-1}^{\mathcal{S}}=0$. Letting $\mathcal{S}_\infty$ be the direct limit, then $\mathbb{R}^M$ is the symmetric reals of $\mathcal{S}_\infty$ for some $g\subseteq Coll(\omega, <\lambda)$, where $\lambda$ is the supremum of the Woodin cardinals of $\mathcal{Y}_\infty$.\label{ftn:gen_iter}} Let $\S^*$ be the resulting model. The derived model of $\S^*$ at $\delta^{\S^*}$ satisfies:
\begin{center}
$L(\Gamma(\S^*,\Psi_{\S^*})|\theta_\beta) \models ``\M$ is a sound $\oplus_{\beta<\lambda^\R} \Lambda_{\R(\beta)}$-mouse such that $\rho_\omega(\M) = \delta^\R$ but $\M$ is not in $\R|\gamma$".
\end{center}
On the other hand, $\S^* \models ``\R|\gamma$ is full with respect to sound $\oplus_{\beta<\lambda^\R} \Lambda_{\R(\beta)}$-mice projecting to $\delta^\R$". This contradicts the displayed line above.

We have shown that $\pi$ is continuous and that $\R$ is full ``at the top", i.e. for every $\M\lhd \textrm{Lp}^{\oplus_{\beta<\lambda^\R} \Lambda_{\R(\beta)}}(\R|\delta^\R)$, there is $\alpha < \gamma$ such that $\M \lhd \R|\alpha$. The remaining clause of fullness is proved in an almost identical manner. Suppose there is a strong cut point $\xi$ such that letting $\alpha < \lambda^{\R}$ be the largest such that $\delta^{\R}_\alpha \leq \gamma$, then in $j(\Gamma)$, there is a mouse $\M\lhd \textrm{Lp}^{\Sigma_{\Q(\alpha)}}(\Q|\gamma)$ such that $\M\notin \R$. The argument given above can be carried out verbatim to obtain a contradiction.

\end{proof}

\begin{definition}\label{dfn:good}
In $M$, suppose $X\prec (H_{\mathfrak{c}^+}, \in)$ is countable.\footnote{Sometimes, we just write $H_{\mathfrak{c}^+}$ for $(H_{\mathfrak{c}^+}, \in)$ for brevity. Also, note that $\mathfrak{c}^+=\omega_2$ in $M$ by elementarity.} $X$ is \textit{good} if letting $\pi_X: M_X\rightarrow X$ be the uncollapse map, 
\begin{enumerate}[(a)]
\item $j[\mathcal{H}^+]\cup \{j(\mathcal{H}^+)\}\subset \textrm{rng}(\pi_X)$;
\item $\mathcal{H}^+\cup \{\mathcal{H}^+\}\subset M_X$;
\item letting $\P_X = \pi_X^{-1}(j(\mathcal{H}^+))$, then $\P_X$ is $j(\Gamma)$-full and for any $\alpha < \lambda^{\P_X}$, $\pi_X\rest \P_X(\alpha) = i^{\Lambda^X_\alpha}_{\P_X(\alpha),\infty}$, where $\Lambda^X_\alpha$ is a tail of $\Lambda$ for some (equivalently any) hod pair $(\Q,\Lambda)\in j(\mathcal{F})\cap X$ such that $\Lambda$ is $j(\Gamma)$-fullness preserving and has branch condensation and $(\M_\infty(\Q,\Lambda))^{M_X} = \P_X(\alpha)$.
\end{enumerate}
\end{definition}
\begin{remark}\label{rmk:good}
\begin{enumerate}[(a)]
\item Note that if $X$ is good, then $\P_X$ is the transitive collapse of $Hull^{j(\mathcal{H}^+)}(j[\mathcal{H}^+] \cup \oplus_{\alpha < \lambda^{\P_X}}i^{\Lambda^X_\alpha}_{\P_X(\alpha),\infty})$. 
\item Letting $X^* = Hull^{H^V_{\mathfrak{c}^+}}(\mathcal{H}^+)$ and $X = j[X^*]$, then $X$ is good.
\item Any good $X$ is cofinal in $o(j(\mathcal{H}^+))$ by Lemma \ref{fullness}.
\end{enumerate}
\end{remark}

\begin{lemma}\label{lem:good_hulls}
In $M$, the set $\{X\cap \mathbb{R} : X \textrm{ is good}\}$ is in $j(\mathcal{F}_\mathcal{I})$ and the set of good $X$ is closed and unbounded.
\end{lemma}
\begin{proof}
Let $X$ be as in Remark \ref{rmk:good}(b) and let $Y\prec (H_{\mathfrak{c}^+},\in)$ be countable in $M$, $X\prec Y$, and $\mathcal{H}^+\cup\{\mathcal{H}^+\}\subset Y$. Since $\mathcal{H}^+$ is countable in $M$, there is a club of such $Y$. Clearly, (a) and (b) in Definition \ref{dfn:good} hold for $Y$. For (c), using the notation above and Lemma \ref{fullness}, we have that $\P_Y$ is $j(\Gamma)$-full. Furthermore, for all $\alpha<\lambda^{\P_Y}, \pi_Y\rest \P_Y(\alpha) = i^{\Lambda^Y_\alpha}_{\P_Y(\alpha),\infty}$ by elementarity of $\pi_Y$.
\end{proof}

Suppose $X$ is a good hull, we let $j_X: \mathcal{H}^+\rightarrow \P_X$ be $j_X = \pi_X^{-1}\circ j$. We let $\Lambda_X$ be the strategy for $\P_X$ defined from $\pi_X$ the same way $\Lambda$ is defined from $j$ for $\mathcal{H}^+$ (again, we take $\Lambda_X$ with $\Gamma(\P_X,\Lambda_X)$ minimal). By Lemma \ref{fullness} and the fact that $X$ is good, $\Lambda_X$ is $j(\Gamma)$-fullness preserving. By \cite{ATHM}, there is an iterate $(\T_X, \Q_X)$ of $(\P_X,\Lambda_X)$ such that letting $\Psi_X = (\Lambda_X)_{\T_X,\Q_X}$, $\Psi_X$ has branch condensation, and is commuting (see \cite{ATHM}).  Let now $\M^X_\infty = \M_\infty(\Q_X,\Psi_X)$. Note that $\M_\infty^X = j(\mathcal{H}^+)(\gamma)$ for some $\gamma < j(\lambda^\mathcal{H})$ and $\M^X_\infty$ does not depend on the choice of $(\Q_X,\Psi_X)$.

By construction of $\Lambda_X$, there is a map $m_X: \M^X_\infty \rightarrow j(\mathcal{H}^+)$ such that
\begin{center}
$\pi_X\rest \P_X = m_X \circ i^{\Psi_X}_{\Q_X,\infty}\circ i^{\T_X}$.\footnote{Recall we assume $j$ is discontinuous at $\lambda^{\mathcal{H}^+}$. Othewise, $\M^X_\infty = j(\mathcal{H}^+)$ and $m_X$ is the identity.}
\end{center}

We need a strong form of condensation to show $\mathcal{H}^+\models ``\Theta$ is regular"; basically, this form of condensation will imply that if $m_X$ is nontrivial, then
\begin{center}
crt$(m_X) = \delta^{\M^X_\infty}$.\footnote{It could be that $M^X_\infty = j(\mathcal{H}^+)$ and $m_X$ is the identity map. In which case, we cannot conclude $\Theta$ is regular in $\mathcal{H}^+$. In this case, $\Gamma(\mathcal{H}^+,\Lambda) = j(\Gamma)$. We then simply continue the core model induction. See Section \ref{sec:outline}.} 
\end{center}
Therefore, $\M^X_\infty \vDash ``\delta^{\M^X_\infty}$ is a regular cardinal which is a limit of Woodin cardinal." This easily implies $\Theta$ is regular in $\mathcal{H}^+$.

The following definition originates from \cite[Definition 11.14]{sargsyanCovering2013}. Let $\mathfrak{S}$ be the set of good hulls. For each $X\in \mathfrak{S}$, let $\Theta_X = j_X(\Theta)$. 

\begin{definition}
\label{ACondensation}
Suppose $X\in \mathfrak{S}$ and $A\in \P_X \cap \powerset(\Theta_X)$. We say that $\pi_X$ has \textit{$A$-condensation} if whenever there are elementary embeddings $\upsilon: \P_X \rightarrow \Q$, $\tau:\Q\rightarrow j(\mathcal{H}^+)$ such that $\Q$ is countable in $M$ and $\pi_X = \tau\circ \upsilon$, then 
\begin{center}
$\upsilon(T_{\P_X,A}) = T_{\Q,\tau,A}$, 
\end{center}
where  
\begin{center}
$T_{\P_X,A} = \{(\phi,s) \ | \ s\in [\Theta_X]^{<\omega}\wedge \P_X \vDash \phi[s,A]\}$,
\end{center}
and
\begin{center}
$T_{\Q,\tau,A} = \{(\phi,s) \ | \ s\in [\delta_\alpha^\Q]^{<\omega} \textrm{ for some } \alpha<\lambda_\Q \wedge j(\mathcal{H}^+) \vDash \phi[i^{\Sigma^{\tau,-}_{\Q}}_{\Q(\alpha),\infty}(s),\pi_X(A)]\}$,
\end{center}
where $\Sigma^{\tau,-}_\Q$ is the $\tau$-pullback strategy of $j(\Sigma)$.\footnote{$\Sigma^{\tau,-}_\Q = \oplus_{\alpha<\lambda^\Q} j(\Sigma)^{\tau}_{\Q(\alpha)}$.} 

We say $\pi_X$ has \textit{condensation} if it has $A$-condensation for every $A\in \P_X \cap \powerset(\Theta_X)$. 
\end{definition}


\begin{theorem}[$j$-condensation lemma]\label{thm:cond} 
Let $X^* = Hull^{H^V_{\mathfrak{c}^+}}(\mathcal{H}^+)$ and $X = j[X^*]$; so $\P_X = \mathcal{H}^+$, $\Theta_X=\Theta$, and $\pi_X\rest \P_X = j\rest \P_X$. Then $\pi_X$ has condensation.
\end{theorem}
\begin{proof}
Fix $A\in \P_X\cap\powerset(\Theta_X)$. We show that $\pi_X$ has $A$-condensation. Suppose not. 

We first claim that if $Y\in \mathfrak{S}$ is such that $X\prec Y$ and $\pi_Y$ has $\pi_{X,Y}(A)$-condensation, then $\pi_X$ has $A$-condensation. Fix such a $Y$. Note that $k(\pi_X) = k(\pi_Y)\circ \pi_{X,Y}$ and $k(\pi_Y) = k\rest j(\P_X)\circ \pi_Y$. By elementarity, $k(\pi_Y)$ has $\pi_{X,Y}(A)$-condensation in $N$ and hence $k\rest j(\P_X)$ has $j(A)$-condensation in $N$, by the following calculations: for any countable $\R$ in $N$, suppose there are embeddings $i:j(\P_X)\rightarrow \R$ and $\tau: \R\rightarrow k(j(\P_X))$ such that $k\rest j(\P_X) = \tau \circ i$, then
\begin{align*}
i(T_{j(\P_X),j(A)}) & = i(\pi_Y(T_{\P_Y,\pi_{X,Y}(A)})) \\
& = T_{\R,\tau, \pi_{X,Y}(A)} \\
&= T_{\R,\tau, j(A)};
\end{align*}
the second equality uses the fact that $k(\pi_Y)$ has $\pi_{X,Y}(A)$-condensation in $N$ and $k(\pi_Y) = \tau\circ i\circ \pi_Y$. Therfore, $\pi_X$ has $A$-condensation (in $M$) by the elementarity of $j$.

Suppose now for every $Y\in \mathfrak{S}$ such that $X\prec Y$, $\pi_Y$ does not have $\pi_{X,Y}(A)$-condensation. Recall that if $(\P,\Sigma)$ is a hod pair such that $\delta^\P$ has measurable cofinality then we let $\Sigma^- = \oplus_{\alpha < \lambda^{\P}} \Sigma_{\P(\alpha)}$. We say that a tuple $\{\langle \P_i,\Q_i,\tau_i,\xi_i,\pi_i,\sigma_i \ | \ i<\omega \rangle, \M^Y_{\infty}\}$ is a \textbf{bad tuple} (see Figure \ref{fig:diamond}) if
\begin{enumerate}
\item $Y\in \mathfrak{S}$;
\item $\P_i = \P_{X_i}$ for all $i$, where $X_i\in \mathfrak{S}$;
\item $X_0 = X$ and for all $i < j$, $X_i \prec X_j \prec Y$;
\item for all $i$, $\xi_i:\P_i\rightarrow \Q_i$, $\sigma_i:\Q_i \rightarrow \M^Y_{\infty}$,  $\tau_i: \P_{i+1}\rightarrow \M^Y_{\infty}$, and $\pi_i: \Q_i \rightarrow \P_{i+1}$;
\item for all $i$, $\tau_i = \sigma_i\circ \xi_i$, $\sigma_i = \tau_{i+1} \circ \pi_i$, and $\pi_{X_i,X_{i+1}}\rest \P_i=_{\textrm{def}} \phi_{i,i+1} = \pi_i \circ \xi_i$;
\item $\phi_{i,i+1}(A_i) = A_{i+1}$, where $A_i = \pi_{X,X_i}(A)$;
\item for all $i$, $\xi_i(T_{\P_i,A_i}) \neq T_{\Q_i,\sigma_i,A_i}$. 
\end{enumerate}
In (7), $T_{\Q_i,\sigma_i,A_{X_i}}$ is computed relative to $\M^Y_{\infty}$, that is
\begin{center}
$T_{\Q_i,\sigma_i,A_i} = \{(\phi,s) \ | \ s\in [\delta_\alpha^{\Q_i}]^{<\omega} \textrm{ for some } \alpha<\lambda^{\Q_i} \wedge \M^Y_{\infty} \vDash \phi[i^{\Sigma^{\sigma_i,-}_{\Q_i}}_{\Q_i(\alpha),\infty}(s),\tau_i(A_i)]\}$
\end{center}

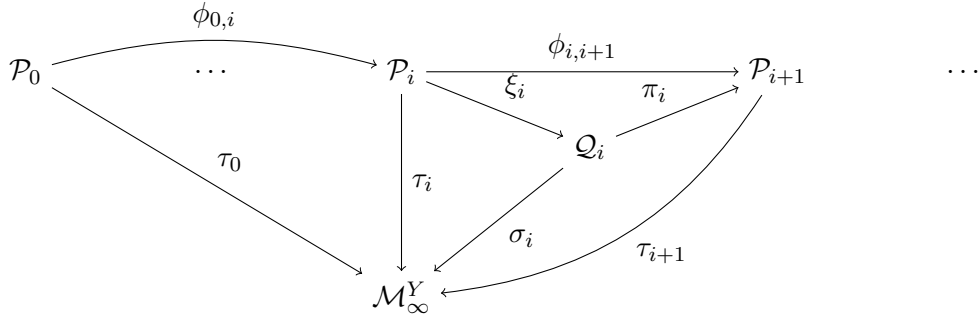
\begin{figure}
\centering
\begin{tikzpicture}[node distance=2.5cm, auto]
  \node (A) {$\P_0$};
  \node (B) [right of=A] {$\dots$};
  \node (C) [right of=B] {$\P_i$};
  \node (D) [right of=C] {};
  \node (E) [node distance=1cm, below of=D] {$\Q_i$};
  \node (G) [right of=D] {$\P_{i+1}$};
  \node (H) [right of=G] {$\dots$};
  \node (I) [node distance=3cm, below of=C]{$\M^Y_\infty$};
  \draw[->,bend left=15] (A) to node {$\phi_{0,i}$}(C);
  \draw[->] (A) to node {$\tau_0$}(I); 
  \draw[->] (C) to node {$\tau_i$}(I);
  \draw[->, bend left=25] (G) to node {$\tau_{i+1}$}(I);
  \draw[->] (C) to node  {$\xi_i$} (E);
  \draw[->] (E) to node {$\pi_i$} (G);
  \draw[->] (C) to node {$\phi_{i,i+1}$} (G);
  \draw[->] (E) to node {$\sigma_i$} (I);
  \end{tikzpicture}
\caption{A bad tuple}
\label{fig:diamond}
\end{figure}

\noindent \textbf{Claim: } There is a bad tuple.
\begin{proof}
For brevity, we first construct a bad tuple $\{\langle \P_i,\Q_i,\tau_i,\xi_i,\pi_i,\sigma_i \ | \ i<\omega \rangle, j(\mathcal{H}^+)\}$ with $j(\mathcal{H}^+)$ playing the role of $\M^Y_{\infty}$. We then simply choose a sufficiently large $Y\in \mathfrak{S}$ and let $i_Y:\P_Y\rightarrow \M^Y_{\infty}$ be the direct limit map, $m_Y: \M^Y_{\infty}\rightarrow \mathcal{H}^+$ be the natural factor map, i.e. $m_Y\circ i_Y = \pi_Y$. It's easy to see that for all sufficiently large $Y$, the tuple $\{\langle \P_i,\Q_i,m_Y^{-1}\circ \tau_i,m_Y^{-1}\circ \xi_i,m_Y^{-1}\circ \pi_i,m_Y^{-1}\circ \sigma_i \ | \ i<\omega \rangle, \M^Y_{\infty}\}$ is a bad tuple. But the existence of such a tuple $\{\langle \P_i,\Q_i,\tau_i,\xi_i,\pi_i,\sigma_i \ | \ i<\omega \rangle, j(\mathcal{H}^+)\}$ follows from our assumption.
\end{proof}

By essentially the same proof as in Claim \ref{claim:witness}, we have a $\Sigma^-_{\P_0}$-hod pair $(\P_0^+,\Pi)$\footnote{$\Sigma^-_{\P_0}$ is just $\Sigma$ since $\P_0 = \mathcal{H}^+$.} such that
\begin{enumerate}[(a)]
\item $\lambda^{\P_0^+}$ is limit ordinal of the form $\alpha' + \omega$, and  such that $\Lambda_{Y} \leq_w \Pi_{\P_0^+(\alpha')}$ (so $\Lambda_{X_i}\leq_w \Pi_{\P_0^+(\alpha')}$ for all $i$).
\item $(\P_0^+,\Pi\rest V)\in V$.
\item In $M$, $\P_0^+$ is countable and $\Gamma(\P_0^+(\alpha'),\Pi_{\P_0^+(\alpha')}) \vDash \mathcal{A}$ is a bad tuple.
\item $\Pi$ has branch condensation, strong hull condensation, is $j(\Gamma)$-fullness preserving. 
\item $\Pi = j(\Pi\rest V)^j$.
\end{enumerate}
The properties above for $(\P_0^+,\Pi)$ can be obtained by a proof similar to that of Claim \ref{claim:witness}, with the last clause coming from Lemma \ref{lem:pullback}.

This type of reflection is possible because we replace $j(\mathcal{H}^+)$ by $\M^Y_{\infty}$. If $\mathcal{Z}$ is the result of iterating $\P_0^+$ via $\Pi$ above $\delta_{\alpha'}^{\P_0^+}$ to make $\mathbb{R}^{M}$ generic (see Footnote \ref{ftn:gen_iter}), then letting $h$ be $\mathcal{Z}$-generic for the Levy collapse of the supremum of $\mathcal{Z}$'s Woodin cardinals such that $\mathbb{R}^{M}$ is the symmetric reals of $\mathcal{Z}[h]$, then in $\mathcal{Z}(\mathbb{R}^{M})$,

\begin{center}
$\Gamma(\P_0^+(\alpha'),\Pi_{\P_0^+(\alpha')}) \vDash \mathcal{A}$ is a bad tuple.
\end{center}

Now we define by induction $\xi_i^+: \P_i^+ \rightarrow \Q_i^+$, $\pi_i^+: \Q_i^+ \rightarrow \P_{i+1}^+$, $\phi_{i,i+1}^+: \P_i^+\rightarrow \P_{i+1}^+$ as follows. $\phi_{0,1}^+: \P_0^+\rightarrow \P_{1}^+$ is the ultrapower map by the extender of length $\Theta_{X_1}$ derived from $\pi_{X_0,X_1}$. Note that $\phi_{0,1}^+$ extends $\phi_{0,1}$. Let $\xi_0^+: \P_0^+ \rightarrow \Q_0^+$ extend $\xi_0$ be the ultrapower map by the $(\textrm{crt}(\xi_0),\delta^{\Q_0})$-extender derived from $\xi_0$. Finally let $\pi_0^+ = (\phi^+_{0,1})^{-1}\circ \xi_0^+$. The maps $\xi_i^+, \pi_i^+, \phi_{i,i+1}^+$ are defined similarly. Let also $\M_Y = \textrm{Ult}(\P_0^+,F)$, where $F$ is the extender of length $\Theta_Y$ derived from $\pi_{X,Y}$. There are maps $\epsilon_{2i}: \P_i^+ \rightarrow \M_Y$, $\epsilon_{2i+1}:\Q_i^+\rightarrow \M_Y$ for all $i$ such that $\epsilon_{2i} = \epsilon_{2i+1}\circ \xi^+_i$, $\epsilon_{2i} = \epsilon_{2i+2}\circ \phi^+_{i,i+1}$, and $\epsilon_{2i+1} = \epsilon_{2i+2}\circ \pi_i^+$. Let $\pi: \M_Y\rightarrow j(\P_0^+)$\footnote{$\pi = \sigma_1\circ \sigma_0$, where $\sigma_0: \M_Y\rightarrow \pi_E(\P_0^+)$ is given by $\sigma_0(\pi_{X,Y}(f)(a)) = \pi_E(f)(\pi_Y(a))$ for $f\in \P_0^+$ and $a\in [\Theta_Y]^{<\omega}$ and $\sigma_1: \pi_E(\P_0^+)\rightarrow j(\P_0^+)$ is defined as: $\sigma_1(\pi_E(f)(a)) = j(f)(a)$ for $f\in \P_0^+$ and $a\in [\pi_E(\Theta)]^{<\omega}$.} be the factor map. When $i = 0$, $\epsilon_0$ is simply $\pi_F$, the ultrapower map by $F$. That these maps are well-defined and the objects $\P_i^+, \Q_i^+$ end-extend $\P_i,\Q_i$ respectively come from the fact that $j\rest \mathcal{H}^+$ is continuous (see Lemma \ref{fullness}).

Letting $\Sigma_i = \Sigma^-_{\P_i}$ and $\Psi_i=\Sigma_{\Q_i}^-$, there is a finite sequence of ordinals $t$ and a formula $\theta(u,v)$ such that in $\Gamma(\P_0^+,\Pi)$
\begin{enumerate}
\setcounter{enumi}{7}
\item for every $i<\omega$, $(\phi,s)\in T_{\P_i,A_i} \Leftrightarrow \theta[i^{\Sigma_i}_{\P_i(\alpha),\infty}(s),t]$, where $\alpha$ is least such that $s\in [\delta_\alpha^{\P_i}]^{<\omega}$;
\item for every $i$, there is $(\phi_i,s_i)\in T_{\Q_i,\xi_i(A_i)}$ such that $\neg \theta[i^{\Psi_i}_{\Q_i(\alpha),\infty}(s_i),t]$ where $\alpha$ is least such that $s_i\in [\delta_\alpha^{\Q_i}]^{<\omega}$.
\end{enumerate}
The pair $(\theta,t)$ essentially defines a Wadge-initial segment of $\Gamma(\P_0^+,\Pi)$ that can define the pair $(\M^Y_{\infty}, A^*)$, where $\tau_i(A_i)=A^*$ for some (any) $i$. In fact, these parameters are inside $\Gamma(\P_0^+(\alpha),\Pi)$.

Let $\Pi_i$ be the $\pi\circ \epsilon_i$-pullback of $j(\Pi)$. Hence, 
\begin{center}
$\Sigma_Y \leq_w \Pi_0=\Pi = j(\Pi\rest V)^j \leq_w \Pi_1 \dots \leq_w j(\Pi\rest V)^{\pi}$. 
\end{center}

%

We can use the strategies $\Pi_i$'s to simultaneously execute a $\mathbb{R}^{M}$-genericity iterations. We outline the process here. First we rename $\langle\P^+_i,\Q^+_i,\xi^+_i,\phi^+_{i,i+1},\pi^+_i$$\ | \ i<\omega\rangle$ to $\langle\P_i^0,\Q_i^0,\xi_i^0,\phi_i^0,\pi_i^0 \ | \ i<\omega\rangle$. We fix in $M^{Col(\omega,\mathbb{R})}$, $\langle x_i \ | \ i<\omega\rangle$, a generic enumeration of $\mathbb{R}^M$. We get $\langle \mathcal{P}^n_i,\mathcal{Q}^n_i, \xi^n_i, \phi^n_i, \pi^n_i, \tau^n_i,k^n_i \ | \ n \leq \omega \wedge i<\omega\rangle$ such that
\begin{enumerate}[(i)]
\item $\mathcal{P}^\omega_i$ is the direct limit of the $\P^n_i$'s under maps $\tau^n_i$'s for all $i < \omega$.
\item $\mathcal{Q}^\omega_i$ is the direct limit of the $\mathcal{Q}^n_i$'s under maps $k^n_i$'s for all $i<\omega$.
\item $\P^n_\omega$ is the direct limit of the $\P^n_i$'s under maps $\pi^n_i$'s.
\item for all $n\leq \omega$, $i<\omega$, $\phi^n_i: \mathcal{P}^n_i \rightarrow \mathcal{P}^n_{i+1}$; $\xi^n_i: \mathcal{P}^n_i \rightarrow \mathcal{Q}^n_i$; $\pi^n_i: \mathcal{Q}^n_i \rightarrow \mathcal{P}^n_{i+1}$ and $\phi^n_i = \pi^n_i\circ\xi^n_i$.
\end{enumerate}
Then we start by iterating $\mathcal{P}^0_0$ above $\delta_\alpha^{\P^0_0}$ to $\P^1_0$ to make $x_0$-generic at $\delta_{\alpha+1}^{\P^1_0}$; say the tree is $\T_0$. We let $\tau^0_0:\P^0_0\rightarrow \P^1_0$ be the iteration map. During this process, we lift $\T_0$ to all $\P^0_n, \Q^0_n$ for $n<\omega$ using the maps $\xi^0_i, \phi^0_i$. We pick branches for the trees on $\P^0_i, \Q^0_i$ according to the strategies $\Pi_i$. We describe this process for the models $\Q^0_0, \P^0_1$. Let $\W$ be the end model of the lift-up tree $\xi^0_0 \T$ on $\Q^0_0$. Note that the tree $\xi^0_0\T_0$ is according to $\Pi_1$. We then iterate $\W$ to $\Q^1_0$ (using $(\Pi_1)_\W$) to make $x_0$ generic at $\delta_{\alpha+1}^{\Q^1_0}$. Let $\xi^1_0: \P^1_0\rightarrow \Q^1_0$ be the natural embedding. Let $\T_1$ be the $x_0$-genericity iteration tree on $\W$ just described and $\W^*$ be the last model of $\phi^0_0\T_0 ^\smallfrown \xi \T_1$, where $\xi$ is the natural map from $\W$ to the last model of $\phi^0_0 \T_0$.  We then iterate the end model of the lifted stack $\phi^0_0\T_0 ^\smallfrown \xi \T_1$ on $\Q^0_1$, noting that this stack is according to $\Pi_2$, to $\Q^1_1$ to make $x_0$ generic at $\delta_{\alpha+1}^{\Q^1_1}$. Let $k^0_0:\Q^0_0 \rightarrow \Q^1_0$, $\tau^0_1: \P^0_1\rightarrow \P^1_1$ be the iteration embeddings, $\pi^1_0: \Q^1_0 \rightarrow \P^1_1$ be the natural map, and $\phi^1_0 = \pi^1_0\circ\xi^1_0$. Continue this process of making $x_0$ generic for the later models $\Q^0_n$'s and $\P^0_n$'s for $n < \omega$. We then start at $\P^1_0$ and repeat the above process, iterating above $\delta_{\alpha+1}^{\P^1_0}$ to make $x_1$ generic at images of $\delta_{\alpha+2}^{\P^1_0}$ etc. This whole process defines models and maps $\langle \mathcal{P}^n_i,\mathcal{Q}^n_i, \xi^n_i, \phi^n_i, \pi^n_i, \tau^n_i,k^n_i \ | \ n \leq \omega \wedge i<\omega\rangle$ as described above.

The process yields a sequence of models $\langle\P^+_{i,\omega}=\P^\omega_i,\Q_{i,\omega}^+ = \Q^\omega_i \ | \ i<\omega\rangle$ and maps $\xi^+_{i,\omega}=\xi^\omega_i:\P^+_{i,\omega}\rightarrow \Q^+_{i,\omega}$, $\pi^+_{i,\omega}=\pi^\omega_i:\Q^+_{i,\omega}\rightarrow \P^+_{i+1,\omega}$, and $\phi^+_{i,i+1,\omega}=\phi^\omega_{i} = \pi^+_{i,\omega}\circ \pi^+_{i,\omega}$. Furthermore, each $\P^+_{i,\omega}, \Q^+_{i,\omega}$ embeds into a $j(\Pi\rest V)^{\pi}$-iterate of $\M_Y$ and hence the direct limit $\P_\infty$ of $(\P^+_{i,\omega},\Q^+_{j,\omega} \ | \ i,j<\omega)$ under maps $\pi^+_{i,\omega}$'s and $\xi^+_{i,\omega}$'s is wellfounded. See Figure \ref{3Ddiagram}.

\begin{figure}
\centering
\begin{tikzpicture}[description/.style={fill=white,inner sep=1pt}]
\matrix (m) [matrix of math nodes, row sep=1.5em,
column sep=1.5em, text height=1ex, text depth=0.15ex]
{ \P^0_0 & & \P^0_1 & & \P^0_2 & & \ .\ \ . \ \ .  \\
  & \Q^0_0 & & \Q^0_1 & & & & \\
\P^1_0 & & \P^1_1 & & \P^1_2 & & \ .\ \ . \ \ .  \\
  & \Q^1_0 & & \Q^1_1 & & & & \\
\P^2_0 & & \P^2_1 & & \P^2_2 & & \ .\ \ . \ \ . \\
  & \Q^2_0 & & \Q^2_1 & & & & \\
. &   & . &   & . &   & & . \\
. & . & . & . & . &   & & . \\
. & . & . & . & . &   & & . \\
  & . &   & . &   &   & &   \\
\P^\omega_0 & & \P^\omega_1 & & \P^\omega_2 & & \ .\ \ . \ \ . & \P_\infty\\
  & \Q^\omega_0 & & \Q^\omega_1 & & & & \\
};

\path[->,font=\scriptsize]
(m-1-1) edge node[description] {$ \phi^0_0 $} (m-1-3)
edge node[description] {$ \xi^0_0 $} (m-2-2)
(m-2-2) edge node[description] {$ \pi^0_0 $} (m-1-3)
(m-1-3) edge node[description] {$ \phi^0_1 $} (m-1-5)
edge node[description] {$ \xi^0_1 $} (m-2-4)
(m-2-4) edge node[description] {$ \pi^0_1$} (m-1-5)
(m-1-5) edge node[description] {$ \phi^0_2$} (m-1-7)
(m-3-1) edge node[left][description] {$ \xi^1_0 $} (m-3-3)
edge node[description] {$ \pi^1_0 $} (m-4-2)
(m-4-2) edge node[description] {$ j^1_0 $} (m-3-3)
(m-3-3) edge node[left][description] {$ \phi^1_1 $} (m-3-5)
edge node[description] {$ \xi^1_1 $} (m-4-4)
(m-4-4) edge node[description] {$ \pi^1_1$} (m-3-5)
(m-3-5) edge node[description] {$ \phi^1_2$} (m-3-7)
(m-1-1)edge node[description]{$ \tau^0_0 $}(m-3-1)
(m-2-2)edge node[above right]{$ k^0_0 $}(m-4-2)
(m-1-3)edge node[description]{$ \tau^0_1$}(m-3-3)
(m-2-4)edge node[above right]{$ k^0_1 $}(m-4-4)
(m-1-5)edge node[description]{$ \tau^0_2$}(m-3-5)
(m-5-1) edge node[left][description] {$ \phi^2_0 $} (m-5-3)
edge node[description] {$ \xi^2_0 $} (m-6-2)
(m-6-2) edge node[description] {$ \pi^2_0 $} (m-5-3)
(m-5-3) edge node[left][description] {$ \phi^2_1 $} (m-5-5)
edge node[description] {$ \xi^2_1 $} (m-6-4)
(m-6-4) edge node[description] {$ \pi^2_1$} (m-5-5)
(m-5-5) edge node[description] {$ \phi^2_2$} (m-5-7)
(m-3-1)edge node[description]{$ \tau^1_0 $}(m-5-1)
(m-4-2)edge node[above right]{$ k^1_0 $}(m-6-2)
(m-3-3)edge node[description]{$ \tau^1_1$}(m-5-3)
(m-4-4)edge node[above right]{$ k^1_1 $}(m-6-4)
(m-3-5)edge node[description]{$ \tau^1_2$}(m-5-5)
(m-5-1)edge node[description]{$ \tau^2_0$}(m-7-1)
(m-6-2)edge node[description]{$ k^2_0$}(m-8-2)
(m-5-3)edge node[description]{$ \tau^2_1$}(m-7-3)
(m-6-4)edge node[description]{$ k^2_1$}(m-8-4)
(m-5-5)edge node[description]{$ \tau^2_2$}(m-7-5)
(m-11-1) edge node[description] {$ \phi^\omega_0 $} (m-11-3)
edge node[description] {$ \xi^\omega_0 $} (m-12-2)
(m-12-2) edge node[description] {$ \pi^\omega_0 $} (m-11-3)
(m-11-3) edge node[description] {$ \phi^\omega_1 $} (m-11-5)
edge node[description] {$ \xi^\omega_1 $} (m-12-4)
(m-12-4) edge node[description] {$ \pi^\omega_1$} (m-11-5)
(m-11-5) edge node[description] {$ \phi^\omega_2 $}(m-11-7)
;
\end{tikzpicture}
\caption{The $(x_n: n<\omega)$ genericity iteration process}
\label{3Ddiagram}
\end{figure}
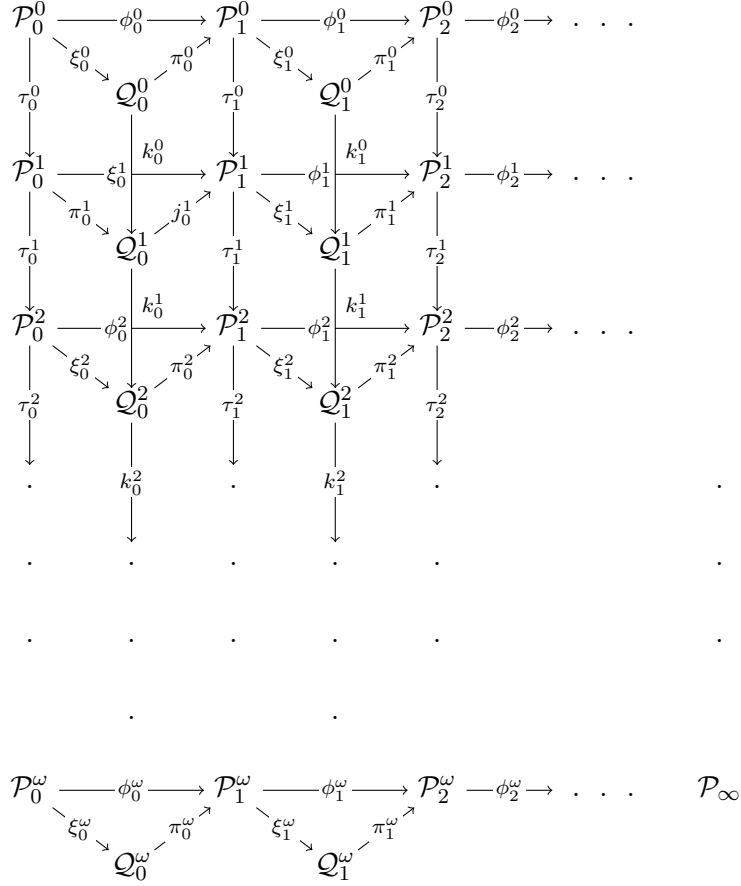


Let $C_i$ be the derived model of $\P^+_{i,\omega}$, $D_i$ be the derived model of $\Q^+_{i,\omega}$ (at the sup of the Woodin cardinals of each model), then $\mathbb{R}^{N} = \mathbb{R}^{C_i} = \mathbb{R}^{D_i}$. Furthermore, $C_i\cap \powerset(\mathbb{R})\subseteq D_i\cap \powerset(\mathbb{R})\subseteq C_{i+1}\cap \powerset(\mathbb{R})$ for all $i$.

(8), (9) and the construction above give us that there is a $t\in [\textrm{OR}]^{<\omega}$, a formula $\theta(u,v)$ such that
\begin{enumerate}
\setcounter{enumi}{9}
\item for each $i$, in $C_i$, for every $(\phi,s)$ such that $s\in \delta^{\P_i}$, $(\phi,s)\in T_{\P_i,A_i}\Leftrightarrow \theta[i^{\Sigma_i}_{\P_i(\alpha),\infty}(s),t]$ where $\alpha$ is least such that $s\in [\delta_\alpha^{\P_i}]^{<\omega}$.
\end{enumerate}
Let $n$ be such that for all $i\geq n$, $\xi^+_{i,\omega}(t) = t$. Such an $n$ exists because the direct limit $\P_\infty$ is wellfounded.\footnote{We can arrange that $\P_\infty$ embeds into a $j(\Pi)^+$-iterate of $j(\P_0^+)$, where $j(\Pi)^+$ is the canonical extension of $j(\Pi)$ in $N$.} By elementarity of $\xi^+_{i,\omega}$ and the fact that $\xi^+_{i,\omega}\rest \P_i = \xi_i$,
\begin{enumerate}
\setcounter{enumi}{10}
\item for all $i\geq n$, in $D_i$, for every $(\phi,s)$ such that $s\in \delta^{\Q_i}$, $(\phi,s)\in T_{\Q_i,\xi_i(A_i)}\Leftrightarrow \theta[i^{\Psi_i}_{\Q_i(\alpha),\infty}(s),t]$ where $\alpha$ is least such that $s\in [\delta_\alpha^{\Q_i}]^{<\omega}$.
\end{enumerate}
However, using (9), we get
\begin{enumerate}
\setcounter{enumi}{11}
\item for every $i$, in $D_i$, there is a formula $\phi_i$ and some $s_i\in [\delta^{\Q_i}]^{<\omega}$ such that $(\phi_i,s_i)\in T_{\Q_i,\xi_i(A_i)}$ but $\neg \phi[i^{\Psi_i}_{\Q_i(\alpha),\infty}(s_i),t]$ where $\alpha$ is least such that $s\in [\delta_\alpha^{\Q_i}]^{<\omega}$.
\end{enumerate}
Clearly (11) and (12) give us a contradiction. This shows that $\pi_X$ has $A$-condensation. Since $A$ is arbitrary, $\pi_X$ has condensation. This completes the proof of the theorem.
\end{proof}

From the above theorems, we obtain a nice, $j$-realizable iteration strategy  $\Lambda$ with the following property:
\begin{adjustwidth}{1cm}{1cm}
letting $\M_\infty(\mathcal{H}^+,\Lambda)$ be the direct limit of (all countable) $\Lambda$-iterates of $\mathcal{H}^+$ in $M$, then there is an elementary map $\tau:\M_{\infty}(\mathcal{H}^+,\Lambda)\rightarrow j(\mathcal{H}^+)$ such that $\tau \circ \pi^{\Lambda}_{\mathcal{H}^+,\infty} = j\rest \mathcal{H}^+$ and if $\tau$ is nontrivial, then $\textrm{crt}(\tau) = \delta^{\M_\infty(\P^+,\Lambda)}$. 
\end{adjustwidth}

The map $\tau$ is defined as follows: for any $x\in \M_\infty(\mathcal{H}^+,\Lambda)$, let $\R\in M$ be a $\Sigma$-iterate of $\mathcal{H}^+$ such that there is some $y\in \R$ such that $\pi^{\Lambda_\R}_{\R,\infty}(y)=x$. Now by construction of $\Lambda$, there is a map $\tau_\R:\R\rightarrow j(\mathcal{H}^+)$ such that $j\rest \P^+ = \tau_\R\circ \pi^\Lambda_{\mathcal{H}^+,\R}$ and $\tau_\R\rest \delta^\R$ agrees with the iteration map by $\Lambda$. We then let $\tau(x)=\tau_\R(y)$. $\tau$ is well-defined by the fact that some iterate of $\Lambda$ has branch condensation and is commuting.

The reason $\Lambda$ is nice is because by construction, whenever $i: \mathcal{H}^+\rightarrow \R$ is according to $\Lambda$, letting $\tau_\R:\R\rightarrow j(\mathcal{H}^+)$ be given by the construction of $\Lambda$, then $\tau_\R\rest \delta^\R = \pi^{\Lambda_\R}_{\R,\infty}\rest \delta^\R$ and $\Lambda_{\R|\delta^\R} = j(\Sigma)^\tau$. From this and standard theorems in the theory of hod mice, see \cite[Theorem 3.26]{ATHM}, we get that for all $\alpha < \delta^\R$, $\Lambda_{\R(\alpha)}$ satisfies (i) Definition \ref{NiceStrategy}. The other two clauses are also clear. Furthermore, if $\tau:\M_{\infty}(\mathcal{H}^+,\Lambda)\rightarrow j(\mathcal{H}^+)$ is as above and is nontrivial, then since $\tau$ is the ``direct limit" of the $\tau_\R$'s for non-dropping $\Lambda$-iterates $\R$ of $\mathcal{H}^+$, crt$(\tau) = \delta^{\M_\infty(\P^+,\Lambda)}$.

There are two cases. The first case is when $\tau$ is non-trivial, we then have that $\delta^{\M_\infty(\P^+,\Lambda)}$ is a regular cardinal which is a limit of Woodin cardinals of $\M_\infty(\P^+,\Lambda)$. Furthermore, by fullness preservation of $\Lambda$, $\delta^{\M_\infty(\P^+,\Lambda)} = \theta_\alpha^{j(\Gamma)}$ for some $\alpha$ and hence $L(\M_\infty(\P^+,\Lambda), j(\Gamma)|\theta_\alpha)\models ``\sf{AD}^+ + $$\Theta$ is regular." Contradiction to our smallness assumption.

The remaining case is when $\tau$ is trivial. In other words, $\Gamma(\mathcal{H}^+,\Lambda) = j(\Gamma)$. By elementarity, there is a reasonable pair $(\P,\Lambda)$ in $V$ such that $\Lambda$ is fullness preserving, has hull and branch condensation, is pullback consistent, commuting, and $\Gamma(\P,\Lambda)=\Gamma$. We need to show.
\begin{lemma}\label{lem:limitUB}
There is a reasonable hod pair $(\Q,\Psi)$ such that $\Q\in V$ is countable, $\Gamma = \Gamma(\Q,\Psi)$, $\Psi$ has a unique extension $\Psi^+$ that acts on stacks in $H_{\omega_2}^V$ and $\Psi$ is $\omega_1$-UB.
\end{lemma} 
\begin{proof}
Let $(\P,\Lambda)$ be a reasonable pair in $V$ such that $\Lambda$ is fullness preserving, has hull and branch condensation, is pullback consistent, commuting, and $\Gamma(\P,\Lambda)=\Gamma$.  Let $(\R,\Lambda')$ be the result of boolean comparing all ``finite variations" of $\Lambda$ i.e. for a $g\subset Coll(\omega,\omega_1)$, for a condition $q$, let $g_q = g - g\rest dom(q) \cup q$, let $\Lambda_q = j_{g_q}(\Lambda)$ and compare in $V[g]$ all pairs $(\P,\Lambda_q)$ (see Section \ref{sec:booleanComp}). 

\begin{claim}\label{claim:terminate}
The boolean comparisons outlined above succeeds and hence $(\R,\Lambda')$ above exists.
\end{claim}
\begin{proof}

The argument is basically from \cite[Theorem 2.47]{ATHM}. Suppose the comparison doesn't succeed. We can then build a ``diamond sequence" of length $\omega_1$. More precisely, we have a sequence $\mathcal{B} = (\R_\alpha,\S_\alpha,\P_\alpha^q, \vec{\T}_\alpha, \vec{\U}_\alpha, \vec{\W}_\alpha^q, b^q_\alpha, i_{\alpha}^q, j^q_\alpha, \xi_\alpha : \alpha < \beta < \omega_1 \wedge q\in Coll(\omega,\omega_1^V))$, where
\begin{enumerate}
\item $\R_0 = \P$. 
\item $\vec{\T}_0^\smallfrown \vec{\U}_0$ is a minimal disagreement between $\{(\R_0, \Lambda_q) : q\in Coll(\omega,\omega_1^V) \}$. $\vec{\T}_0$ is according to all $\Lambda_q$ with last model $\S_0$. For each $q$, $b^q_0 = (\Lambda_q)_{\vec{\T}_0}(\vec{\U}_0)$ and $i_{0}^q: \S_0\rightarrow \P^q_\alpha$ is the iteration embedding according to $(\Lambda_q)_{\vec{\T}_0,\S_0}$, i.e. $i^q_0 = i^{\vec{\U_0}}_{b^q}$. We write $\Psi_{0,q}$ for $\Lambda_q$.

\item For $\alpha>0$, $\vec{\T}_\alpha^\smallfrown \vec{\U}_\alpha$ is a minimal disagreement between $\{(\R_\alpha, \Psi_{\alpha,q}) : q\in Coll(\omega,\omega_1^V) \}$. $\vec{\T}_\alpha$ is according to all $\Psi_{\alpha,q}$ with last model $\S_\alpha$, where $\Psi_{\alpha,q}$ is the appropriate tail of $\Lambda_q$ on $\R_\alpha$ via the stack $\oplus_{\beta<\alpha}\vec{\T}_\beta^\smallfrown \vec{\U}_\beta^\smallfrown \vec{\W}^q_\beta$. For each $q$, $b^q_\alpha = (\Psi_{\alpha,q})_{\vec{\T}_\alpha,\S_\alpha}(\vec{\U}_\alpha)$ and $i_{\alpha}^q: \S_\alpha \rightarrow \P^q_\alpha$ is the corresponding iteration embedding according to $(\Psi_{\alpha,q})_{\vec{\T}_\alpha}$, i.e. $i^q_\alpha = i^{\vec{\U_\alpha}}_{b^q_\alpha}$.
\item $j^q_\alpha: \P^q_\alpha \rightarrow \R_{\alpha+1}$ are iteration maps via stack $\vec{\W}^q_\alpha$ according to $(\Psi_{\alpha,q})_{\vec{\T}_\alpha^\smallfrown \vec{\U}_\alpha^\smallfrown b^q_\alpha}$.
\item For any $\alpha$, for any $\beta < \lambda^{\R_{\alpha+1}}$, for any $p\neq q$, 
\begin{center}
$(\Psi_{\alpha,p})_{\vec{\T}_\alpha^\smallfrown \vec{\U}_\alpha^\smallfrown {b^p_\alpha}^\smallfrown \vec{\W}^p_\alpha, \R_{\alpha+1}(\beta)} = (\Psi_{\alpha,q})_{\vec{\T}_\alpha^\smallfrown \vec{\U}_\alpha^\smallfrown {b^q_\alpha}^\smallfrown \vec{\W}^q_\alpha,\R_{\alpha+1}(\beta)}$.
\end{center}
but for some $p\neq q$,
\begin{center}
$(\Psi_{\alpha,p})_{\vec{\T}_\alpha^\smallfrown \vec{\U}_\alpha^\smallfrown {b^p_\alpha}^\smallfrown \vec{\W}^p_\alpha, \R_{\alpha+1}} \neq (\Psi_{\alpha,q})_{\vec{\T}_\alpha^\smallfrown \vec{\U}_\alpha^\smallfrown {b^q_\alpha}^\smallfrown \vec{\W}^q_\alpha,\R_{\alpha+1}}$.
\end{center}
\item For each $\beta<\omega_1$, $\xi_\beta$ is the least $\xi\in (\delta(\vec{\T}_\beta),\lambda^{\S_\beta})$\footnote{$\delta(\vec{\T}_\beta)$ is the supremum of generators used along $\vec{\T}_\beta$.} such that $\vec{\U}_\beta$ is a stack on $\S_\beta(\xi+1)$ and there are $p\neq q$ such that $(\Psi_{\alpha,p})_{\vec{\T}_\beta,\S_\beta(\xi+1)}\neq(\Psi_{\alpha,q})_{\vec{\T}_\beta,\S_\beta(\xi+1)} $ but for all $p, q$ $(\Psi_{\alpha,p})_{\vec{\T}_\beta,\S_\beta(\xi)} = (\Psi_{\alpha,q})_{\vec{\T}_\beta,\S_\beta(\xi)} $; so $b^p_\alpha\neq b^q_\alpha$.
\end{enumerate}

Clause (6) explains the term ``minimal disagreement" used in (2) and (3). By our assumption, for each $\alpha$, there are $p\neq q$ such that $b^p_\alpha\neq b^q_\alpha$, equivalently $\vec{\U}_\alpha$ witnesses $(\Psi_q)_{\vec{\T}_\alpha,\S_\alpha}\neq (\Psi_p)_{\vec{\T}_\alpha,\S_\alpha}$. For each $\alpha, q$, let $\Sigma_q$ be the appropriate tail of $\Lambda_q$ on $\P^q_\alpha$ and $\lambda^{\alpha,q}$ be the order type of the Woodin cardinals of $\P^q_\alpha$. The maps $j^q_\alpha$ (in (4)) exist by the process of simultaneously comparing all $(\P^q_\alpha, (\oplus_{\alpha<\lambda^{\alpha,q}} \Sigma_q(\alpha))$ into a common hod pair construction inside $j(\Gamma)$.  Furthermore, the common model of the comparison exists and is called $\R_{\alpha+1}$, see \cite[Theorem 2.47]{ATHM}.  The main point is $\{q: q\in Coll(\omega,\omega_1^V)\}$ is countable in $M$ and the supremum of the Wadge ranks of $\{\oplus_{\alpha<\lambda^{\alpha,q}} \Sigma_q(\alpha) : q\in Coll(\omega,\omega_1^V)\}$ is bounded in $j(\Gamma)$. That is why we can find a coarse $\Omega$-Woodin mouse $(N, \Psi, \delta^N)$ that Suslin captures $\{\oplus_{\alpha<\lambda^{\alpha,q}} \Sigma_q(\alpha) : q\in Coll(\omega,\omega_1^V)\}$ (and a universal $\Omega$-set, for $\Omega\subsetneq j(\Gamma)$, a Suslin co-Suslin pointclass containing all $\{\oplus_{\alpha<\lambda^{\alpha,q}} \Sigma_q(\alpha) : q\in Coll(\omega,\omega_1^V)\}$) and performs the above comparison with the hod pair construction done inside $N$ to guarantee (5); this process is further explained in Section \ref{sec:booleanComp}. The comparison succeeds for each $\alpha < \omega_1$. So the sequence is of length $\omega_1$.

Now, the proof of \cite[Theorem 2.49]{ATHM} gives us a contradiction. We sketch the proof here for the reader's convenience. Let $\mathcal{B}$ be the sequence above and let $X_0\prec X_1 \prec H_{\omega_2}$  be countable and contain all relevant objects (recall we work in $V[g]$). Let $\pi_i: H_i \rightarrow X_i$ be the uncollapse map, $\kappa_i = crt(\pi_i)$ for $i\in \{0,1\}$ and let $\pi: H_0\rightarrow H_1$ be the map $\pi_1^{-1}\circ \pi_0$. For each $p\in Coll(\omega,\omega_1^V)$, let $j^p_{\kappa_0,\kappa_1}$ be the iteration embedding from $\R_{\kappa_0}$ to $\R_{\kappa_1}$ by $\Psi{\kappa_0,p}$. It is easy to see that (see  \cite[Theorem 2.49]{ATHM} for the simple calculations) for each such $p$:
\begin{center}
$j^p_{\kappa_0,\kappa_1} = \pi\rest \R_{\kappa_0}$.
\end{center}
Let then $j^p: \S_{\kappa_0}\rightarrow \R_{\kappa_1}$ be the embeddings according to $(\Psi_{\kappa_0,p})_{\vec{\T}_{\kappa_0},\S_{\kappa_0}}$. For each $x\in \S_{\kappa_0}$, let $f\in \R_{\kappa_0}$ and $a\in \delta(\vec{\T}_{\kappa_0})^{<\omega}$ such that $x = \pi^{\vec{\T}_{\kappa_0}}(f)(a)$, it is easy to see that
\begin{center}
$j^p(x) = \pi(f)(j^p(a))$.
\end{center}
But note that the maps $j^p\rest \delta(\vec{\T}_{\kappa_0})$ agree (by property (6)), so indeed, the maps $j^p$ agree on $\S_{\kappa_0}$. Using this and pullback consistency, an argument just as in \cite[Theorem 2.48]{ATHM} shows that for all $p, q$,
\begin{equation}\label{eqn:equal_tails}
(\Psi_{\kappa_0,p})_{\vec{\T}_{\kappa_0},\S_{\kappa_0}(\xi_{\kappa_0}+1)}(\vec{\U}_{\kappa_0})=(\Psi_{\kappa_0,q})_{\vec{\T}_{\kappa_0},\S_{\kappa_0}(\xi_{\kappa_0}+1)}(\vec{\U}_{\kappa_0}).
\end{equation}
This clearly contradicts (6).

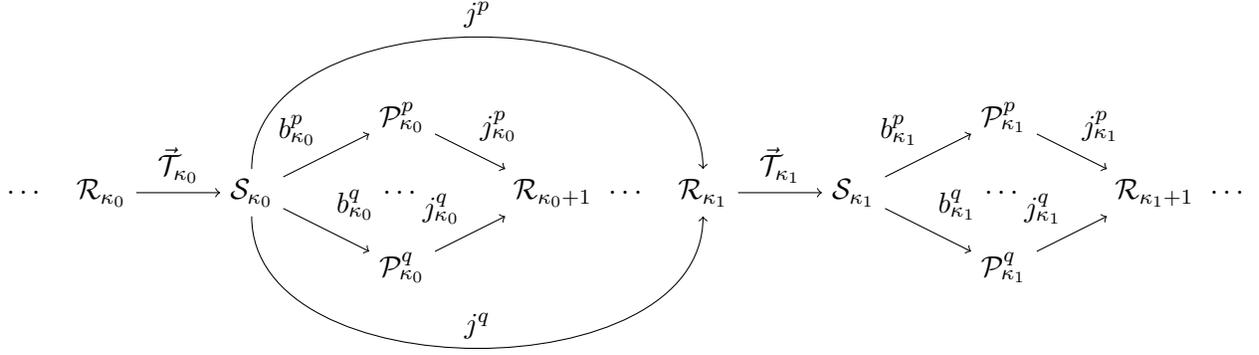
\begin{figure}
\centering
\begin{tikzpicture}[node distance=2cm, auto]
  \node (A) {$\dots$};
  \node (B) [node distance=1cm, right of=A] {$\R_{\kappa_0}$};
  \node (C) [right of=B] {$\S_{\kappa_0}$};
  \node (D) [right of=C] {$\dots$};
  \node (E) [node distance=1cm, above of=D] {$\P^p_{\kappa_0}$};
  \node (F) [node distance=1cm, below of=D] {$\P^q_{\kappa_0}$};
  \node (G) [right of=D] {$\R_{\kappa_0+1}$};
  \node (H) [node distance=1cm, right of=G] {$\dots$};
  \node (I) [node distance=1cm, right of=H]{$\R_{\kappa_1}$};
  \node (J) [right of=I] {$\S_{\kappa_1}$};
  \node (K) [right of=J] {$\dots$};
  \node (L) [node distance=1cm, above of=K] {$\P^p_{\kappa_1}$};
  \node (M) [node distance=1cm, below of=K] {$\P^q_{\kappa_1}$};
  \node (N) [right of=K] {$\R_{\kappa_1+1}$};
  \node (P) [node distance=1cm, right of=N] {$\dots$};

 \draw[->] (B) to node {$\vec{\T}_{\kappa_0}$}(C);
 \draw[->] (C) to node {$b^p_{\kappa_0}$}(E);
  \draw[->] (C) to node {$b^q_{\kappa_0}$}(F);
   \draw[->] (E) to node {$j^p_{\kappa_0}$}(G);
   \draw[->] (F) to node {$j^q_{\kappa_0}$}(G);   
   
    \draw[->] (I) to node {$\vec{\T}_{\kappa_1}$}(J);
 \draw[->] (J) to node {$b^p_{\kappa_1}$}(L);
  \draw[->] (J) to node {$b^q_{\kappa_1}$}(M);
   \draw[->] (L) to node {$j^p_{\kappa_1}$}(N);
   \draw[->] (M) to node {$j^q_{\kappa_1}$}(N);   
  \draw[->,bend left=90] (C) to node {$j^p$}(I);
  \draw[->,bend right=90] (C) to node {$j^q$}(I); 
    \end{tikzpicture}
\caption{A Diamond sequence}
\label{fig:diamond}
\end{figure}

The equality (\ref{eqn:equal_tails}) above holds because, by pullback consistency, for each $p$
\begin{center}
$(\Psi_{\kappa_0,p})_{\vec{\T}_{\kappa_0},\S_{\kappa_0}(\xi_{\kappa_0}+1)}(\vec{\U}_{\kappa_0}) = (\Psi_{\kappa_1,p})_{\R_{\kappa_1}(j^p(\xi_{\kappa_0}+1))}(j^p \vec{\U}_{\kappa_0})$,
\end{center}
and by (5) and the fact that the maps $j^p$'s agree on $\S_{\kappa_0}$, for any $p,q$,
\begin{center}
$ (\Psi_{\kappa_1,p})_{\R_{\kappa_1}(j^p(\xi_{\kappa_0}+1))}(j^p \vec{\U}_{\kappa_0}) = (\Psi_{\kappa_1,p})_{\R_{\kappa_1}(j^q(\xi_{\kappa_0}+1))}(j^q \vec{\U}_{\kappa_0})$.
\end{center}
This easily implies (\ref{eqn:equal_tails}).
\end{proof}

By the claim, $\R\in V$ and $\Lambda' \rest V\in V$.  By Lemma \ref{lem:pullback}, $\Lambda' = j(\Lambda')^j$. By elementarity, there is a hod pair $(\Q,\Psi)$ such that $\Q\in V$ is countable, an elementary embedding $\pi: \Q\rightarrow \R$ such that $\Psi = (\Lambda')^\pi$ and $\Gamma(\Q,\Psi) = \Gamma$. $\Psi$ is an $(\omega_2,\omega_2)$-strategy with branch condensation is $\Gamma$-fullness preserving.

\begin{claim}\label{claim:UB}
$\Psi$ is $\omega_1$-UB.
\end{claim}
\begin{proof}
Let $\M_\infty$ be the direct limit of all non-dropping iterates of $\Psi$ and $i: \Q\rightarrow \M_\infty$ be the direct limit map. Let $a = (\Q,i,\M_\infty, \Psi)$. We want to find a formula $\varphi[x,a]$ such that for a club of countable $X\prec H_{\omega_3}$ containing all relevant objects, letting $\pi_X:M_X\rightarrow X$ be the uncollapse and $(\omega_1^X, a^X) = \pi_X^{-1}(\omega_1, a)$, for any $M_X$-generic $h\subset Coll(\omega,\omega_1^X)$ in $V$, for any $\vec{\T}\in M_X[h]$ according to $\Psi$ and $b\in M_X[h]$ a cofinal branch of $\vec{\T}$,
\begin{equation}\label{eqn:UB_equiv}
M_X[h]\models \varphi[(\vec{\T},b), a^X] \Leftrightarrow V\models \varphi[(\vec{\T},b), a].
\end{equation}

We describe the formula $\varphi$. $\varphi[x, a]$ says:
\begin{itemize}
\item If $x_0^\smallfrown x_1$ is non-dropping, letting $\vec{\T} = x_0$ and $b = x_1$, then there is a map $\sigma: \M^{\vec{\T}}_b\rightarrow a_2$ such that $\sigma \circ i^{\vec{\T}}_b = a_1$. Here we think of $x$ as an ordered pair $(x_0,x_1)$ and $a$ as an ordered tuple $(a_0,a_1,a_2,a_3)$.
\item If $x_0^\smallfrown x_1$ drops, letting $\vec{\T} = x_0$, $b = x_1$, and $(\vec{\T}_\alpha, \P_\alpha, \xi_\alpha: \alpha\leq \nu \wedge \xi_\alpha < \lambda^{\P_\alpha})$ be the essential components of $\vec{\T}$,\footnote{See \cite[Definition 2.37]{ATHM}.} then for $\alpha < \nu$, $\vec{\T}_\alpha$ does not drop and is based on $\P_\alpha(\xi_\alpha)$, $\vec{\T}_\nu$ is a stack based on $\P_\nu(\xi_\nu)$, $\xi_\nu$ is a successor ordinal,  with cofinal branch $b$. There is a $\Q'$, a nondropping iterate of $a_0$ according to $a_3$, with iteration map $\tau:a_0\rightarrow \Q'$, and a $\sigma_1: \P_\nu\rightarrow \Q'$ such that $\sigma_1\circ i^{\vec{\T}_{<\nu}} = \tau$, where $\vec{\T}_{<\nu}=\oplus_{\alpha<\nu} \vec{\T}_\alpha$, and letting $\Psi' = (a_3)_{\Q}^{\sigma_1}$\footnote{We suppress from the notation the stack $\vec{\U}$ with iteration map $\tau$ and last model $\Q$, technically we should write $(a_3)^\sigma_{\vec{\U},\Q}$. This case includes the case $\nu=0$ and $\tau,\sigma_1$ are the identity maps.}, then $b = \Psi'(\vec{\T}_\nu)$. 
\end{itemize}

Now we show \ref{eqn:UB_equiv}, suppose $\vec{\T}^\smallfrown b$ does not drop, then the equivalence follows easily. This is because if $M_X[h]\models \varphi[(\vec{\T},b), a^X]$, then $\pi_X\circ \sigma: \M^{\vec{\T}}_b\rightarrow \M_\infty$ is such that $\pi_X\circ \sigma\circ i^{\vec{\T}}_b = \pi_X\circ a_1$. $\pi_X\circ a_1 = \pi_X(a_1): \Q\rightarrow \M_\infty$ is the direct limit map. By branch condensation, $b=\Psi(\vec{\T})$. If $V\models \varphi[(\vec{\T},b), a]$, then again by branch condensation, $b = \Psi(\vec{\T})$. $M_X[h]\models \varphi[(\vec{\T},b), a^X]$ by boolean comparisons done inside $M_X[h]$.

Suppose $b$ drops. Then clearly, $\Q(b,\vec{\T}_\nu)$ exists. If  $M_X[h]\models \varphi[(\vec{\T},b), a^X]$ then clearly $V\models \varphi[(\vec{\T},b), a]$. Conversely, by boolean comparison, we can find a $\tau,\Q'\in M_X$ and $\sigma_1\in M_X[h]$ that satisfy the second clause above. Letting $\Psi' = (a_3)_{\Q}^{\sigma_1}$ and $c = \Psi'(\vec{\T}_\nu)$, then since $\Psi'$ is fullness preserving (see \cite[Theorem 3.26]{ATHM}), $\Q(\vec{\T},c)$ must exist, and therefore $c = b$ as shown in the previous section. In both cases, $b=\Psi(\vec{\T})$.

\end{proof}
\end{proof}
We can then proceed with the CMI and show Lp$^{\Lambda^+}(\mathbb{R})\models \sf{AD}^+$ and go on with the induction.

\section{OUTLINE OF THE PROOF OF THEOREM \ref{cor:equiconsistency}}\label{sec:lower_bound_2}

We outline the argument constructing models of ``$\sf{AD}_\mathbb{R} + $$\Theta$ is regular" from the assumption that the non-stationary ideal on $\powerset_{\omega_1}(\mathbb{R})$ is strong and pseudo-homogeneous. We let $\mathcal{I}$ be the non-stationary ideal on $\powerset_{\omega_1}(\mathbb{R})$. Let $G\subseteq \mathbb{P}_{\mathcal{I}}$ be $V$-generic and $j=j_G: V \rightarrow M = \textrm{Ult}(V,G)\subseteq V[G]$ be the generic embedding. Let $k: M\rightarrow N$ be the generic embedding given by an $M$-generic $H\subset j(\mathbb{P}_\mathcal{I})$. We note that
\begin{itemize}
\item $j(\omega_1) = \mathfrak{c}^+$ (by the strength of the ideal).
\item The properties in Lemma \ref{lem:quasi_hom} hold for $j$.
\item Letting $M = \textrm{Ult}(V,G)$. $M$ need not be closed under $\omega$-sequences in $V[G]$. In particular, $\mathbb{R}^M$ may differ from $\mathbb{R}^{V[G]}$. Also, $\mathfrak{c}^+$ may be $> \omega_2^V$.
\end{itemize}

We let $\Gamma$ be defined as in Section \ref{sec:outline} and operate under the smallness assumption $(\ddag)$ as before. Our inductive hypothesis in this case is: 
\begin{center}
(**): \ \ \ \ if $J$ is a $\Sigma$-cmi operator for some reasonable hod pair $(\P,\Sigma)$ such that $\Sigma$ is definable in $V$ from a countable sequence of ordinals, then $J$ is definable in $V$ from a countable sequence of ordinals.
\end{center}
The core model induction is very similar to the one given in the previous section; however, instead of maintaining the inductive hypothesis $(\dag)$, we maintain $(**)$. We mention some key points below. The details are left to the reader. We fix the pair $(\P,\Sigma)$ as in $(**)$. $(\P,\Sigma)$ is allowed to be $(\emptyset,\emptyset)$.
\begin{itemize}
\item If $J$ is a $\Sigma$-cmi operator on (a cone above some $a$ in) $H^V_{\omega_1}$ that satisfies $(**)$, then by pseudo-homogeneity, we can show $j(J)\rest V\in V$ and by strongness, $j(J)\rest V\in V$ has domain the cone above $a$ in $H_{\mathfrak{c}^+}^V$. The definability calculations are done in $M$ and $V[G]$ plays no role in the argument. For instance, one can show using pseudo-homogeneity (as the base case) that $j(\Sigma)\rest V\in V$.
\item One can then show the existence of $\M_1^{\F,\sharp}$ whenever $\F$ is a $\Sigma$-cmi operator that satisfies $(**)$. Using pseudo-homogeneity again, one shows the operator $H: x \mapsto \M_1^{\F,\sharp}(x)$ has the property that $j(H)\rest V\in V$ and $(**)$ holds for $H$. This is the analog of Theorem \ref{thm:m1sharp}.

\item Theorem \ref{thm:sjs} can be proved by a similar argument, though much simpler as Claim \ref{claim:independent} follows easily from pseudo-homogeneity. The proof of Lemmata \ref{lem:II_wins} and \ref{lem:scale} is also given in \cite{wilson2012contributions}. This gives also that $o(\textrm{Lp}^\Sigma(\mathbb{R})) < j(\omega_1) =\mathfrak{c}^+$.

\item The above gives an analog of Theorem \ref{thm:omega1UB}, namely the existence of a hod pair $(\P',\Sigma')$ such that $\Sigma'$ is Lp$^\Sigma(\mathbb{R})$-fullness preserving, $\Sigma'\notin \textrm{Lp}^\Sigma(\mathbb{R})$, and $\Sigma'$ is definable in $V$ from a countable sequence of ordinals.

\item In the limit case, we can define in $M$ the model $\mathcal{H}^+$ (see \ref{eqn:cases}) from $j\rest\mathcal{H}$. Since $j\rest \mathcal{H}$ is independent of $G$ and hence $j\rest \mathcal{H}\in V$, $\mathcal{H}^+\in V$ by pseudo-homogeneity. 
\item By an argument similar to that of Proposition \ref{cor:H_small}, $\mathcal{H}^+$ is countable in $M$. We can argue $j$ is continuous at $o(\mathcal{H}^+)$ as follows.
\begin{claim}\label{claim:jcont}
Let $\gamma = o(\mathcal{H}^+)$. Then $j(\gamma) = \textrm{sup}_{\alpha < \gamma} j(\alpha)$.
\end{claim}
\begin{proof}
We first claim $j\rest \mathcal{H}^+\in V$. Let $\prec$ be the canonical well-order of $\mathcal{H}^+$; $\prec$ is definable over $\mathcal{H}^+$. We think of $\prec$ as a bijection from $o(\mathcal{H}^+)$ onto $\mathcal{H}^+$. Note that $j(\mathcal{H}^+)\in V$ (equivalently $j(\prec)\in V$) and $j\rest o(\mathcal{H}^+)\in V$ (this follows from the above discussion). $j\rest\mathcal{H}^+$ can be easily computed from $j\rest o(\mathcal{H}^+), j(\mathcal{H}^+), j(\prec)$. Therefore, $j\rest \mathcal{H}^+\in V$.

Suppose for contradiction that $j(\gamma) >  \textrm{sup}_{\alpha < \gamma} j(\alpha)$. Let $\nu =  \textrm{sup} \ j[\gamma]$. Let $\vec{C} = (C_\alpha : \alpha < \gamma)$ be the canonical $\square_\Theta$-sequence defined over $\mathcal{H}^+$ (see \cite{schimmerling2004characterization} for a construction of such a sequence). Let $D = j(\vec{C})_\nu$. Since $\nu < j(\gamma)$, $D$ is defined and is club in $\nu$. Furthermore, since $j\rest\mathcal{H}^+\in V$, 
\begin{center}
cof$^V(\nu) = \textrm{cof}^V(\gamma) > \omega$.
\end{center}
Since $j(\mathcal{H}^+)\in V$, cof$^{j(\mathcal{H}^+)}(\nu) > \omega$. This, in particular, implies that the set of limit points of $D$ is non-empty and in fact a club in $\nu$. By the property of $\square$-sequences, for each limit point $\alpha\in D$, 
\begin{center}
$D\cap \alpha = j(\vec{C})_\alpha$.
\end{center}

Since $j\rest \mathcal{H}^+\in V$, $E =_{def} j^{-1}[D]\in V$ is an $\omega$-club in $\nu$ with the property: for all limit point $\alpha$ of $E$ with cof$^V(\alpha)=\omega$,
\begin{center}
$E\cap \alpha = C_\alpha$.
\end{center}
By the construction of $\vec{C}$, $E$ induces a $\P\lhd Lp^{\Sigma,j(\Gamma)}(\mathcal{H})$, but also that every $\M \lhd \mathcal{H}^+ = Lp^{\Sigma,j(\Gamma)}(\mathcal{H})$ is an initial segment of $\P$. So $\P \notin Lp^{\Sigma,j(\Gamma)}(\mathcal{H})$. Contradiction.

\end{proof}

\item We can show the corresponding claim in Section \ref{sec:lim} that continuity of $j$ at $\lambda^{\mathcal{H}}$ implies cof$^V(\lambda^\mathcal{H}) = \omega$ as follows. If $\kappa\in [\omega_1,\mathfrak{c}]$ is a successor cardinal or a weakly inaccessible  cardinal, then $j$ is discontinuous at $\kappa$. This is because $j\rest \kappa \in M$ and if $j$ is continuous at $\kappa$, then $j(\kappa)$ is singular in $M$. This contradicts the fact that $j(\kappa)$ is successor or weakly inaccessible, hence regular, in $M$. This implies cof$^V(\lambda^{\mathcal{H}}) = \omega$. The proof that $|\mathcal{H}^+|\leq \mathfrak{c}$, $\Sigma\rest V\in V$ and does not depend on $G$, $\mathcal{H}^+ = \textrm{Lp}^{\Sigma,j(\Gamma)}(\mathcal{H})\models ``\textrm{cof}^V(\lambda^{\mathcal{H}})$ is measurable" (if $j$ is discontinuous at $\lambda^{\mathcal{H}}$) is similar, using pseudo-homogeneity.

\item From this point on, we assume $j$ is discontinuous at $\lambda^{\mathcal{H}}$ and hence $\mathcal{H}^+ = \textrm{Lp}^{\Sigma,j(\Gamma)}(\mathcal{H})\models ``\textrm{cof}^V(\lambda^{\mathcal{H}})$ is measurable". Otherwise, the argument is much easier.

\item Claim \ref{claim:jcont} and the above argument show cof$(o(\mathcal{H}^+)) = \omega$. 

\item By Lemma \ref{lem:quasi_hom}, arguments in Proposition \ref{cor:H_small} and the fact that $j$ is continuous at $o(\mathcal{H}^+)$, we get that $j\rest \mathcal{H}^+\in V\cap M$. 

\item The analog of Lemma \ref{fullness} is the following.

\begin{lemma}\label{lem:fullness2}
$\Lambda$ is $j(\Gamma)$-fullness preserving .
\end{lemma}
\begin{proof}
Suppose not. Let $\vec{\mathcal{T}}$ be according to $\Lambda$ with end model $\Q$ such that $\Q$ is not $j(\Gamma)$-full. This means there is a strong cut point $\gamma$ such that letting $\alpha\leq \lambda^{\Q}$ be the largest such that $\delta^{\Q}_\alpha \leq \gamma$, then without loss of generality, in $j(\Gamma)$, there is a mouse $\M\lhd \textrm{Lp}^{\Sigma_{\Q(\alpha)}}(\Q|\gamma)$\footnote{The case where $\gamma = \delta_\alpha$ and $\M\lhd \textrm{Lp}^{\oplus_{\beta<\alpha}\Sigma_{\Q(\beta)}}(\Q|\gamma)$ is similar.} such that $\M\notin \Q$. Let $l: \Q\rightarrow j(\mathcal{H}^+)$ be such that $j\rest \mathcal{H}^+ = l \circ i^{\vec{\mathcal{T}}}$; here by the above discussions, $j(\mathcal{H}^+)=\pi_E(\mathcal{H}^+ = \{j(f)(a): a\in [j(\Theta)]^{<\omega} \wedge f\in \mathcal{H}^+ \}$ and $l$ is defined as: 
\begin{center}
$l(i^{\vec{\T}}(f)(a)) = j(f)(i^{\Sigma_\Q}_{\Q,\infty}(a))$,
\end{center}
where $f\in \mathcal{H}^+$, $a\in [\delta^\Q]^{<\omega}$. Here $E$ is the (long) extender of length $o(j(\mathcal{H)})$ derived from $j$. We use $i$ to denote $i^{\vec{\mathcal{T}}}$ from now on. 

\begin{claim}\label{claim:witness}
There is a $\Sigma$-hod pair $(\P,\Phi)$ such that 
\begin{enumerate}[(a)]
\item $\P\in V$, $\Phi\rest V\in V$,\footnote{By $\Phi\rest V$, we mean $\Phi\rest H^V_{\mathfrak{c}^+}$.} and $\Phi\in j(\Gamma)$ is fullness preserving and has branch condensation.
\item $\P$ is countable in $M$, $\lambda^{\P}$ is limit and cof$^{\P}(\lambda^{\P})$ is not measurable in $\P$.
\item in $j(\Gamma)$, $\Gamma(\P,\Phi)$ witnesses $\Lambda$ is not fullness preserving. 
\end{enumerate}
\end{claim}
\begin{proof}

First note that in $M$, there is some $\alpha$ such that $\Sigma_\M$, the canonical strategy of $\M$, is in $j(\Gamma)|\delta^{\P^*}_\alpha$, where $\P^* = \textrm{HOD}^{j(\Gamma)}_{\Sigma}(\alpha)$ \footnote{We identify $\textrm{HOD}^{j(\Gamma)}_{\Sigma}$ with the direct limit of $\Sigma$-hod pairs $(\R,\Psi)$ and $\Psi$ is fullness preserving and has branch condensation in $j(\Gamma)$.}  and $\P^*\models \exists \beta \alpha = \beta+\omega$.  Such $\P^*$ and $\alpha$ exists by our assumptions on $\Gamma$. $\P^*\in V$ follows from pseudo-homogeneity. Let $\Psi$ be the strategy of $\P^*$ which is the tail of some (equivalently, all) $\Sigma$-hod pair $(\R,\Psi^*)\in j(\Gamma)$ $\Psi$ is fullness preserving and has branch condensation in $j(\Gamma)$ and $\M_\infty(\R,\Psi^*) = \P^*$. Note that $\Psi$ is fullness preserving and has branch condensation in $k(j(\Gamma))$. It follows that $\Psi\rest V \in V$. From pseudo-homogeneity, we can ordinal define $\Psi\rest V$ in $M$ from $\Sigma$ and $\P$ with the prescription above, using the fact that $j(\Gamma)$ is $OD$ in $M$ and $j(\Theta)$, the Wadge rank of $j(\Gamma)$, doesn't depend on the choice of $G$.  

We also have that $j(\P^*)\in V$. This is because $j(\P^*)$ is definable in $M$ from $\{j(\alpha), j(\mathcal{I}), j(\mathcal{H}^+)\}$, but $j(\mathcal{I})$ and $j(\mathcal{H}^+)$ are both definable in $M$.\footnote{This is one place where we use the ideal $\mathcal{I}$ is the non-stationary ideal, or just that it is definable in $V$. Technically, $j(\mathcal{H}^+)$ is definable in $M$ from $j(\mathcal{H})$ and a countable sequence of ordinals, namely any sequence $(j(\gamma_n) : n<\omega)$, where $(\gamma_n : n<\omega)$ is cofinal in $o(\mathcal{H}^+)$ and there is $j(\mathcal{H}) \lhd \M_n \lhd j(\mathcal{H}^+)$ such that $o(\M_n)=j(\gamma_n)$.} By an argument similar to that of Claim \ref{claim:jcont}, $j\rest \P^*\in V$. We want to find a countable-in-$M$ version of $\P^*$ in $V$.

Let $(\dot{T}, \dot{\Q}, \dot{\M},\dot{\Lambda})$ be $\mathbb{P}_{\mathcal{I}}$-names for $(\vec{\T},\Q,\M,\Lambda)$ and let $p\in\mathbb{P}_{\mathcal{I}}$ force all relevant facts about these objects. Let $X\prec (H_{\lambda},\in)$ where 
\begin{itemize}
\item $\lambda > \mathfrak{c}^+$ is regular, 
\item $X^\omega\subset X$, 
\item $\mathfrak{c}\cup \Gamma\cup\mathcal{H}^+ \cup\{\dot{T}, \dot{\Q}, \dot{\M}, \Gamma, (\P^*,\Psi\rest V), (j(\P^*),j\rest\P^*) \}\subset X$, and \item $|X| \leq \mathfrak{c}$.
\end{itemize}
Let $\pi: M_X\rightarrow X$ be the transitive uncollapse map and for any $x\in X$, let $\bar{x} =\pi^{-1}(x)$. Note that 
\begin{center}
$\overline{\mathcal{H}^+} = \mathcal{H}^+$. 
\end{center}
Let $\mathbb{P} = \mathbb{P}_\mathcal{I}$ and $h\subset \bar{\mathbb{P}}$ be $M_X$-generic such that $h\in M$. Such an $h$ exists by the properties of $X$. \footnote{We do not have a way of lifting $\pi$ to all of $M_X[h]$. This creates complications and forces us to argue as below.}  

Work in $M_X[h]$, let $(\overline{\T}, \overline{\Q}, \overline{\M},\overline{\Lambda})$ be the interpretation of $(\overline{\dot{\T}}, \overline{\dot{\Q}}, \overline{\dot{\M}}, \overline{\dot{\Lambda}})$. Let $\sigma = j\rest \P^*$; so $\overline{\sigma}: \overline{\P^*}\rightarrow \overline{j(\P^*)}$. Let $\overline{R}$ be the image of $\overline{\P^*}$ under the extender $F$ derived from $i^{\bar{\T}}$, i.e. 
\begin{center}
$\overline{R} = \{i^{\bar{\T}}(f)(a) : f\in \overline{\P^*} \wedge a\in [\delta^{\overline{\Q}}]^{<\omega}\}$. 
\end{center}
Let $i_F: \overline{\P^*}\rightarrow \overline{R}$ be the associated ultrapower map, and let $\bar{l}: \overline{\R}\rightarrow \overline{j(\P^*)}$. Let $\tau: \overline{\R}\rightarrow j(\P^*)$ be $\tau  = \pi \circ \bar{l}$. Note that $\sigma\circ \pi = \tau\circ i_F$.

Let $\Upsilon = j(\Psi\rest V)$ and $\Psi^* = \pi^{-1}(\Psi\rest V)$. In $M_X[h]$, $\overline{\Lambda}$ is not full as witnessed by $\overline{T},\overline{Q},\overline{\M}$ inside $\bar{j}(\overline{\Gamma})|\bar{\alpha}$, where $\bar{j}$ is the generic ultrapower induced by $h$. Therefore, letting $j(\Psi\rest V)^{\tau\circ i_F} = \Sigma_1$ and $j(\Psi\rest V)^\tau = \Sigma_2$, we note that 
\begin{center}
$\Sigma_1 \leq_w \Sigma_2$. 
\end{center}
In $M$,
\begin{center}
$\Gamma(\overline{\P^*},\Sigma_1)\subset \Gamma(\overline{\R},\Sigma_2)$,
\end{center}
and letting $\Sigma_3 = j(\Sigma)^\tau$,
\begin{center}
$L(\Gamma(\overline{\P^*},\Sigma_1)) \models ``\overline{\M}$ is a $\Sigma_3$-mouse and $\neg (\overline{\M}\lhd \overline{\Q})$."
\end{center}
Finally, note that $\overline{\T}$ is according to $\Lambda$ as $\overline{\T}$ is $j$-realizable. It is easy then to see that (a),(b), (c) hold for $(\overline{\P^*},\Sigma_1)$. Therefore, the pair $(\overline{\P^*},\Sigma_1)$ is the desired $(\P,\Phi)$. See Figure \ref{fig:witness} for an illustration of the argument above.

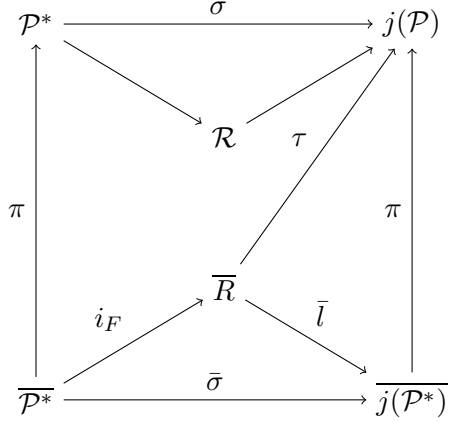
\begin{figure}
\centering
\begin{tikzpicture}[node distance=2.5cm, auto]
  \node (A) {$\P^*$};
  \node (B) [node distance=5cm, below of=A] {$\overline{\P^*}$};
  \node (C) [node distance = 5cm, right of=A] {$j(\P)$};
  \node (D) [node distance=5cm, below of=C] {$\overline{j(\P^*)}$};
  \node (E) [node distance=1.5cm, below of=A] {};
  \node (F) [node distance=2.5cm, right of=E]{$\R$};
  \node (G) [node distance=1.5cm, above of=B]{};
  \node (H) [node distance=2.5cm, right of=G]{$\overline{R}$};
   \draw[->] (A) to node {$\sigma$}(C); 
  \draw[->] (B) to node {$\bar{\sigma}$}(D);
  \draw[->] (B) to node {$\pi$}(A);
  \draw[->] (D) to node  {$\pi$} (C);
  \draw[->] (A) to node {$$}(F);
  \draw[->] (F) to node {$$}(C);
  \draw[->] (B) to node {$i_F$}(H);
  \draw[->] (H) to node {$\bar{l}$}(D);
  \draw[->] (H) to node {$\tau$}(C);
  \end{tikzpicture}
 
\caption{Diagram for the proof of Claim \ref{claim:witness}. Here $\tau = \pi\circ \overline{\l}$.}
\label{fig:witness}
\end{figure}

\end{proof}

Now we proceed to finish the proof of Lemma \ref{lem:fullness2}. Let $(\P,\Phi)$ be as in the claim. We assume that $L(\Gamma(\P,\Phi))$ satisfies the statement: ``$\Q$ is not full as witnessed by $\M$", i.e. we reuse the notation for $\vec{\T}, \Q, \M, l$. By arguments similar to that used in Lemma \ref{lem:no_proj_across}, no levels of $\P$ projects across $\Theta$ and in fact, $o(\mathcal{H}^+)$ is a cardinal of  $\P$. The second clause follows from the following argument. Suppose not and for simplicity, let $\mathcal{H}^+ \unlhd \N\lhd \P$ be least such that $\rho_1(\N) = \Theta$. Let $f: \kappa^* \rightarrow \Theta$ be an increasing and cofinal map in $\mathcal{H}^+$, where $\kappa^* = \textrm{cof}^{\mathcal{H}^+}(\Theta)$. $\N$ is intercomputable with the sequence $g = \langle \N_\alpha \ | \ \alpha < \kappa^*\rangle$, where $\N_\alpha = Th_{\Sigma_1}^\N( \delta^{\mathcal{H}^+}_{f(\alpha)}\cup \{p_\N\})$. Note that $\N_\alpha\in \mathcal{H}^+$ for each $\alpha<\kappa^*$. Now let $\R_0 = \textrm{Ult}_{0}(\mathcal{H}^+,\mu)$, $\R_1= \textrm{Ult}_{1}(\N,\mu)$, where $\mu\in \mathcal{H}^+$ is the (extender on the sequence of $\mathcal{H}^+$ coding a) measure on $\kappa^*$ with Mitchell order $0$. Let $i_0:\mathcal{H}^+ \rightarrow \R_0$, $i_1:\N\rightarrow \R_1$ be the ultrapower maps. Letting $\delta = \delta_{\lambda^{\mathcal{H}^+}} = \Theta$, it's easy to see that $i_0\rest (\kappa^*+1)=i_1\rest(\kappa^*+1)$ and $\powerset(\delta)^{\R_0} = \powerset(\delta)^{\R_1}$. The second equality follows from the fact that $\R_0$ is full in $j(\Gamma)$ (and hence in $k(j(\Gamma))$).

This means $\langle i_1(\N_\alpha) \ | \ \alpha<\kappa^*\rangle \in \powerset(\delta)^{\R_0}$. By fullness of $\mathcal{H}^+$ in $j(\Gamma)$, $\langle i_1(\N_\alpha) \ | \ \alpha<\kappa^*\rangle \in \mathcal{H}^+$.\footnote{Any $A\subset \delta$ in $\R_0$ is $OD^{j(\Gamma)}_{\Sigma}$, this means OD$^{L(j(\mathbb{R}),C)}_{\Sigma}$ for some $C\in j(\Gamma)$) and so by Strong Mouse Capturing ($\sf{SMC}$, see \cite{ATHM}), $A\in \mathcal{H}^+$.}  Similarly, $\langle i_0(\N_\alpha) \ | \ \alpha<\kappa^*\rangle \in \mathcal{H}^+$. Using these and the fact that $i_0\rest \mathcal{H}^+|\Theta = i_1\rest \N|\Theta \in \mathcal{H}^+$, we can get $\N \in \mathcal{H}^+$ as follows. For any $\alpha < \Theta, \beta<\kappa^*$, $\alpha \in \N_\beta$ if and only if $i_0(\alpha) \in i_1(\N_\beta)= i_0(\N_\beta)$. Since $\mathcal{H}^+$ can compute the right hand side of the equivalence, it can compute the sequence $\langle \N_\alpha \ | \ \alpha<\kappa^*\rangle$. Contradiction.

In other words, $\P$ thinks $\mathcal{H}^+$ is full.  Let $\Psi = \Phi\rest V$ and let
\begin{center}
$i^*: \P\rightarrow \R$
\end{center}
be the ultrapower map by the extender induced by $i$ of length $\delta^\Q$. Note that $\Q\lhd \R$ and $\R$ is wellfounded since there is a natural map \begin{center}
$l^*: \R \rightarrow \P_E$
\end{center} 
extending $l$ and $\pi_E\rest \mathcal{P} = l^*\circ i^*$; here $l^*(i^*(f)(a)) = \pi_E(f)(i^{\Sigma_\Q}_{\Q,\infty}(a))$ for $f\in \P$ and $a\in [\delta^\Q]^{<\omega}$ and $\P_E = \{\pi_E(f)(a) : f\in \P \wedge a\in [j(\Theta)]^{<\omega}\}$. We note here that since $\pi_E$ is continuous at $o(\mathcal{H}^+)$, $j(\mathcal{H}^+)$ is a cardinal initial segment of $\P_E$. Furthermore, there is a natural embedding $\sigma^-: \P_E \rightarrow j(\P)$ such that 
\begin{center}
$j\rest \P = \sigma^-\circ l^* \circ i^*$.
\end{center} 
Here $\sigma^-(\pi_E(f)(a)) = j(f)(a)$ for all $f\in \P$ and $a\in [j(\Theta)]^{<\omega}$. The equality above just comes from the fact that $E$ is an extender derived from $j$.

By the choice of $(\P,\Phi)$, $\M$'s unique strategy $\Sigma_\M \leq_w \Phi$ and $\Sigma_\M\in L(\Gamma(\P,\Phi))$; so in particular, $L(\Gamma(\P,\Phi))$ knows $\Q$ is not full as witnessed by $(\M,\Sigma_\M)$.  

Let $\W = \M_\omega^{\Phi,\sharp}$ and $\Lambda^*$ be the unique strategy of $\W$; again $\W\in V$, $\W$ is countable in $M$, and $\Lambda^*\rest V\in V$. Furthermore, by fullness of $\P$, $o(\P)$ is a cardinal of  $\W$. Let $\W^*$ be a $\Lambda^*$-iterate of $\W$ below its first Woodin cardinal that makes $(\Q,\vec{\T})$ generic via the $(\Q,\vec{\T})$-genericity iteration. Letting $K$ be the generic for the extender algebra of $\W^*$ at its first Woodin cardinal such that $(\Q,\vec{\T})\in \W^*[K]$, then the derived model $D(\W^*[K])$ (at the supremum of the Woodin cardinals of $\W^*[K]$) satisfies
\begin{center}
$L(\Gamma(\P,\Phi),\mathbb{R}) \vDash \Q \textrm{ is not full}.\footnote{Here we abuse notations a bit, by using the same notation for $\Phi$ and its various restrictions.}$\footnote{This is because we can continue iterating $\W^*$ above the first Woodin cardinal to $\W^{**}$ such that letting $\lambda$ be the sup of the Woodin cardinals of $\W^{**}$, then there is a $Col(\omega,<\lambda)$-generic $h$ such that $\mathbb{R}^{V[G]}$ is the symmetric reals for $\W^{**}[h]$. And in $\W^{**}(\mathbb{R}^{V[G]})$, the derived model satisfies that $L(\Gamma(\P,\Phi)) \vDash \Q$ is not full. In the above, we have used the fact that the interpretation of the UB-code of the strategy for $\P$ in $\mathcal{W}^{**}$ to its derived model is $\Phi\rest\mathbb{R}^{V[G]}$; this key fact is proved in \cite[Theorem 3.26]{ATHM}.}
\end{center}
So the above fact is forced over $\W^*[K]$.

Now further extend $i^*$ to $i^+: \W \rightarrow \mathcal{Y}$ and extend $l^*$ to $l^+: \mathcal{Y}\rightarrow \W_E$ so that $\pi_E\rest \mathcal{W} = l^+ \circ i^+$; $i^+, l^+,\W_E$ are defined in a similar manner as above. Again, there is a natural map $\sigma: \W_E\rightarrow j(\W)$ such that $\sigma\circ l^+\circ i^+ = j\rest \W$. Note that $(\mathcal{Y},\sigma\circ l^+)$ are countable in $M$; this is the key reason we need $\P$ is countable in $M$. Therefore, it makes sense to pullback in $M$ via $\sigma\circ l^+$. Let
\begin{center}
$\Psi^* = j(\Lambda^*)^{\sigma\circ l^+}$.
\end{center}
Now note that $\Phi = (\pi_E(\Psi)^{l^*})^{i^*}$ and $\Lambda^* = (\Psi^*)^{i^+}$, so 
\begin{equation}\label{eqn:one}
\Gamma(\mathcal{\P},\Phi)\subseteq \Gamma(\mathcal{R},\pi_E(\Psi)^{l^*})
\end{equation}
and 
\begin{equation}\label{eqn:two}
\Lambda \leq_w \Psi^*. 
\end{equation}
Now iterate $\mathcal{Y}$ using $\Psi^*$ to $\mathcal{Y}^*$ above $\Q$ to make $\mathbb{R}^{M}$ generic \footnote{We write $(\delta_i^{\mathcal{Y}}: i<\omega)$ for the Woodin cardinals of $\mathcal{Y}$ and a similar notation applies to iterates of $\mathcal{Y}$. We work in $M[L]$ where $L\subseteq Coll(\omega,\mathbb{R}^M)$. We have a generic enumeration $(x_n : n<\omega)$ of $\mathbb{R}^M$ and we have a sequence of normal trees and models $(\T_n , \mathcal{Y}_n: n<\omega)$ according to $\Psi^*$, where $\T_0$ is on $\mathcal{Y}=\mathcal{Y}_0$, $\T_n$ is a $x_n$-genericity iteration tree on $\mathcal{Y}_n$ on the window $(\delta^{\mathcal{Y}_n}_{n-1},\delta^{\mathcal{Y}_n}_n)$ according to the $\T_{n-1}$-tail of $\Psi^*$, here $\delta_{-1}^{\mathcal{Y}}=0$. Letting $\mathcal{Y}_\infty$ be the direct limit, then $\mathbb{R}^M$ is the symmetric reals of $\mathcal{Y}_\infty$ for some $g\subseteq Coll(\omega, <\lambda)$, where $\lambda$ is the supremum of the Woodin cardinals of $\mathcal{Y}_\infty$.\label{ftn:gen_iter}}. From \ref{eqn:one} and \ref{eqn:two}, we get that in $D(\mathcal{Y}^*)$,
\begin{center}
$L(\Gamma(\mathcal{R},\pi_E(\Psi)^{l^*}))\vDash \Q$ is not full as witnessed by $\M$.
\end{center}
This gives $\M$ is $OD^{D(\mathcal{Y}^*)}_{\Sigma_\Q}$, so $\M\in \mathcal{Y^*}$ and so $\M\in \R$ since $\R$ is a cardinal initial segment of $\mathcal{Y}^*$. This contradicts the internal fullness of $\Q$ inside $\R$ ($\P$ thinks $\mathcal{H}^+$ is full, so by elementarity, $\R$ thinks $\Q$ is full). See Figure \ref{fig:discontinuous} for an illustration of the argument above.

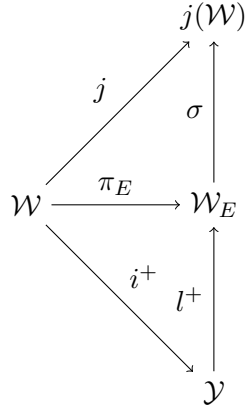
\begin{figure}
\centering
\begin{tikzpicture}[node distance=2.5cm, auto]
  \node (A) {$\W$};
  \node (B) [right of=A] {$\W_E$};
  \node (C) [below of=B] {$\mathcal{Y}$};
  \node (D) [above of=B] {$j(\W)$};
   \draw[->] (A) to node {$i^+$}(C); 
  \draw[->] (C) to node {$l^+$}(B);
  \draw[->] (A) to node {$\pi_E$}(B);
  \draw[->] (A) to node  {$j$} (D);
  \draw[->] (B) to node {$\sigma$} (D);
  \end{tikzpicture}
 
\caption{Diagram for the proof of Lemma \ref{lem:fullness2}.}
\label{fig:discontinuous}
\end{figure}

\end{proof}

\item Regarding the proof of the $j$-condensation lemma (Theorem \ref{thm:cond}), the following are the main changes we need.  Fix a bad tuple $\mathcal{A} = \{\langle \P_i,\Q_i,\tau_i,\xi_i,\pi_i,\sigma_i \ | \ i<\omega \rangle, \M^Y_{\infty}\}$ in $M$ as in the proof of Theorem \ref{thm:cond};  note that $k(\mathcal{A})= \{\langle \P_i,\Q_i,\tau_i,\xi_i,\pi_i,\sigma_i \ | \ i<\omega \rangle, k(\M^Y_{\infty})\}$ is also a bad tuple in $N$ because $k$ fixes all these objects. 

We let $(\P_0^+,\Pi)$ be such that
\begin{enumerate}[(a)]
\item $\P_0^+ = \textrm{HOD}^{j(\Gamma)}_{\Sigma}(\alpha' + \omega)$ for some limit ordinal $\alpha'$ such that $\mathcal{A}\in j(\Gamma) | \theta_{\alpha'}$. Note that $\P_0^+$ is countable in $N$ and $\{k(\P_0^+), k\rest \P_0^+\}\in M$.
\item $\Pi$ is the natural strategy of $\P_0^+$ and is the tail of any $\Sigma$-hod pair $(\R,\Psi)$ such that $\M_\infty(\R,\Psi) = \P_0^+$.
\item $\Pi\rest M\in M$ and $\Pi\rest M \subseteq k(\Pi\rest M)^k$. The latter property follows from the fact that $(\P_0^+,\Pi)$ is a hod pair of limit type, $\Pi$ has branch condensation and is $k(j(\Gamma))$ fullness preserving; therefore, basic theory of hod mice, e.g. the proof of \cite[Theorem 3.26]{ATHM}, implies $\Pi = k(\Pi)^k$. We do not know if the conclusion of Lemma \ref{lem:pullback} holds for all hod pairs constructed in the core model induction here, but fortunately, we do not need it.\footnote{In the context of $\sf{DI}$, we need Lemma \ref{lem:pullback} in situations where the hod pair has successor type. In the case where hod pairs are of limit type, we can argue as above.} We will also write $\Pi$ for $k(\Pi\rest M)^k$ when interpreted in $N$.
\item $\Lambda_{Y} \leq_w \Pi_{\P_0^+(\alpha')}$ (so $\Lambda_{X_i}\leq_w \Pi_{\P_0^+(\alpha')}$ for all $i$) in $N$. Note that we can extend $\Lambda_Y$ (similarly $\Lambda_{X_i}$ for all $i$) in $N$ as the realizable strategy (which we also call $\Lambda_Y$) of $\P_Y$ into $k(j(\mathcal{\mathcal{H}^+}))$ using the map $k\circ \pi_Y$.
\item In $N$, $\P_0^+$ is countable and $\Gamma(\P_0^+(\alpha'),\Pi_{\P_0^+(\alpha')}) \vDash \mathcal{A}$ is a bad tuple.
\end{enumerate}
The rest of the proof is essentially the same as before, but now we run the ``three dimensional argument" using $k$ (instead of $j$) and the argument takes place in $N$ (instead of in $M$). We leave the details to the reader.

\end{itemize}
This completes our outline.

\section{OPEN PROBLEMS AND QUESTIONS}
\label{sec:open}

As mentioned above, there are various important and intriguing questions concerning ideals on $\omega_2$. Woodin has conjectured that (see Theories (a) and (c) in \cite[Question 12]{Woodin})
\begin{conjecture}
The following theories are equiconsistent.
\begin{enumerate}
\item $\sf{ZFC  + MM(\mathfrak{c}) } + $$J_{NS}$ is weakly presaturated.
\item $\sf{ZF + AD_{\mathbb{R}}} + $ ``$\Theta$ is regular". 
\end{enumerate}
\end{conjecture}
In the above $J_{NS}$ is the non-stationary ideal on $\omega_2$ concentrating on ordinals of cofinality $\omega$.  $J_{NS}$ is weakly saturated if for every function $f:\omega_2\rightarrow \omega_2$, for every $S\in \powerset(\omega_2)\slash J_{NS}$, there exists a canonical function $h:\omega_2\rightarrow \omega_2$ such that 
\begin{center}
$\{\alpha\in S : f(\alpha)\leq h(\alpha)\}\notin J_{NS}$
\end{center}
\cite[Theorem 9.137]{Woodin} has established one direction of the conjecture. The converse most likely requires new techniques in the core model induction.

\begin{question}
What is the consistency strength of the theory ``$\sf{ZFC} \ + $ there is a dense ideal on $\omega_2$"?
\end{question}

Finally, as mentioned in the previous section, we do not know the exact consistency strength of the theory ``$\sf{ZFC} \ + $ there is a pseudo-homogeneous ideal on $\powerset_{\omega_1}(\mathbb{R})$", but we conjecture

\begin{conjecture}
The following theories are equiconsistent.
\begin{enumerate}
\item ``$\sf{ZFC} \ + $ there is a pseudo-homogeneous ideal on $\powerset_{\omega_1}(\mathbb{R})$".
\item $\sf{ZF + AD_{\mathbb{R}}} + $ ``$\Theta$ is regular". 
\end{enumerate}
\end{conjecture}  
\bibliographystyle{plain}
\bibliography{Quasi_hom_v2}
\end{document}